\documentclass[review,12pt]{elsarticle}
\usepackage[margin=1in]{geometry}

\usepackage{graphicx} \usepackage{amsmath} \usepackage{amsfonts}
\usepackage{amsmath}
\usepackage{amssymb}
\usepackage{comment}
\usepackage{wrapfig}
\usepackage{floatrow}
\usepackage{breqn}
\usepackage{bm}
\usepackage{multirow}
\usepackage{enumitem}
\usepackage{wasysym}
\usepackage{tikz}
\usepackage{floatrow}
\usepackage{array}
\newcommand*\hexbrace[2]{%
\underset{#2}{\underbrace{\rule{#1}{0pt}}}}
\newcommand{\tikzcircle}[2][red,fill=red]{\tikz[baseline=-0.5ex]\draw[#1,radius=#2](0,0) circle;}

\newcommand\solidrule[1][5mm]{\rule[0.5ex]{#1}{.4pt}}
\newcommand\dashedrule{\mbox{%
  \solidrule[1mm]\hspace{1mm}\solidrule[1mm]}}

\usepackage{hyperref}
\hypersetup{
    colorlinks=false, 
    linktoc=all,     
    linkcolor=black,  
    citecolor=black
}
\def\tsr{\eta}

\def\textw{{\alpha}}
\def\tintw{{\gamma}}

\def\fpt{\emph{findpts}}

\def\dO{\partial \Omega}

\def\dt{ \Delta t }
\def\dx{ \Delta x }

\def\scriptO{{{\it O}\kern -.42em {\it `}\kern + .20em}}
\def\RR{{{\rm l}\kern - .15em {\rm R} }}
\def\PP{{{\rm l}\kern - .15em {\rm P} }}
\def\L2{{{\sf L}^2}}
\def\H1{{{\sf H}^1}}
\def\PN2{{\PP_{N}-\PP_{N-2}}}

\def\complex{{{\rm C} \kern - .53em {\rm l} \kern + .38em}}

\def\a1{{ | \lambda_{\min} |}}

\def\l1{{   \lambda_{\min}  }}

\def\bue{{\underline {\bf e}}}

\def\bu0{{\underline {\bf 0}}}
\def\buu{{\underline {\bf u}}}

\def\bu{{\bf u}}

\def\bx{{\bf x}}

\def\uu{{\underline u}}
\def\uub{{\bar {\underline u}}}

\def\uz{{\underline z}}
\def\u0{{\underline 0}}

\def\ihina{{H_{1}^{-1}}}
\def\ihinb{{H_{2}^{-1}}}

\def\uiiinml{{\uu^{i,n-l}}}
\def\uiiinmlq{{\uu^{i,n-l,Q}}}

\usepackage{lineno}

\journal{Journal}

\begin{document}

\begin{frontmatter}



\title{Stability analysis of a singlerate and multirate predictor-corrector scheme for overlapping grids}

\author{Ketan Mittal\corref{cor1}\fnref{label1}}
\author{Som Dutta\fnref{label2}}
\author{Paul Fischer\fnref{label1,label3}}
\fntext[label1]{Mechanical Science \& Engineering, University of Illinois at Urbana-Champaign, 1206 W. Green St., Urbana, IL 61801}
\fntext[label2]{Mechanical \& Aerospace Engineering, Utah State University, 4130 Old Main Hill, Logan, UT 84332}
\fntext[label3]{Computer Science, University of Illinois at Urbana-Champaign, 201 N. Goodwin Ave., Urbana, IL 61801}
\cortext[cor1]{Corresponding author. Present Address: CASC, Lawrence Livermore National Laboratory, 7000 East Avenue, Livermore, CA 94550}

\begin{abstract}
We use matrix stability analysis for a singlerate and multirate predictor-corrector scheme (PC) used to solve the incompressible Navier-Stokes equations (INSE) in overlapping grids. By simplifying the stability analysis with the unsteady heat equation in 1D, we demonstrate that, as expected, the stability of the PC scheme increases with increase in the resolution and overlap of subdomains. For singlerate timestepping, we also find that the high-order PC scheme is stable when the number of corrector iterations ($Q$) is odd. This difference in the stability of odd- and even-$Q$ is novel and has not been demonstrated in the literature for overlapping grid-based methods. We address the odd-even behavior in the stability of the PC scheme by modifying the last corrector iterate, which leads to a scheme whose stability increases monotonically with $Q$. For multirate timestepping, we observe that the stability of the PC scheme depends on the timestep ratio ($\tsr$).  For $\eta=2$, even-$Q$ is more stable than odd-$Q$. For $\eta\ge3$, even-$Q$ is more stable than odd-$Q$ for a small nondimensional timestep size and the odd-even behavior vanishes as the timestep size increases. The stability analysis presented in this work gives novel insight into a high-order temporal discretization for ODEs and PDEs, and has helped us develop an improved PC scheme for solving the incompressible Navier-Stokes equations.
\end{abstract}

\end{frontmatter}

\section{Introduction}
Numerical solution of partial differential equations (PDEs) is central to much
of today's engineering analysis and scientific inquiry. Techniques such as the
finite element method (FEM), the finite volume method (FVM), and the spectral
element method (SEM) are used to approximate solutions of PDEs on a collection
of volumes or elements whose union constitutes a {\em mesh} that covers the
entire computational domain, $\Omega$.  Construction of an optimal mesh (grid)
is not a trivial task and often becomes a bottleneck for complex domains.
Overlapping Schwarz (OS) based methods circumvent issues posed by conformal
grids by allowing the domain to be represented as the union of simpler
subdomains, each of which can be meshed independently with relatively simple
mesh constructions. The nonconforming union of these meshes allows combinations
of local mesh topologies that are otherwise incompatible, which is a feature of
particular importance for complex 3D domains.  As a results, OS-based methods
have become popular for solving different classes of problems such as
incompressible flow
\cite{mittal2019nonconforming,henshaw1994,merrill2016,rogers1991steady,cd2012v7},
compressible flow
\cite{ahmad1996helicopter,cambier2013onera,eberhardt1985overset,nicholsoverflowmanual,saunier2008third},
electromagnetics \cite{blake1996overset,angel2018}, heat transfer
\cite{meng2017stable,kao1997application,henshaw2009composite}, and particle
tracking \cite{koblitz2017}, and have been implemented using various
discretization approaches such as the finite difference method (FD), FEM, FVM,
and SEM.

The focus of this work is the stability of the temporal discretization used in
the Schwarz-SEM framework for solving the incompressible Navier-Stokes
equations (INSE) on overlapping grids \cite{mittal2019nonconforming,
mittal2020multirate, mittal2020direct, chatterjee2019towards}.  The Schwarz-SEM
framework is based on the spectral element method for monodomain conforming
meshes \cite{patera84, dfm02}, which has demonstrated exponential convergence
of the solution (with the order of the polynomial used for quadrature on each
element) and up to third-order temporal convergence. The Schwarz-SEM framework
uses a high-order spatial interpolation approach \cite{mittal2019nonconforming}
for exchanging overlapping grid solution to maintain the exponential
convergence and a high-order predictor-corrector (PC) timestepping approach to
maintain the temporal convergence of the underlying SEM solver. The spatial and
temporal discretization used in the Schwarz-SEM framework brings forth several
considerations from a stability, accuracy, and computational-cost point of
view. These factors include the extrapolation order used for the interdomain
boundary data at the predictor step, the number of corrector (Schwarz)
iterations at each timestep, and the amount of grid overlap required between
adjacent subdomains.  Here, we present stability analysis to
understand how each of these factors impact the singlerate and multirate
predictor-corrector scheme used in the Schwarz-SEM framework, and develop a
more efficient timestepping approach for the INSE.

Since analysing the stability of the incompressible Navier-Stokes equations in
two- or three-dimensions is not straightforward due to the complexity of the
PDE, we simplify our analysis by considering the unsteady heat equation in 1D
using the finite difference (FD) method.  We specifically choose the unsteady
heat equation due to the similarity in the eigenvalue spectrum of the diffusion
operator and the parabolic unsteady Stokes operator, and use the matrix method
for stability analysis \cite{varga1999matrix} to analyze the
predictor-corrector scheme of the Schwarz-SEM framework. PC schemes have been analyzed
extensively for differential equations with some of the earliest work done
dating back to 1960s \cite{hamming1959stable, chase1962stability,
hall1967stability}, but their understanding in the context of overlapping grids
is limited.  Peet et al. \cite{peet2012} use matrix method for stability
analysis to study a singlerate-based PC scheme for overlapping grids, but their
results fail to describe certain aspects of the stability behavior that we have
observed in the Schwarz-SEM framework.  Mathew et al. \cite{mathew2003maximum}
and Wu et al. \cite{wu2012convergence} use theoretical analysis to understand
stability behavior of their methods for solving PDEs on overlapping grids, but
their method is based on a different temporal discretization compared to the
BDF$k$/EXT$m$-based predictor-corrector scheme that we are interested in.  Meng
et al.'s conjugate heat transfer method \cite{meng2017stable} focuses on
nonoverlapping fluid-solid interfaces and have different boundary conditions
(mixed Dirichlet and Neumann) compared to the Schwarz-SEM framework (purely
Dirichlet).  Love et al. have analyzed a midpoint-based predictor-corrector
scheme for Lagrangian shock hydrodynamics, where they demonstrate a difference
in the stability of odd and even corrector iterates \cite{love2009stability}.
While the temporal discretization used by Love et al. is different from our
method, their result on the difference in stability of odd and even iterate is
similar to what we have observed in the Schwarz-SEM framework.  We have
also found similar evidence of difference in the stability of odd- and even-iterates of
high-order predictor-corrector schemes in the context of ODEs (Fig. 3.2.4 in \cite{dfm02}
and Fig. 1 in \cite{stetter1968improved}). In \cite{stetter1968improved}, Stetter has analyzed an
third-order accurate Adam Bashforth- (AB) and second-order accurate
Adam-Moulton-based (AM) method for ODE of the form
$\frac{\partial y}{\partial t} = \lambda y$ and the stability results therein
indicate that the
high-order PC scheme is relatively more stable for odd-iterates when
$\lambda$ is real-valued. While this observation is not a part of
Stetter's discussion, empirical evidence in our work shows that this behavior
might be applicable to high-order PC schemes in general when $\lambda$ is
real-valued.

The results presented in this paper provide novel insight into the stability
characteristics of a singlerate and multirate predictor-corrector scheme for
overlapping grids. We also develop an improved PC scheme for singlerate
timestepping that addresses the shortcomings in the stability characteristics
of the existing PC scheme.  The remainder of the paper is organized as follows.
Section 2 introduces the matrix method for stability analysis, the overlapping
Schwarz method for solving PDEs, and the FD-based discretization for solving
the unsteady heat equation in a 1D monodomain grid.  In Section 3, we extend
the solution of the unsteady heat equation in a monodomain to overlapping grids
using the singlerate and multirate predictor-corrector scheme from the
Schwarz-SEM framework. In Section 4, we analyze the stability of the singlerate
PC scheme and demonstrate that our simplified analysis qualitatively captures
the stability behavior that we have observed in the Schwarz-SEM framework.
Here, we also developed an improved PC scheme for singlerate timestepping.  In
Section 5, we analyze the multirate PC scheme for an arbitrary timestep ratio.
Finally in Section 6, we summarize our findings and discuss possible directions
for future work.

\section{Preliminaries} \label{sec:prelim}
In this section, we introduce the matrix method for stability analysis,
summarize the overlapping Schwarz method for solving PDEs in
overlapping subdomains, and describe our temporal and spatial discretization
for solving the unsteady heat equation in a single conforming grid.

\subsection{Matrix method for stability analysis}
In the matrix method for stability analysis \cite{varga1999matrix}, if the
system of equations to advance the solution in time can be represented as
\begin{eqnarray}
\label{eq:G}
u^n = G u^{n-1},
\end{eqnarray}
where $u^{n}$ denotes the solution $u$ at time $t^n$,
a sufficient and necessary condition for stability is that the spectral radius
of the propogation operator is
$\rho(G) < 1$.  For our stability analysis, we represent our
predictor-corrector scheme into a system of this form \eqref{eq:G} and use the
spectral radius of $G$ to understand how various
factors of interest impact the stability of the method. We also use this approach to design
a novel singlerate PC scheme that significantly improves the stability
characteristics of the current method (Section \ref{sec:singlestability}).

\subsection{Overlapping Schwarz method}
\begin{figure}[b!] \begin{center}
$\begin{array}{cc}
\includegraphics[height=43mm]{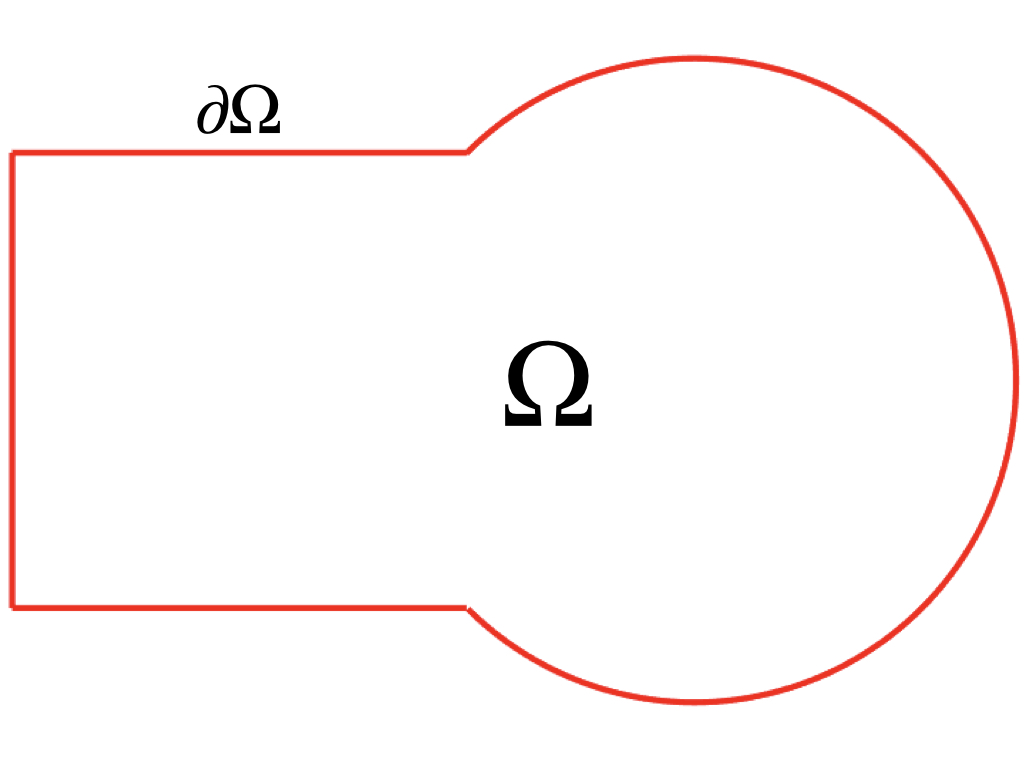} & \hspace{2mm}
\includegraphics[height=45mm]{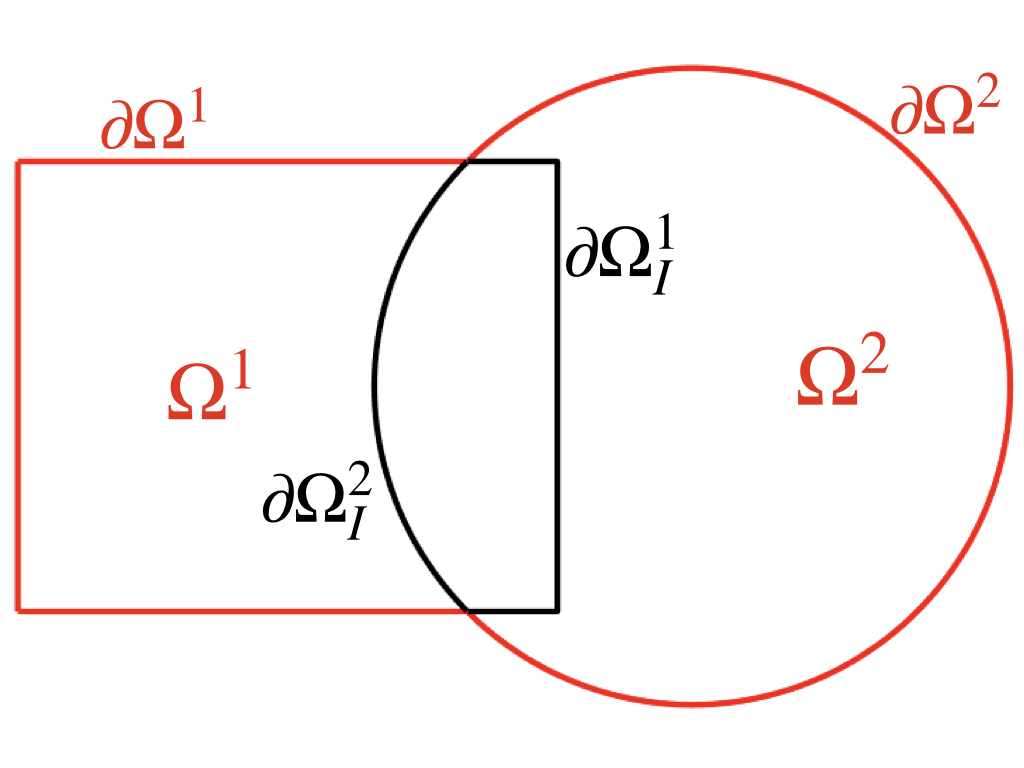}
\end{array}$
\end{center}
\vspace{-7mm}
\caption{(left to right) (a) Composite domain $\Omega$ (b) modeled by overlapping
rectangular ($\Omega^1$) and circular ($\Omega^2$) subdomains. $\dO^{s}_{I}$
denotes the segment of the subdomain boundary $\dO^s$ that is interior to
another subdomain $\Omega^r.$}
\label{fig:schwarz}
\end{figure}

The OS method for solving a PDE in overlapping domains was introduced by
Schwarz in 1870 \cite{schwarz1870}. Figure \ref{fig:schwarz} shows the
composite domain $\Omega$ used in Schwarz's initial model problem, which is
partitioned into two overlapping subdomains: a rectangle ($\Omega^1$) and a
circle ($\Omega^2$). We use $\dO^s_I$ to denote the ``interdomain boundary'',
namely the segment of the subdomain boundary $\dO^s$ that is interior to
another subdomain, and these interdomain boundaries $\dO^1_{I}$ and $\dO^2_I$
are highlighted in Fig. \ref{fig:schwarz}(b).

There are two popular approaches for solving a PDE using the OS method.  In the
Schwarz alternating method, the PDE is solved sequentially in overlapping
subdomains while using the most recent solution to obtain the boundary
condition at interdomain boundaries. As an example, consider the Poisson
equation $-\nabla^2 u=f$ with Dirichlet boundary conditions
($u\rvert_{\dO}=u_b$) on a domain partitioned into $S=2$ overlapping
subdomains. Using $u^{s,q}$ to denote the solution $u$ in $\Omega^s$ at the
$q$th Schwarz iteration, the Schwarz alternating method for solving the Poisson
equation is
\begin{equation}
\begin{aligned}
\label{eq:as1a}
-\nabla^2u^{1,q}&=f \textrm{ in } \Omega^1, \quad  u^{1,q}=u_b \textrm{ on
} \dO^1 \backslash \dO_I^1, \quad u^{1,q}=u^{2,q-1} \textrm{ on } \dO_I^1, \\
-\nabla^2u^{2,q}&=f \textrm{ in } \Omega^2, \quad u^{2,q}= u_b \textrm{
on } \dO^2 \backslash \dO_I^2, \quad u^{2,q}=u^{1,q} \textrm{ on }
\dO_I^2,
\end{aligned}
\end{equation}
with $q=1\dots Q$ for $Q$ Schwarz iterations.  Starting with an initial
condition $u^{s,[0]}$, the Poisson equation is solved sequentially (first in
$\Omega^1$ and then in $\Omega^2$) with interdomain boundary data exchange before
each Schwarz iteration.  The primary drawback of the alternating Schwarz method
is that it does not scale with the number of subdomains ($S$). The sequential
dependencies of the alternating method are overcome by the simultaneous Schwarz
method where the PDE is solved simultaneously in all subdomains with
interdomain boundary data exchange prior to each Schwarz iteration:
\begin{equation}
\begin{aligned}
\label{eq:as2a}
-\nabla^2u^{s,q}&=f \textrm{ in } \Omega^s, \quad  u^{s,q}=u_b \textrm{ on
} \dO^s \backslash \dO_I^s, \quad u^{s,q}=u^{r,q-1} \textrm{ on } \dO_I^s.
\end{aligned}
\end{equation}
Naturally, the advantage of simultaneous Schwarz is its parallelism, which
makes it well suited for large scale problems in an arbitrary number of
overlapping subdomains.  The PC scheme used in the Schwarz-SEM framework is
based on the simultaneous Schwarz method, which we describe in Section
\ref{sec:unsteadyheatms}.

\subsection{Unsteady heat equation in a monodomain grid}
\begin{figure}[b!]
\includegraphics[width=140mm]{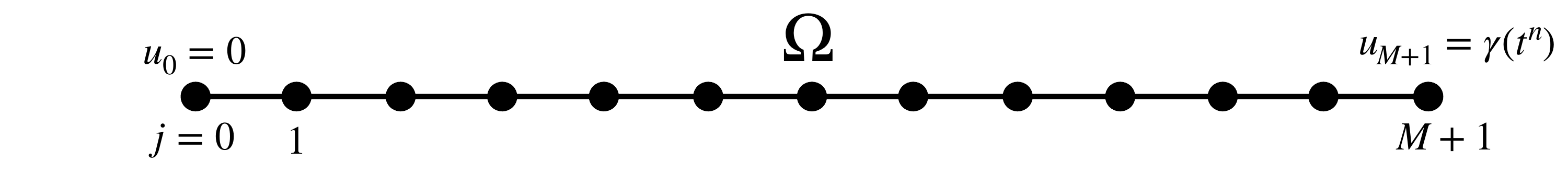}
\vspace{-4mm}
\caption{Monodomain grid with $M=11$.}
\label{fig:monodom} \end{figure}

Since we are interested in the stability analysis of the unsteady diffusion
equation, we first introduce our spatial and temporal discretization in the
context of a monodomain grid. Consider the solution $u(x,t)$ for the unsteady diffusion
equation
\begin{eqnarray}
\label{eq:unsteadyheatmono}
u_t = \nu u_{xx}, \quad x \in [0,1], \quad \nu > 0,
\end{eqnarray}
with homogeneous boundary condition on the left boundary, $u(0,t)=0$ and a
time-dependent inhomogeneous boundary condition on the right boundary,
$u(1,t)=\gamma(t^n)$.

We model the domain with $M+2$ uniformly-spaced grid points (Fig.
\ref{fig:monodom}), such that $\dx = \frac{1}{M+1}$. For notational purposes,
we introduce $u_j^n$ to represent the solution $u$ at $j$th grid point at time
$t^n$, $R$ as the standard restriction matrix (\ref{eq:rmat}) that casts the
solution from the $M+1$ grid points $\uub^n = \lbrack u_1^n,u_2^n, \dots
u_{M+1}^n\rbrack^T$ to the $M$ degrees of freedom $\uu^n = \lbrack u_1^n,u_2^n,
\dots u_{M}^n\rbrack^T$, and $\uub_{b}^n = \lbrack 0,0, \dots ,0,\gamma(t^n)
\rbrack^T$ as a vector of length $M+1$ with zeros and the inhomogeneous
boundary condition $\gamma(t^n)$.  We have omitted $u_0^n$ from our vectors
($\uub^n, \uu^n$, and $\uub_{b}$) due to the homogeneous boundary condition at
$x=0$. The $M \times M+1$ restriction operator $R$ is an identity matrix, with
a column of zeros appended to it:
\begin{eqnarray}
\label{eq:rmat}
R &=& \begin{bmatrix}
1 & 0      & \ddots & \ddots & \ddots & 0 & 0 \\
0 & 1      & \ddots & \ddots & \ddots & 0 & 0 \\
\ddots & \ddots & \ddots & \ddots & \ddots & 0 & 0 \\
\ddots & \ddots & \ddots & \ddots & \ddots & 0 & 0 \\
\ddots & \ddots & \ddots & \ddots & \ddots & 0 & 0 \\
0 & \ddots & \ddots & \ddots & \ddots & 1 & 0 \\
\end{bmatrix}.
\end{eqnarray}
We use the standard 2nd-order accurate central finite difference operator for
$u_{xx}$ that is of size $M+1 \times M+1$, and is defined as $\bar{A}_{ii} =
2/\dx^2$ and $\bar{A}_{i-1,i}=\bar{A}_{i+1,i}=-1/\dx^2$, i.e.,
\begin{eqnarray}
\label{eq:amat}
\bar{A} &=&
\frac{1}{\dx^2} \begin{bmatrix}
2 & -1      & \ddots & \ddots & \ddots & 0 \\
-1 & 2      & -1 & \ddots & \ddots & 0 \\
\ddots & \ddots & \ddots & \ddots & \ddots & 0 \\
\ddots & \ddots & \ddots & \ddots & \ddots & 0 \\
\ddots & \ddots & \ddots & \ddots & \ddots & -1 \\
0 & \ddots & \ddots & \ddots & -1 & 2 \\
\end{bmatrix}.
\end{eqnarray}
Using a $k$th-order accurate BDF$k$ scheme for discretizing $u_t$, it is
straightforward to show that the solution of the unsteady heat equation is
time-advanced from $t^{n-1}$ to $t^n$ at the degrees of freedoms as
\begin{eqnarray}
\label{eq:discundif2}
\uu^n = -\sum_{l=1}^{k} \beta_{l} H^{-1} \uu^{n-l} - \nu \Delta t H^{-1} R
\bar{A} \uub_b^n,
\end{eqnarray}
where $\dt$ is the timestep size (assumed to be the same at all timesteps),
$\beta_{l}$ is the coefficient for the BDF$k$ scheme, $H = \beta_0 I + \nu
\Delta t A$ is the Helmholtz matrix, $I$ is a $M \times M$ identity matrix, and
$A=R\bar{A}R^T$.  The system in (\ref{eq:discundif2}) also depends on the
time-dependent inhomogeneous boundary condition, $\gamma(t^n)$.

\section{Unsteady heat equation in overlapping grids} \label{sec:unsteadyheatms}
\begin{figure}[t!]
\includegraphics[width=140mm]{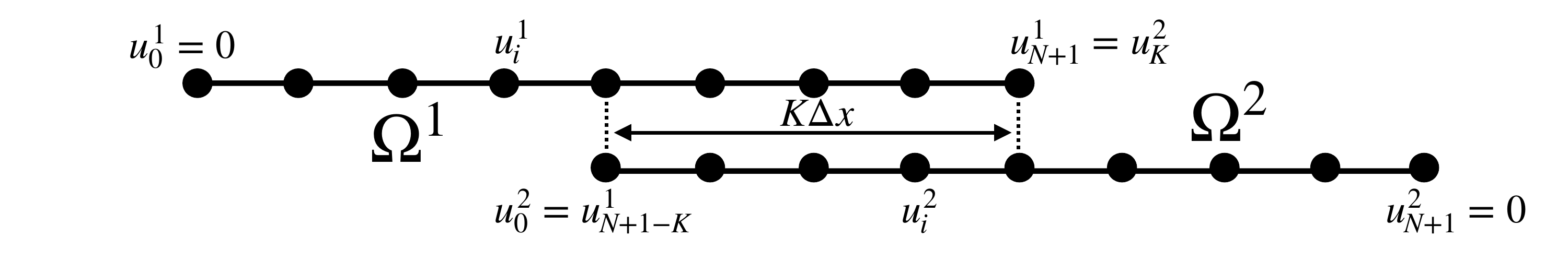}
\vspace{-4mm}
\caption{Overlapping grids with $N=7$ and $K=4$.}
\label{fig:ovdom} \end{figure}

For solving the unsteady diffusion equation with overlapping grids, we split
the monodomain grid ($\Omega$) into two grids ($\Omega^1$ and $\Omega^2$) with
equal number of grid points ($N+2$) and overlap width $K\Delta x$, such that
the grid points in the overlap $\Omega^1 \cap \Omega^2$ coincide. Note that we
use coincident grids because our focus is the temporal discretization of the
method and not the interpolation of interdomain boundary data.

Figure \ref{fig:ovdom} shows the overlapping grids obtained from the monodomain
grid of Fig. \ref{fig:monodom}, which are setup such that $M=2N-K+1$.
Algebraically, this decomposition is realized through restriction matrices,
$R_i$ that extract $N$ of the $M$ values from a given vector $\uu$ on $\Omega$
as $\uu^i = R_i \uu$:
\begin{eqnarray}
\label{rmat1}
R_1 =
\overbrace{\left[ \begin{array}{ccccccccc}
1 & 0      & \ddots & \ddots & \ddots & 0 & 0 & \dots & 0\\
0 & 1      & \ddots & \ddots & \ddots & 0 & 0 & \dots & 0\\
\ddots & \ddots & \ddots & \ddots & \ddots & 0 & 0 & \dots & 0\\
\ddots & \ddots & \ddots & \ddots & \ddots & 0 & 0 & \dots & 0\\
\ddots & \ddots & \ddots & \ddots & \ddots & 0 & 0 & \dots & 0\\
0 & \ddots & \ddots & \ddots & \ddots & 1 & 0 & \dots & 0
\end{array} \right]}^{M},\qquad
R_2 =
\overbrace{\left[ \begin{array}{ccccccccc}
0 & \dots & 0 & 1 & 0      & \ddots & \ddots & \ddots & 0\\
0 & \dots & 0 & 0 & 1      & \ddots & \ddots & \ddots & 0\\
0 & \dots & 0 & \ddots & \ddots & \ddots & \ddots & \ddots & 0\\
0 & \dots & 0 & \ddots & \ddots & \ddots & \ddots & \ddots & 0 \\
0 & \dots & 0 & \ddots & \ddots & \ddots & \ddots & \ddots & 0 \\
0 & \dots & 0 & 0 & \ddots & \ddots & \ddots & \ddots & 1
\end{array} \right]}^{M}. \\
\nonumber
\hexbrace{5cm}{N}\hexbrace{2.2cm}{M-N}\hspace{2.4cm}
\hexbrace{2.1cm}{M-N} \hexbrace{4.9cm}{N}
\end{eqnarray}
For simplicity, we impose homogeneous boundary conditions at the left boundary
of $\Omega^1$ ($u^{1,n}_0=0$) and right boundary of $\Omega^2$
($u^{2,n}_{N+1}=0$).  Here, we have used $\uu^s$ to represent the
solution at the $N$ degrees of freedom of $\Omega^s$, and $u^{s,n}_j$ to
represent the solution at $j$th grid point of $\Omega^s$ at $t^n$.  The
homogeneous boundary conditions on $\dO\backslash\dO^s_I$ allow us to use the
method developed for the monodomain grid (\ref{eq:discundif2}), with the
difference that the boundary data for the interdomain boundary grid points
($u^{1}_{N+1}$ and $u^{2}_{0}$) is obtained from the corresponding overlapping
grid in each subdomain. To effect this interdomain exchange, we define an
interpolation operator, $B_{ij}=(I-R_i^TR_i)R_j^T$, that extracts the value
from $\Omega^j$ at $\dO^i_I \cap \Omega^j$ and maps it to $\Omega^i$. The
operator $B_{ij}$ serves the same purpose as the high-order interpolation
functionality provided by \fpt\ in the Schwarz-SEM framework for interpolating
the interdomain boundary data \cite{gslibrepo}.

\subsection{Singlerate timestepping}
The singlerate timestepping scheme described in this section was originally
developed and analyzed by Peet and Fischer \cite{peet2012}, and has shown shown
to preserve the temporal convergence of the Schwarz-SEM framework for the
incompressible Navier-Stokes equations \cite{mittal2019nonconforming}.  In
singlerate timestepping, we use the same timestep size $\dt$ across all
overlapping subdomains to integrate the PDE of interest.

For notational purposes, we introduce $\uu^{s,n,q}$ to denote the solution $u$
in $\Omega^s$ at the $q$th Schwarz iteration at time $t^n$. Thus, assuming that
the solution is known up to $t^{n-1}$ and has been converged using $q= 1 \dots
Q$ Schwarz iterations, $\uu^{s,n-1,Q}$ denotes the converged solution at
$t^{n-1}$.  Similar to \eqref{eq:discundif2}, the solution of the unsteady heat
equation is advanced in time by using a BDF$k$ scheme to discretize the
time-derivative and replacing the fixed inhomogeneous boundary condition for
the monodomain case with the interpolated interdomain boundary data for the
overlapping grids.  Since the solution is only known up to $t^{n-1}$, the
initial Schwarz iterate ($q=0$, the {\em predictor step}) uses interdomain
boundary data based on $m$th-order extrapolation in time from the solution at
previous timesteps in the overlapping subdomain. The $Q$ subsequent Schwarz
iterations (the {\em corrector steps}) directly interpolate the interdomain
boundary data from the most recent iteration. Thus, the boundary data for
$\dO^s_I$ at the predictor and corrector steps is
\begin{eqnarray}
\label{eq:ibdstsq0}
q=0: \uu^{i,n,0}|_{\dO^i_I} &=& B_{ij}\sum_{l=1}^{m} \textw_{l} \uu^{j,n-l,Q}, \\
\label{eq:ibdstsqall}
q=1 \dots Q: \uu^{i,n,q}|_{\dO^i_I} &=& B_{ij}\uu^{j,n,q-1},
\end{eqnarray}
where $\textw_{l}$ are coefficients for the EXT$m$ scheme and $B_{ij}$ is the
interpolation operator. Using \eqref{eq:ibdstsq0}-\eqref{eq:ibdstsqall},
the PC approach for solving the
unsteady heat equation in overlapping subdomains is
\begin{eqnarray}
\label{eq:uq1}
q=0: \uu^{i,n,0} &=& - \sum_{l=1}^{k} \beta_{l} H_i^{-1} \uiiinmlq
 + \sum_{l=1}^{m} \textw_{l} H_i^{-1} J_{ij} \uu^{j,n-l,Q}, \\
\label{eq:uqall}
q=1 \dots Q: \! \uu^{i,n,q} &=& - \sum_{l=1}^{k} \beta_{l} H_i^{-1} \uiiinmlq
 + H_i^{-1} J_{ij} \uu^{j,n,q-1},
\end{eqnarray}
where
\begin{eqnarray}
\label{eq:discundif3}
\! H_{i} = \beta_{0} I_{i} + \nu \Delta t A_{i}, \quad
J_{ij} = - \nu \Delta t R_{i} A B_{ij}, \quad
A_i = R_i A R_i^T.
\end{eqnarray}
All the matrices in (\ref{eq:discundif3}) are of size $N \times N$ except the
restriction operator $R_i$ (\ref{rmat1}).  In Section
\ref{sec:singlestability}, we analyze the stability of this singlerate
timestepping scheme \eqref{eq:uq1}-\eqref{eq:uqall} using the matrix method for
stability analysis.

We note that the BDF$k$/EXT$m$ scheme in \eqref{eq:uq1}-\eqref{eq:uqall} will be
$k$th-order even if $m=1$, provided that sufficient corrector iterations are done.
By using $m$th order extrapolation at the predictor step, we significantly reduce
$Q$. Numerical experiments in the Schwarz-SEM framework show that a BDF3/EXT$1$
scheme can take as many as 50 corrector iterations at each time step to achieve
third-order temporal convergence. In contrast, a BDF3/EXT3 typically requires
only 3-5 corrector iterations.

\subsection{Multirate timestepping for the unsteady heat equation in overlapping grids} \label{sec:multirate}
In \cite{mittal2020multirate}, we extended the singlerate timestepping scheme
described in the previous section to support an arbitrary timestep ratio in an
arbitrary number of overlapping subdomains for the incompressible Navier-Stokes
equations. This novel multirate timestepping scheme is parallel-in-space and
allows each subdomain to integrate the PDE with a time-step based on the local
Courant–Friedrichs–Lewy (CFL) number of its mesh. We
have demonstrated that the MTS maintains the temporal accuracy of the underlying
BDFk/EXTk-based timestepper and accurately models complex turbulent flow and
heat transfer phenomenon even when the timestep ratio is as high as 50.
Here, we use the MTS scheme developed for the
Schwarz-SEM framework to solve the unsteady heat equation in overlapping grids
and analyze it using the matrix method for stability analysis.

\begin{figure}[t!] \begin{center}
$\begin{array}{c}
\includegraphics[width=100mm]{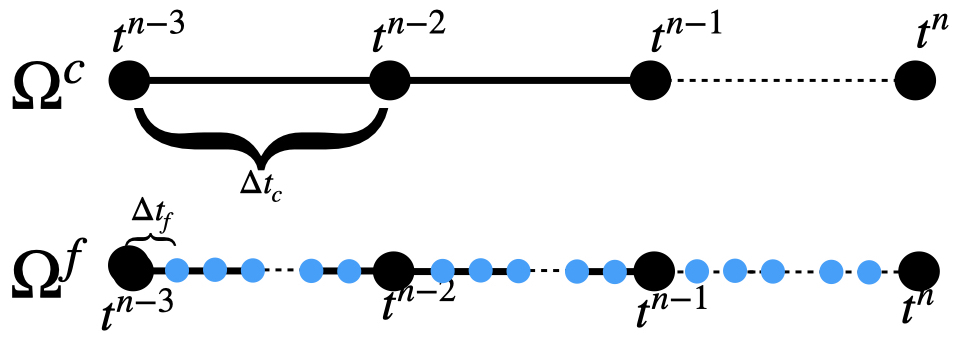}
\end{array}$
\end{center}
\caption{Schematic showing discrete time-levels for the multirate timestepping
with an arbitrary timestep ratio.}
\label{fig:arbitraryschematic}
\end{figure}

For MTS, we consider only integer timestep ratios,
\begin{eqnarray}
\tsr := \frac{\dt_c}{\dt_f} \in \mathbb{Z}^{+},
\end{eqnarray}
where $\dt_c$ corresponds to the subdomain ($\Omega^c$) with larger/coarser
time-scales and $\dt_f$ corresponds to the subdomain ($\Omega^f$) with
smaller/faster time-scales.  Figure \ref{fig:arbitraryschematic} shows a
schematic of the discrete time-levels for the MTS scheme with an arbitrary
timestep ratio. Here, the black circles ($\CIRCLE$) indicate the timestep
levels for both $\Omega^f$ and $\Omega^c$ and the blue circles
(\tikzcircle[cyan,fill=cyan]{2.5pt}) indicate the sub-timestep levels for
$\Omega^f$.

\begin{figure}[b!] \begin{center}
$\begin{array}{cc}
\includegraphics[width=70mm]{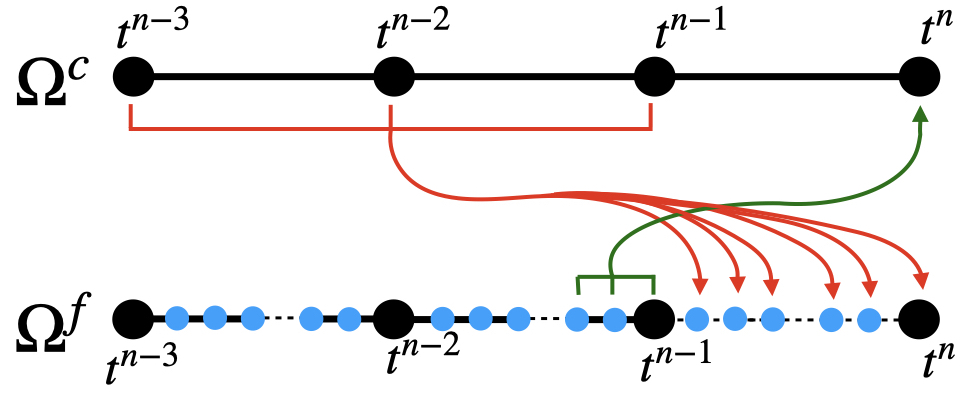} &
\includegraphics[width=70mm]{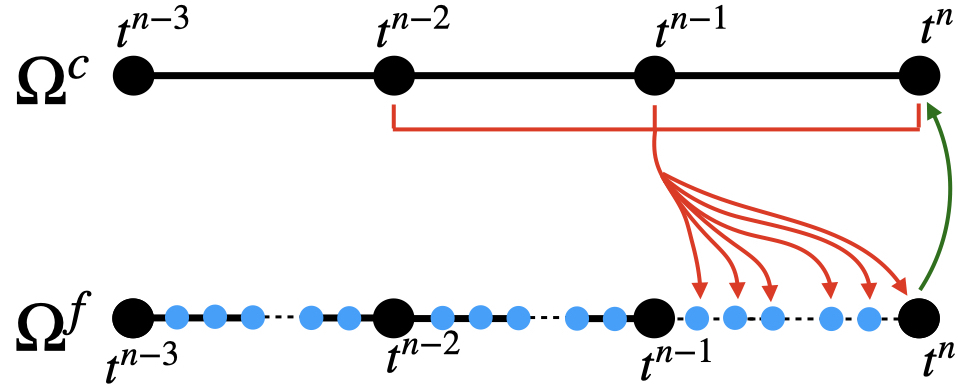}
\end{array}$
\end{center}
\caption{Schematic showing the dependence of the interdomain boundary data for
the (left) predictor and (right) corrector iterations for an arbitrary timestep ratio.}
\label{fig:schematic_ibd_pred_cor_eta_dep}
\end{figure}

Due to the difference in the timestep size, $\Omega^f$ has $\tsr$ sub-timesteps
and $\Omega^c$ has a single timestep to integrate the solution from $t^{n-1}$
to $t^n$.  Similar to the singlerate timestepping scheme, high-order temporal
accuracy is achieved in the multirate setting by extrapolating the interdomain
boundary data obtained from the solution at previous (sub-) timesteps. For the
predictor step, the interdomain boundary data dependency is shown in Fig.
\ref{fig:schematic_ibd_pred_cor_eta_dep}(left):
\begin{eqnarray}
\label{eq:usmultietacp}
i=1:\uu^{c,n,0} \rvert_{\dO^c_I} &=& B_{cf}\bigg(\sum_{j=1}^{m} \textw_{1j}\,
\uu^{f,n-1-\frac{j-1}{\tsr},Q}\bigg).\\
\label{eq:usmultietafp}
i=1\dots \tsr:\uu^{f,n-1+\frac{i}{\tsr},0} \rvert_{\dO^f_I} &=& B_{fc}\bigg(\sum_{j=1}^{m} \textw_{ij}\,
\uu^{c,n-j,Q}\bigg).
\end{eqnarray}
Note that the extrapolation coefficients $\textw_{ij}$ required for
extrapolating the interdomain boundary data for each sub-timestep $\Omega^f$ are determined
using the routines described in \cite{fornberg1998practical}.

After the predictor iteration, $Q$ corrector iterations are done where the
solution in $\Omega^f$ is obtained from $\Omega^c$ using the converged solution
at previous timesteps ($\uu^{c,n-1,Q}$ and $\uu^{c,n-2,Q}$) and the solution
from the most recent iteration at the current timestep ($\uu^{c,n,q-1}$).  For
the only timestep of $\Omega^c$, the solution is directly interpolated from the
most recent iteration ($\uu^{f,n,q-1}$). The interdomain boundary data dependency
for the $q=1\dots Q$ corrector iterations is shown in Fig.
\ref{fig:schematic_ibd_pred_cor_eta_dep}(right):
\begin{eqnarray}
\label{eq:usmultietacc}
i=1:\uu^{c,n,q} \rvert_{\dO^c_I} &=& B_{cf}\uu^{f,n,q-1},\\
\label{eq:usmultietafc}
i=1\dots \tsr:\uu^{f,n-1+\frac{i}{\tsr},q} \rvert_{\dO^f_I} &=& B_{fc}
\bigg(\tintw_{i1}\uu^{c,n,q-1} + \tintw_{i2}\uu^{c,n-1,Q} +
\tintw_{i3}\uu^{c,n-2,Q}\bigg).
\end{eqnarray}
Note that we compute the interpolation coefficients $\tintw_{ij}$ assuming
linear interpolation when $m=1$ or $2$, and quadratic interpolation when $m=3$.
This approach ensures that the desired temporal accuracy $\mathcal{O}(\dt^m)$
is maintained for the BDF$k$/EXT$m$-based scheme.

Equations \eqref{eq:usmultietacp}-\eqref{eq:usmultietafc} are used to solve the
unsteady heat equation with multirate timestepping in two overlapping grids
with a formulation similar to \eqref{eq:uq1}-\eqref{eq:uqall}.  The predictor
step ($q=0$) for the only timestep of $\Omega^c$ and $\tsr$ sub-timesteps of
$\Omega^f$ is:
\begin{eqnarray}
\label{eq:uqmspc}
\uu^{c,n,0} &=& - \sum_{l=1}^{k} \beta_{l} H_c^{-1} \uu^{c,n-l,Q}
 + \sum_{l=1}^{m} \textw_{l} H_c^{-1} J_{cf} \uu^{f,n-1-\frac{l-1}{\tsr},Q}, \\
 \label{eq:uqmspf}
 i=1\dots \tsr: \uu^{f,n-1+\frac{i}{\tsr},0} &=& - \sum_{l=1}^{k} \beta_{l} H_f^{-1} \uu^{f,n-1+\frac{i-l}{\tsr},Q}
  + \sum_{l=1}^{m} \textw_{il} H_f^{-1} J_{fc} \uu^{c,n-l,Q},
\end{eqnarray}
and the $Q$ corrector iterations are
\begin{eqnarray}
\label{eq:uqmscc}
\uu^{c,n,q} &=& - \sum_{l=1}^{k} \beta_{l} H_c^{-1} \uu^{c,n-l,Q}
 +  H_c^{-1} J_{cf} \uu^{f,n,q-1}, \\
 \label{eq:uqmscf}
 i=1\dots \tsr: \uu^{f,n-1+\frac{i}{\tsr},q} &=& - \sum_{l=1}^{k} \beta_{l} H_f^{-1} \uu^{f,n-1+\frac{i-l}{\tsr},Q} \\\nonumber
  &&+ H_f^{-1} J_{fc} \bigg(\tintw_{i1}\uu^{c,n,q-1} + \tintw_{i2}\uu^{c,n-1,Q} +
  \tintw_{i3}\uu^{c,n-2,Q}\bigg).
\end{eqnarray}
The $H$ and $J$ operators used here for multirate timestepping are similar to
the operators of singlerate timestepping \eqref{eq:discundif3}, with the only
difference that the timestep size used in their construction is different for
each subdomain.

We analyze the stability of the multirate timestepping scheme described here
\eqref{eq:uqmspc}-\eqref{eq:uqmscf} in Section \ref{sec:multistability} using
the matrix method for stability analysis.

\section{Stability of singlerate PC method} \label{sec:singlestability}

In this section, we analyze the stability of the singlerate PC scheme by
casting it as $\uz^n = G \uz^{n-1}$, where $G=C^{Q}P$ is
a product of the matrix $P$ corresponding to the predictor step with the EXT$m$
scheme (\ref{eq:uq1}) and matrix $C$ corresponding to the $Q$ corrector steps
(\ref{eq:uqall}).

Introducing the notation $\uz^{n,q}$ to denote the vector of solutions in both
overlapping subdomains
\begin{eqnarray}
\label{eq:uzsingle}
\uz^{n,q} = \big[\uu^{1,n,q^T} \,
\uu^{2,n,q^T} \, \uu^{1,n-1,Q^T} \, \uu^{2,n-1,Q^T} \, \uu^{1,n-2,Q^T} \,
\uu^{2,n-2,Q^T} \, \uu^{1,n-3,Q^T} \, \uu^{2,n-3,Q^T} \big]^T,
\end{eqnarray}
the predictor step \eqref{eq:uq1} is
\small
\begin{eqnarray}
\label{eq:pmat}
\underbrace{\left[ \begin{array}{c} \uu^{1,n,0} \\ \uu^{2,n,0} \\
\uu^{1,n-1,Q} \\ \uu^{2,n-1,Q} \\
\uu^{1,n-2,Q} \\ \uu^{2,n-2,Q} \\
\uu^{1,n-3,Q} \\ \uu^{2,n-3,Q} \\
\end{array} \right]}_{ \let\scriptstyle\textstyle
    \substack{\uz^{n,0}}} =
\underbrace{
\setlength\arraycolsep{2pt}
\begin{pmatrix}
-\beta_1 \ihina & \textw_1 \ihina J_{12} & -\beta_2 \ihina & \textw_2 \ihina J_{12} & -\beta_3 \ihina
& \textw_3 \ihina J_{12} & 0 & 0\\
\textw_1 \ihinb J_{21} & -\beta_1 \ihinb & \textw_2 \ihinb J_{21} & -\beta_2 \ihinb &
\textw_3 \ihinb J_{21} & -\beta_3 \ihinb & 0 & 0 \\
I_1 & 0   & 0   & 0   & 0   & 0   & 0 & 0\\
0   & I_2 & 0   & 0   & 0   & 0   & 0 & 0\\
0   & 0   & I_1 & 0   & 0   & 0   & 0 & 0\\
0   & 0   & 0   & I_2 & 0   & 0   & 0 & 0\\
0   & 0   & 0   & 0   & I_1 & 0   & 0 & 0\\
0   & 0   & 0   & 0   & 0   & I_2 & 0 & 0\\
\end{pmatrix}}_{\let\scriptstyle\textstyle
    \substack{P}}
\underbrace{\left[ \begin{array}{c}
\uu^{1,n-1,Q} \\ \uu^{2,n-1,Q} \\
\uu^{1,n-2,Q} \\ \uu^{2,n-2,Q} \\
\uu^{1,n-3,Q} \\ \uu^{2,n-3,Q} \\
\uu^{1,n-4,Q} \\ \uu^{2,n-4,Q} \end{array}
\right]}_{ \let\scriptstyle\textstyle \substack{\uz^{n-1,Q}}}
\end{eqnarray}
\normalsize
where the interdomain boundary data is extrapolated from the previous timesteps
in the overlapping subdomain.  Similarly, the matrix $C$ corresponding to
the $Q$ corrector iterations from \eqref{eq:uqall} is given by \eqref{eq:qmat}
where the solution $z^{n,q}$ depends on the interdomain boundary data interpolated
from the most recent Schwarz iteration ($z^{n,q-1}$). Using (\ref{eq:pmat}) and
(\ref{eq:qmat}), the solution is advanced in time as $z^{n,Q} = G z^{n-1,Q}$,
$G=C^QP$, and the spectral radius of $G$ is used to determine the stability of
the PC scheme for different grid sizes, overlap widths, extrapolation orders,
and corrector iterations. We note that using the parameters $N$ and $K$, the
grid overlap width is determined as $K\dx=K/(2N-K+2)$.
\small
\begin{eqnarray}
\label{eq:qmat}
\underbrace{\left[ \begin{array}{c}
\uu^{1,n,q} \\ \uu^{2,n,q} \\
\uu^{1,n-1,Q} \\ \uu^{2,n-1,Q} \\
\uu^{1,n-2,Q} \\ \uu^{2,n-2,Q} \\
\uu^{1,n-3,Q} \\ \uu^{2,n-3,Q} \\ \end{array} \right]}_{ \let\scriptstyle\textstyle
    \substack{\uz^{n,q}}} =
\underbrace{
\setlength\arraycolsep{2pt}
\begin{pmatrix}
0 & \ihina J_{12} & -\beta_1 \ihina & 0 & -\beta_2 \ihina & 0 & -\beta_3 \ihina
& 0 \\
\ihinb J_{21} & 0 & 0 & -\beta_1 \ihinb & 0 & -\beta_2 \ihinb  &
0 & -\beta_3 \ihinb  \\
0 & 0 & I_1 & 0   & 0   & 0   & 0   & 0   \\
0 & 0 & 0   & I_2 & 0   & 0   & 0   & 0   \\
0 & 0 & 0   & 0   & I_1 & 0   & 0   & 0   \\
0 & 0 & 0   & 0   & 0   & I_2 & 0   & 0   \\
0 & 0 & 0   & 0   & 0   & 0   & I_1 & 0   \\
0 & 0 & 0   & 0   & 0   & 0   & 0   & I_2 \\
\end{pmatrix}}_{\let\scriptstyle\textstyle \substack{C}}
\underbrace{\left[ \begin{array}{c}
\uu^{1,n,q-1} \\ \uu^{2,n,q-1} \\
\uu^{1,n-1,Q} \\ \uu^{2,n-1,Q} \\
\uu^{1,n-2,Q} \\ \uu^{2,n-2,Q} \\
\uu^{1,n-3,Q} \\ \uu^{2,n-3,Q} \\
\end{array} \right]}_{ \let\scriptstyle\textstyle \substack{\uz^{n,q-1}}}
\end{eqnarray}
\normalsize
\setlength{\abovedisplayskip}{3pt}
\setlength{\belowdisplayskip}{3pt}

\subsection{Stability Results for different BDF$k$/EXT$m$ schemes} \label{subsec:bdfk/extm}
In this section, we look at the stability of different BDF$k$/EXT$m$ schemes
with $Q$. We set $N=32$ and $K=5$, and in each case, we vary the nondimensional
timestep size  ($\nu\dt/\dx^2$) to see how the spectral radius of the propogation
operator changes.

\begin{figure}[t!]
\begin{center}
$\begin{array}{ccc}
\includegraphics[height=40mm]{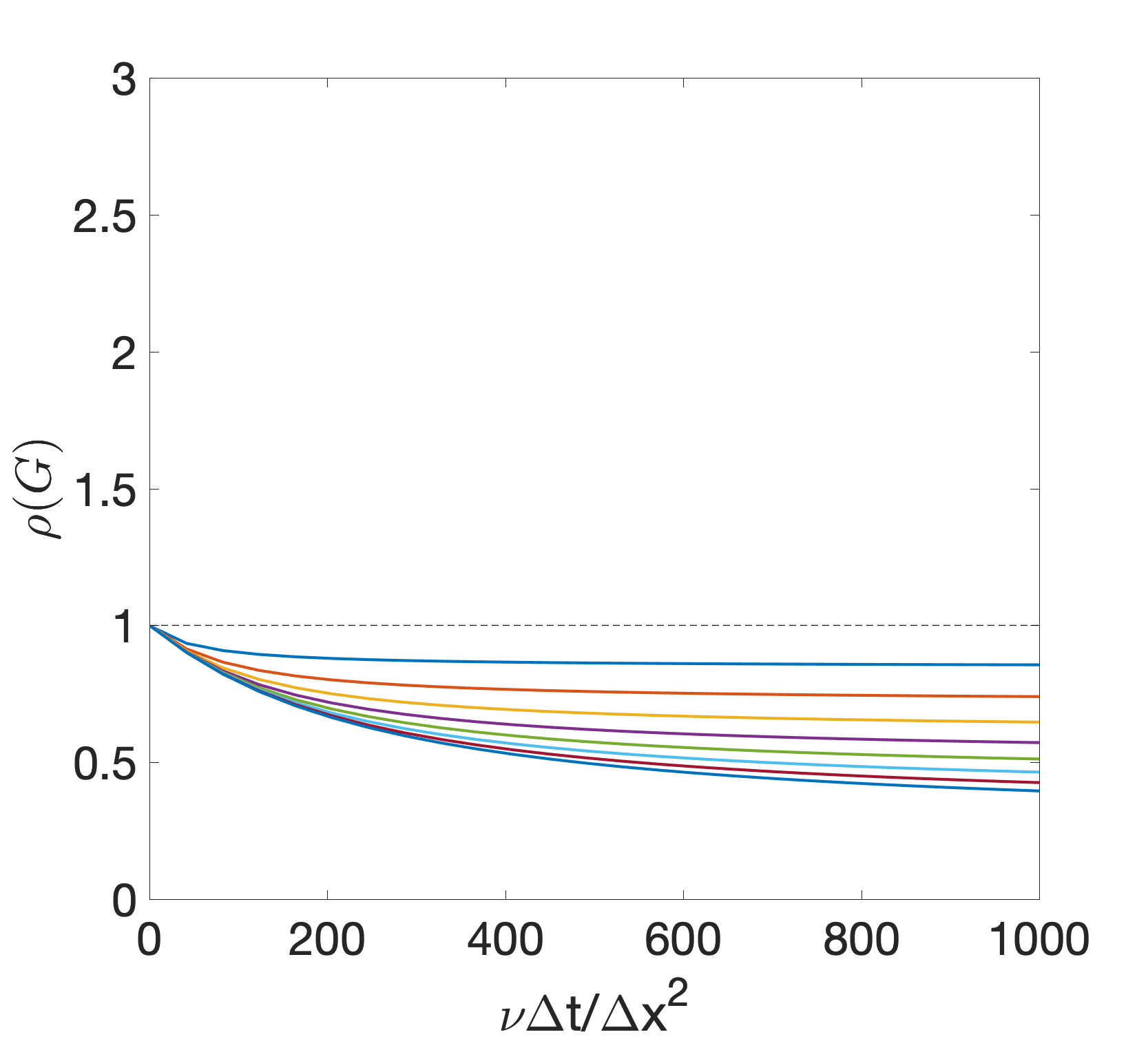} &
\includegraphics[height=40mm]{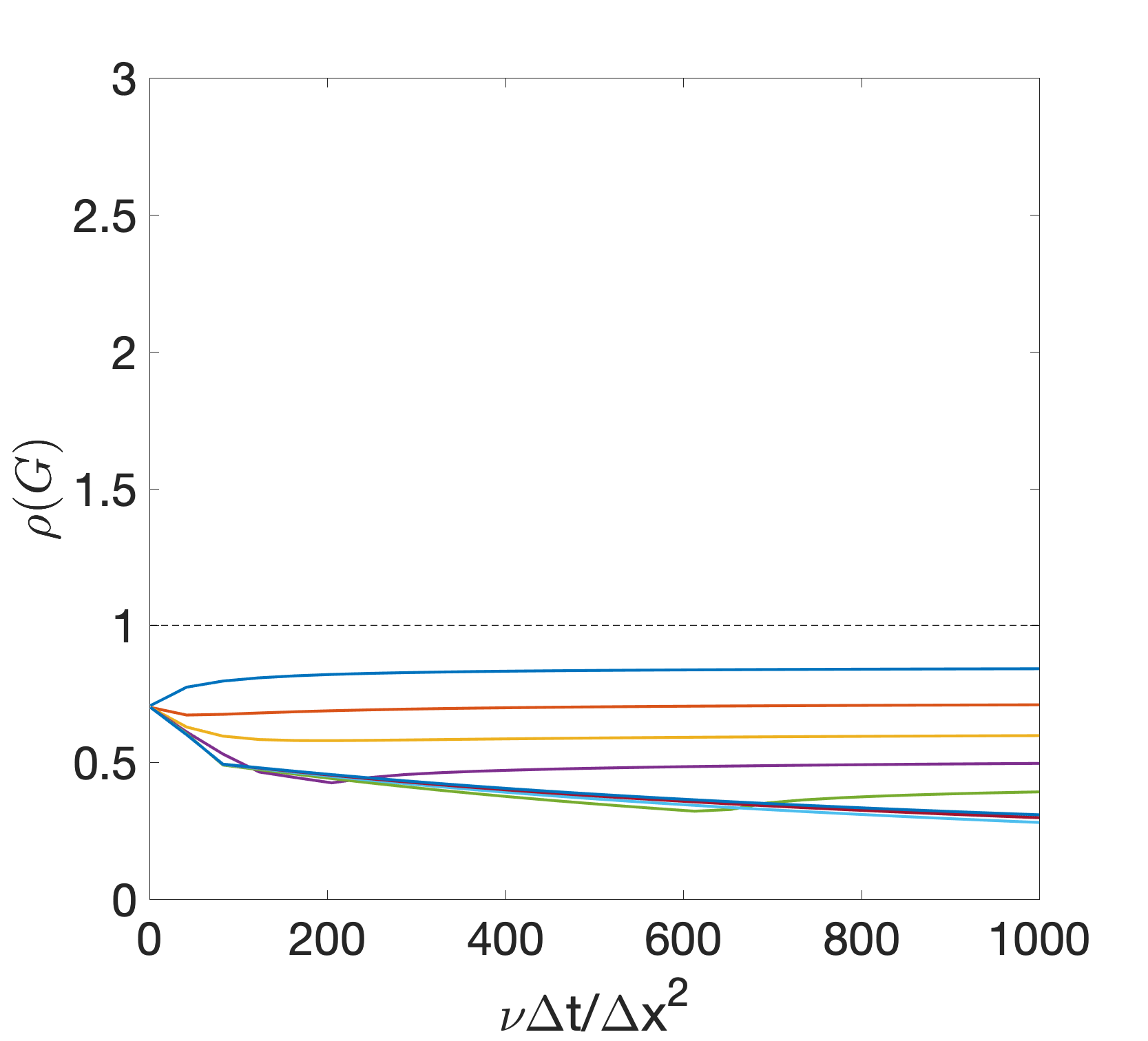} &
\includegraphics[height=40mm]{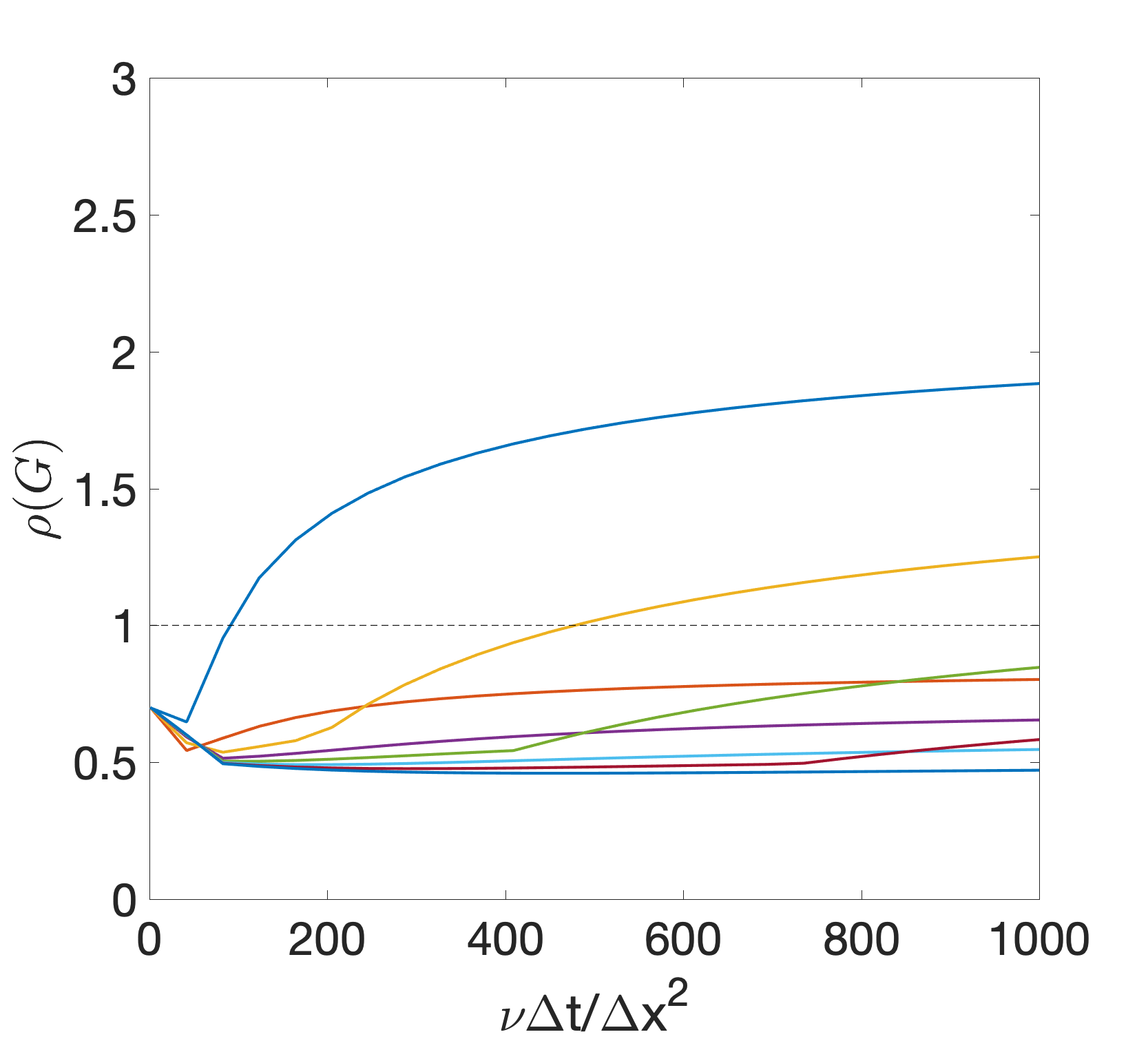} \\
\textrm{(a) BDF$1$/EXT$1$} &
\textrm{(b) BDF$2$/EXT$1$} &
\textrm{(c) BDF$2$/EXT$2$} \\
\includegraphics[height=40mm]{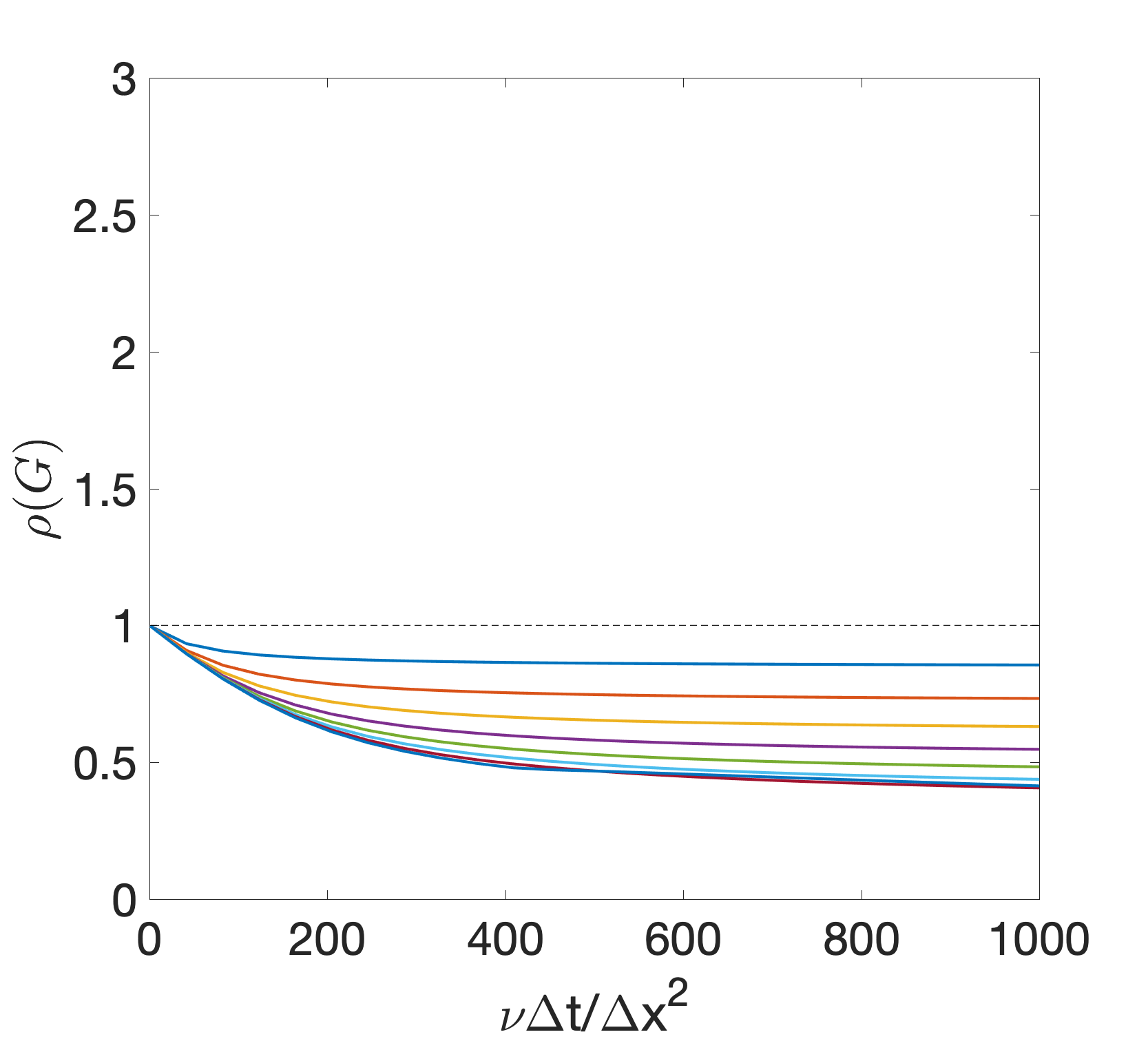} &
\includegraphics[height=40mm]{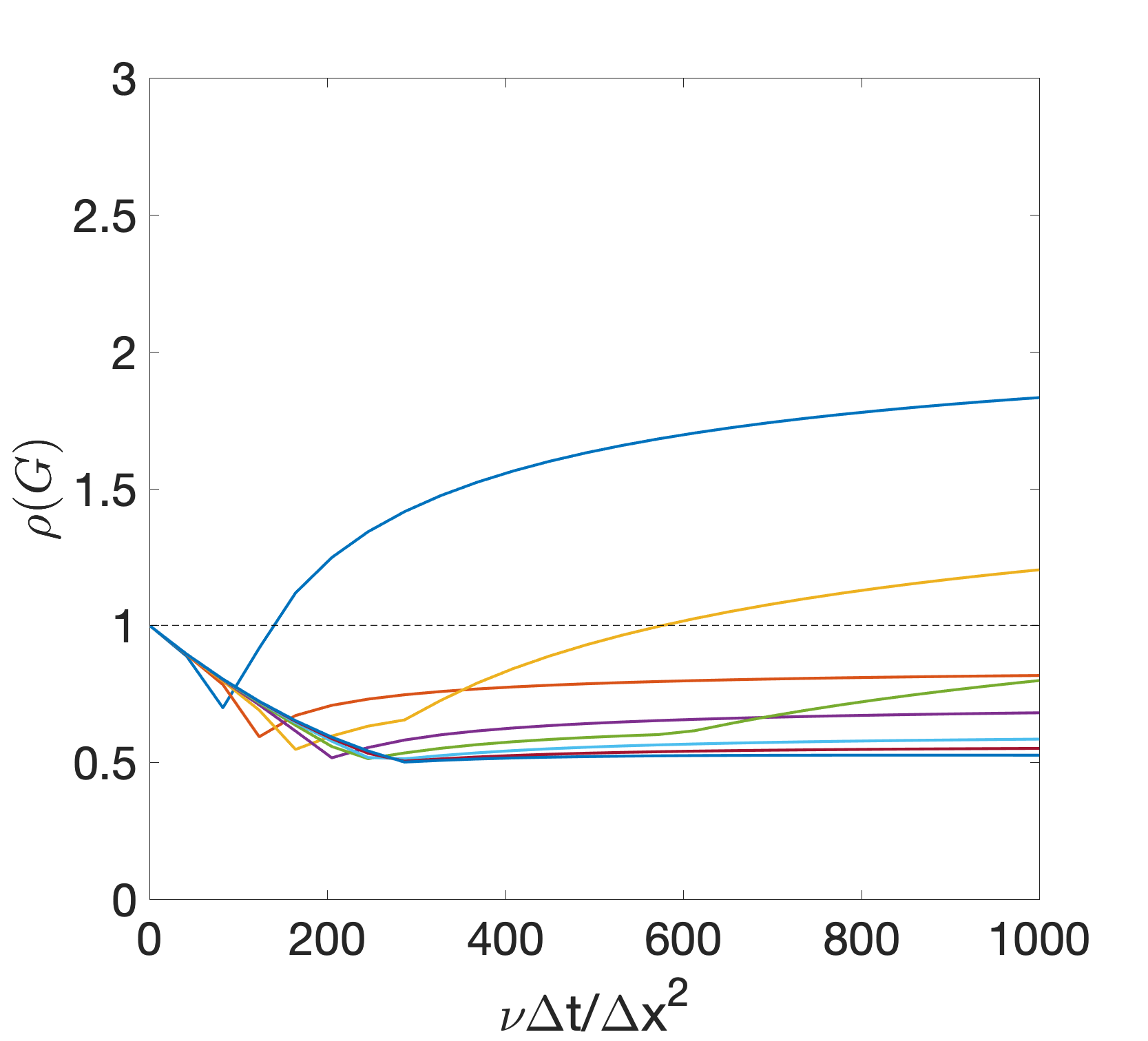} &
\includegraphics[height=40mm]{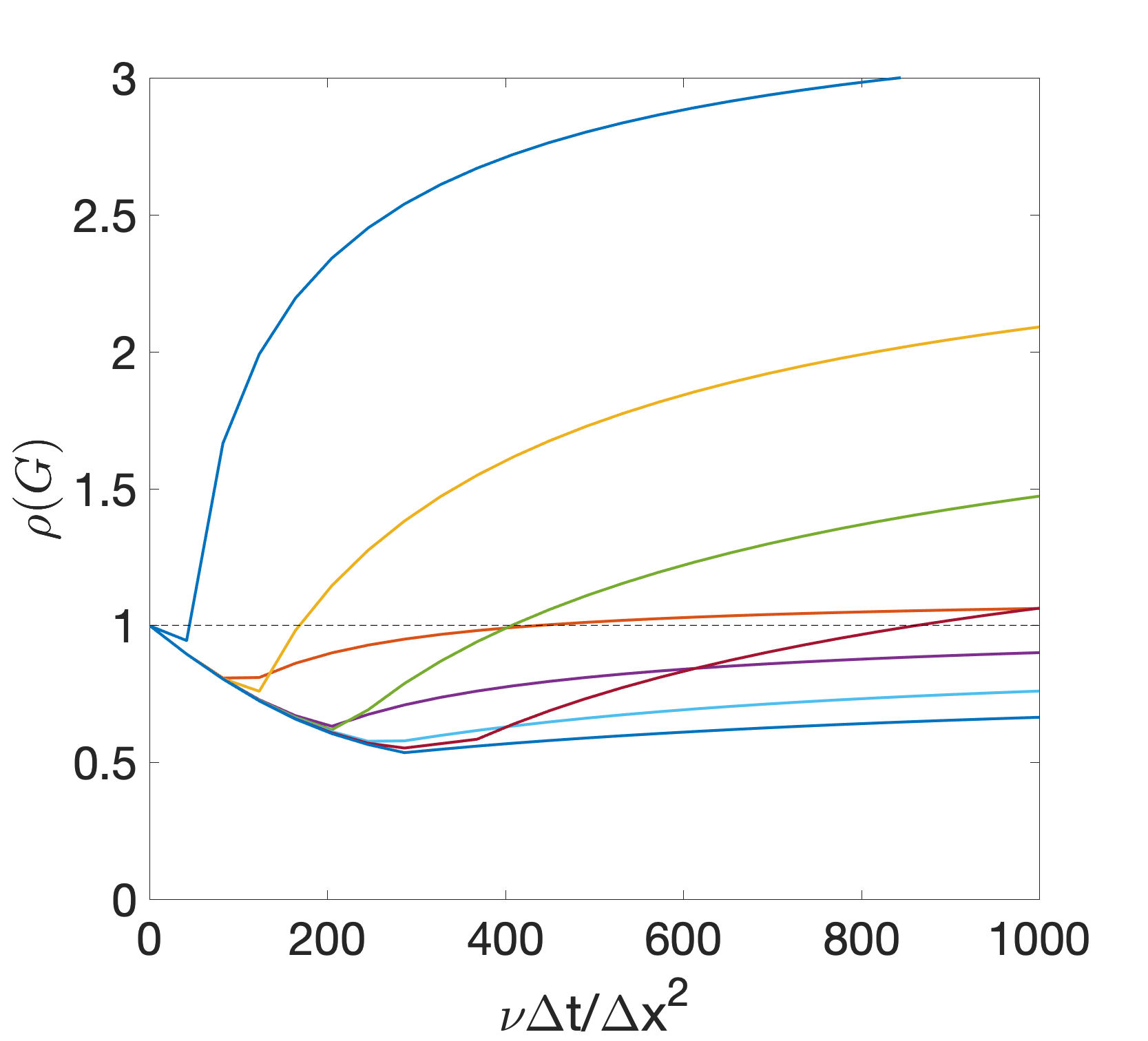} \\
\textrm{(d) BDF$3$/EXT$1$} &
\textrm{(e) BDF$3$/EXT$2$} &
\textrm{(f) BDF$3$/EXT$3$} \\
\multicolumn{3}{c}{\includegraphics[width=150mm]{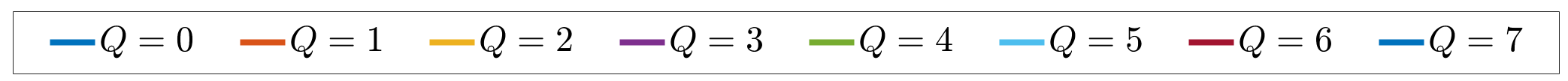}} \\
 \end{array}$
\end{center}
\vspace{-6mm}
\caption{Spectral radius $\rho(G)$ versus nondimensional time $\frac{\nu
\dt}{\dx^2}$ for different BDF$k$/EXT$m$ schemes with $Q=0\dots7$, $N=32$ and
$K=5$.} \label{fig:rhovss}
\end{figure}

Figure \ref{fig:rhovss} shows the spectral radius for the BDF$k$/EXT$m$ schemes
for $k=1\dots 3$ and $m \le k$. We observe that the
BDF$k$/EXT$1$ schemes are unconditionally stable, and the use of high-order
extrapolation ($m>1$) for interdomain boundary data requires correct iterations
for stability. These results are also indicated in \cite{peet2012}.  A novel
result that comes forth from Fig. \ref{fig:rhovss} is that for
high-order extrapolation, i.e., BDF$k$/EXT$2$ and BDF$k$/EXT$3$ schemes, odd
$Q$ is relatively more stable than even $Q$.  To the best of our knowledge,
this behavior has not been observed in the current literature for overlapping
grids. In terms of general ODEs and PDEs for a single domain, there is a sparse
evidence of such behavior. We will discuss this in the context of existing
literature and connect it to the Schwarz-SEM framework for INSE in Section
\ref{sec:oddevendiscussion}.
\clearpage

\subsubsection{Effect of increasing subdomain overlap}
The subdomain overlap width has an impact on the convergence of Schwarz-based
methods \cite{smith2004domain}.  For practical purposes, one would like to
minimize the overlap width to minimize the total number of elements needed for
modeling a domain.  Thus, we look at the impact of overlap width on the
stability of the PC scheme.  Since we are mainly interested in the high-order
extrapolation scheme ($m=1$ is unconditionally stable as shown in the previous
section), we look at the results for the BDF2/EXT2, BDF3/EXT2, and BDF3/EXT3
schemes, with $N=32$, and change the grid overlap parameter, $K$.

\begin{figure}[t!]
\begin{center}
$\begin{array}{ccc}
\includegraphics[height=40mm]{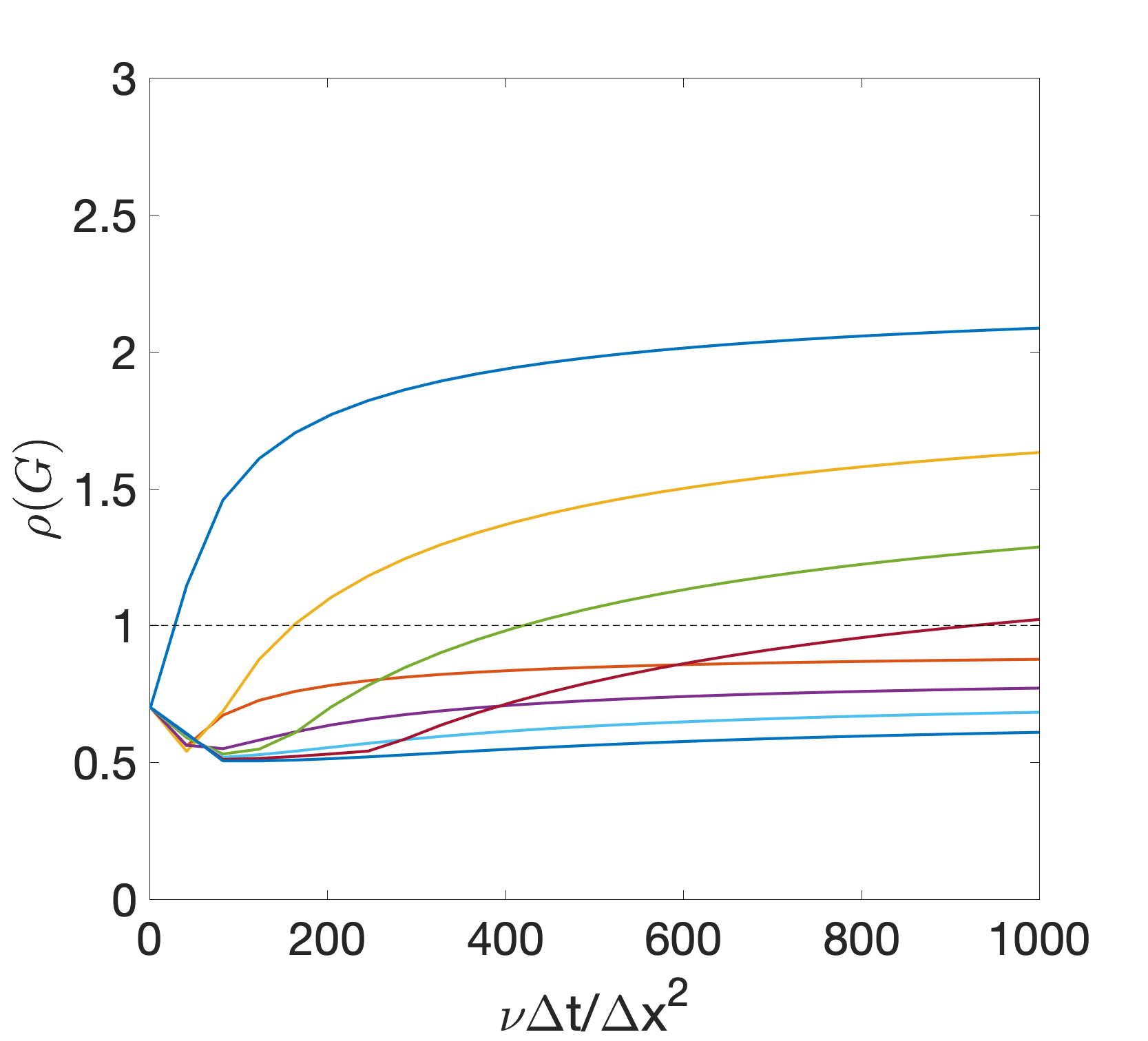} &
\includegraphics[height=40mm]{figures/rho_vs_s_n132_k5_BDF2_EXT2} &
\includegraphics[height=40mm]{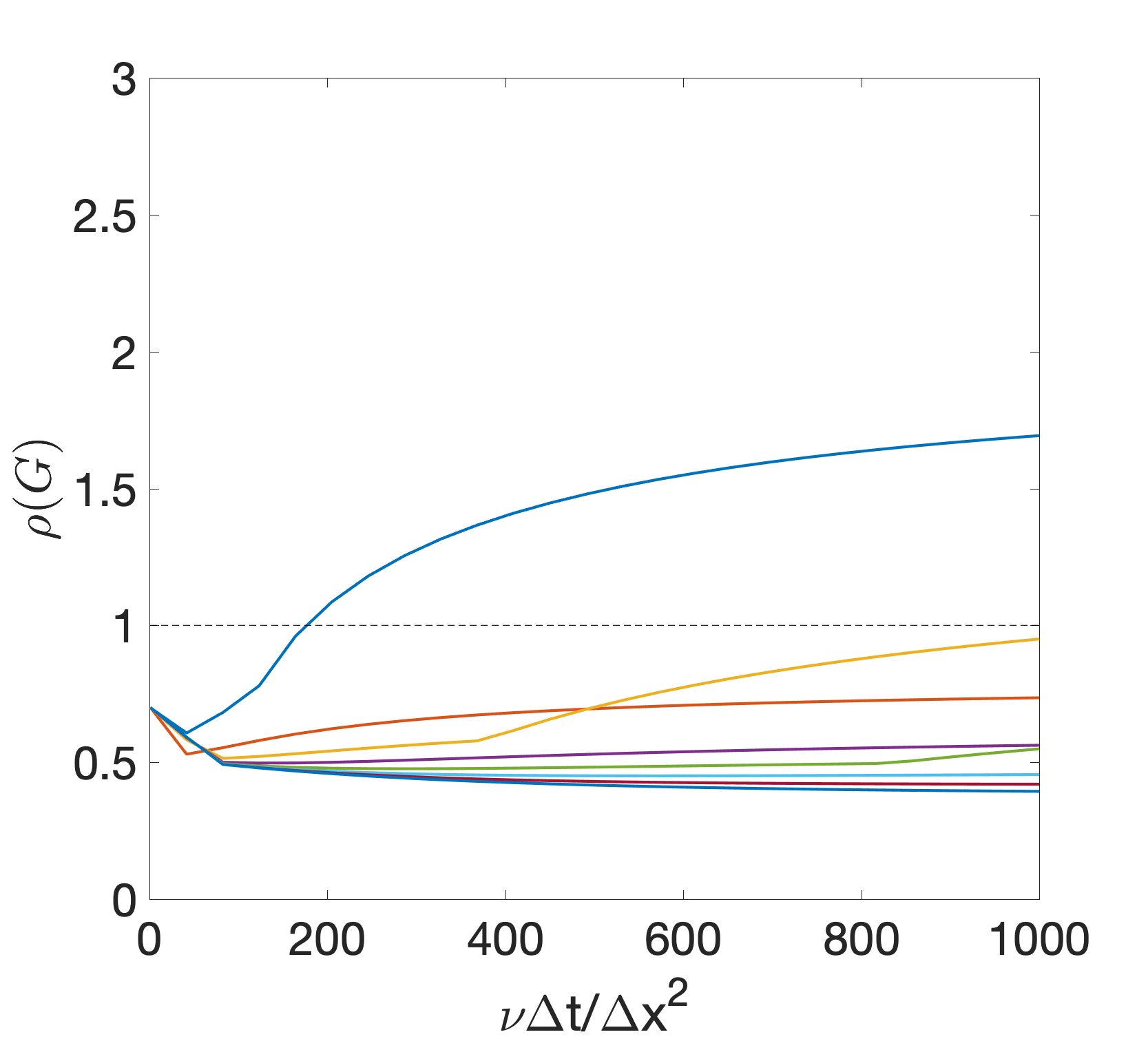} \\
\multicolumn{3}{c}{\textrm{(a) BDF2/EXT2}} \\
\includegraphics[height=40mm]{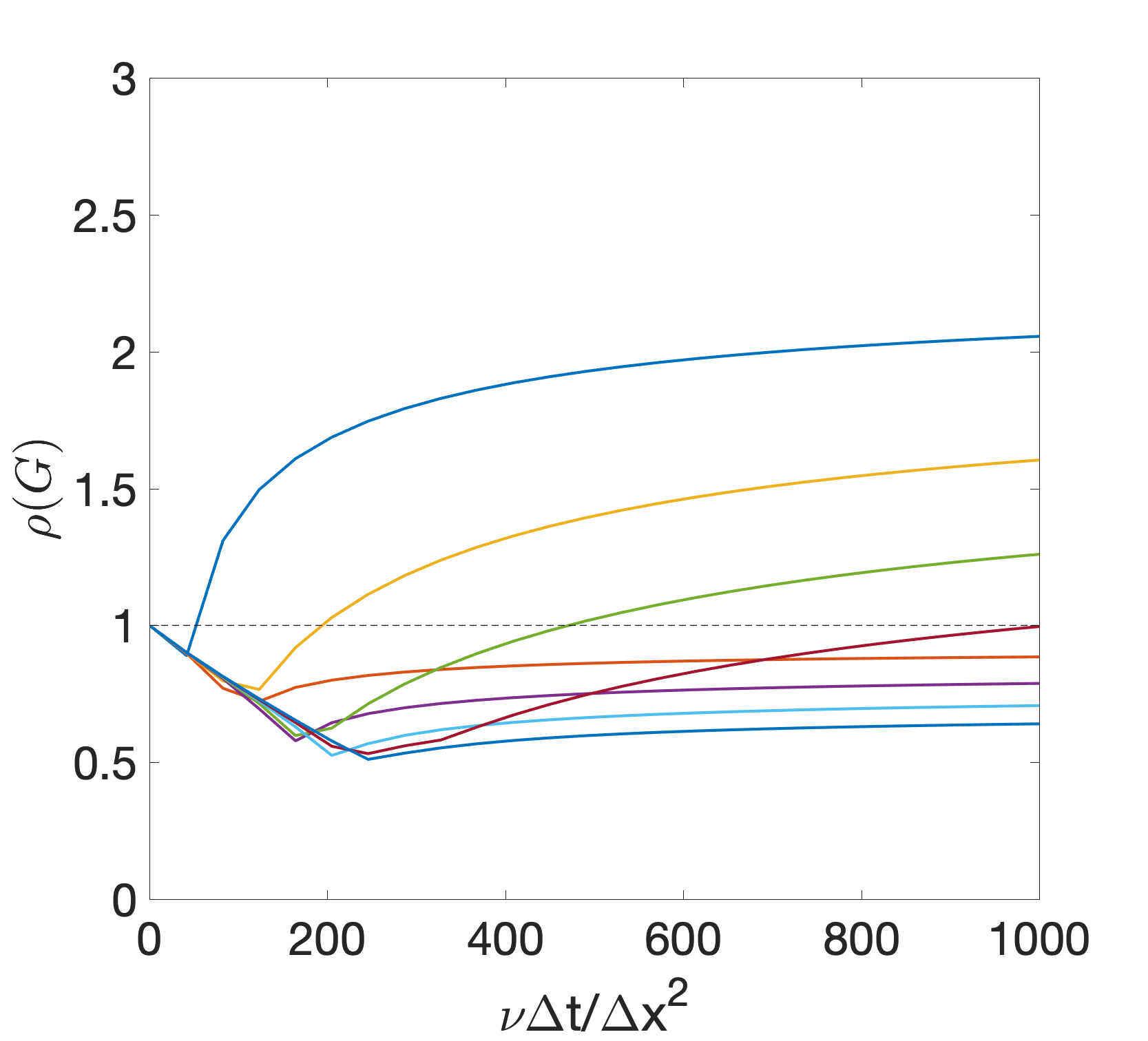} &
\includegraphics[height=40mm]{figures/rho_vs_s_n132_k5_BDF3_EXT2} &
\includegraphics[height=40mm]{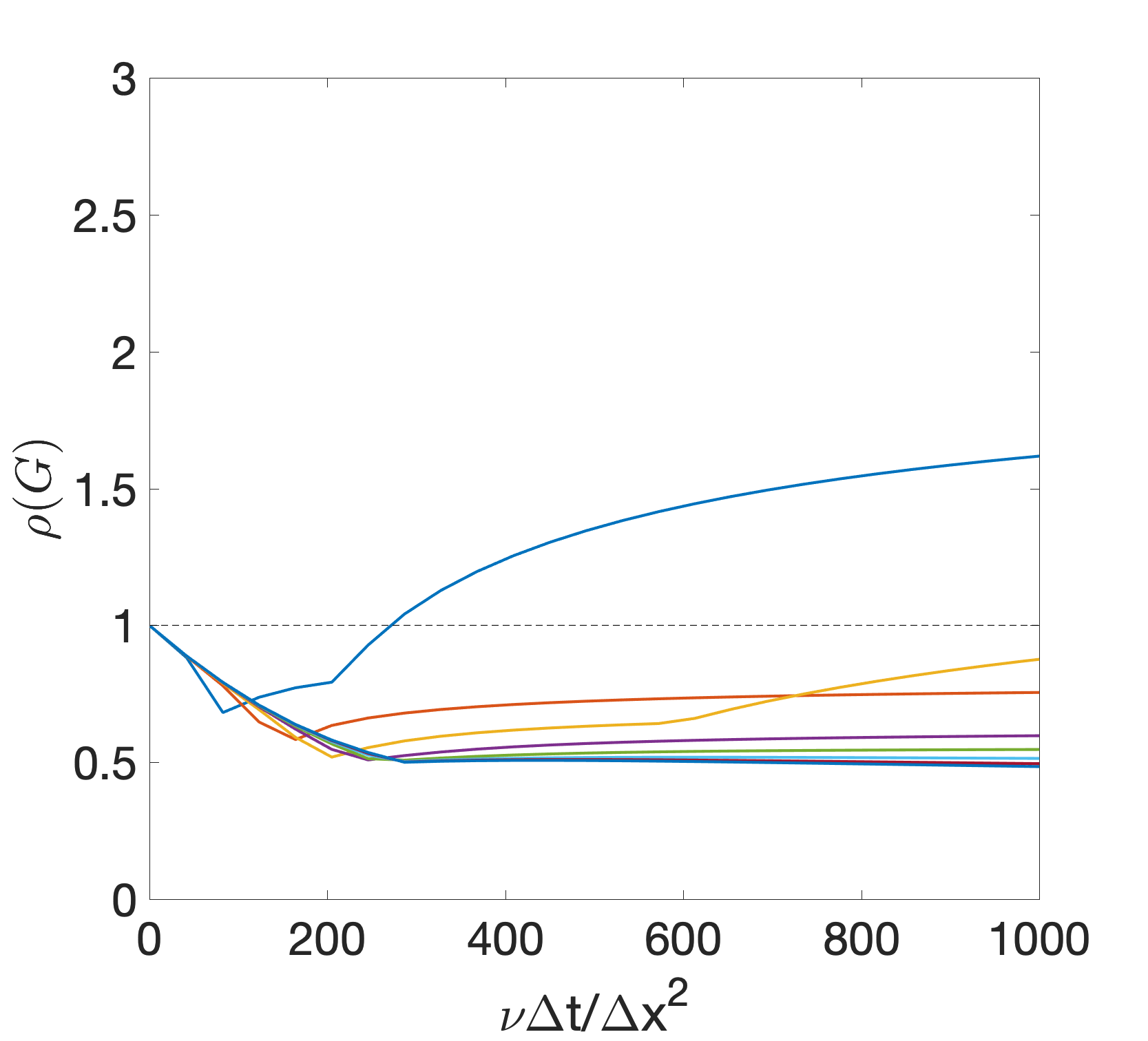} \\
\multicolumn{3}{c}{\textrm{(b) BDF3/EXT2}} \\
\includegraphics[height=40mm]{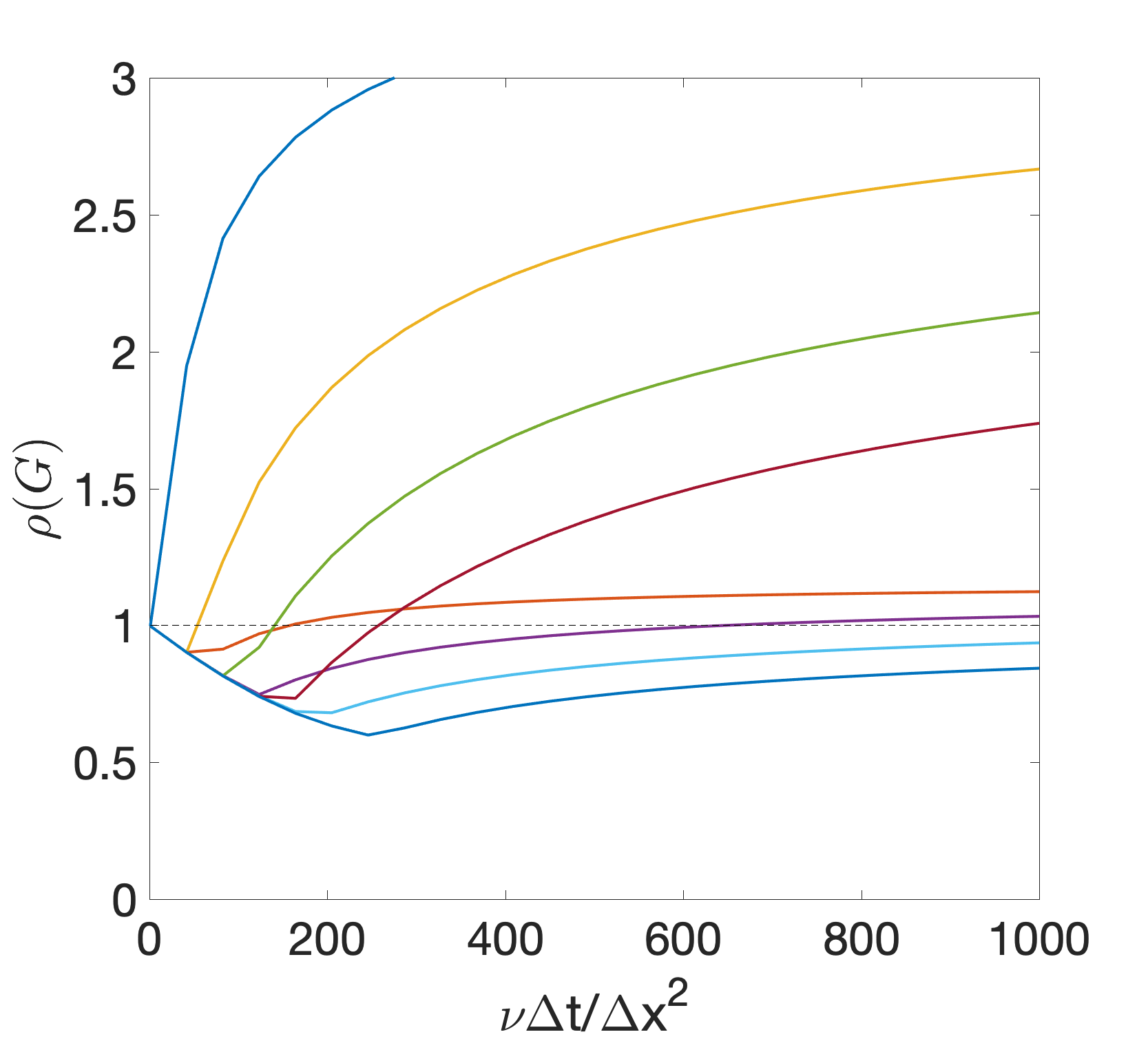} &
\includegraphics[height=40mm]{figures/rho_vs_s_n132_k5_BDF3_EXT3} &
\includegraphics[height=40mm]{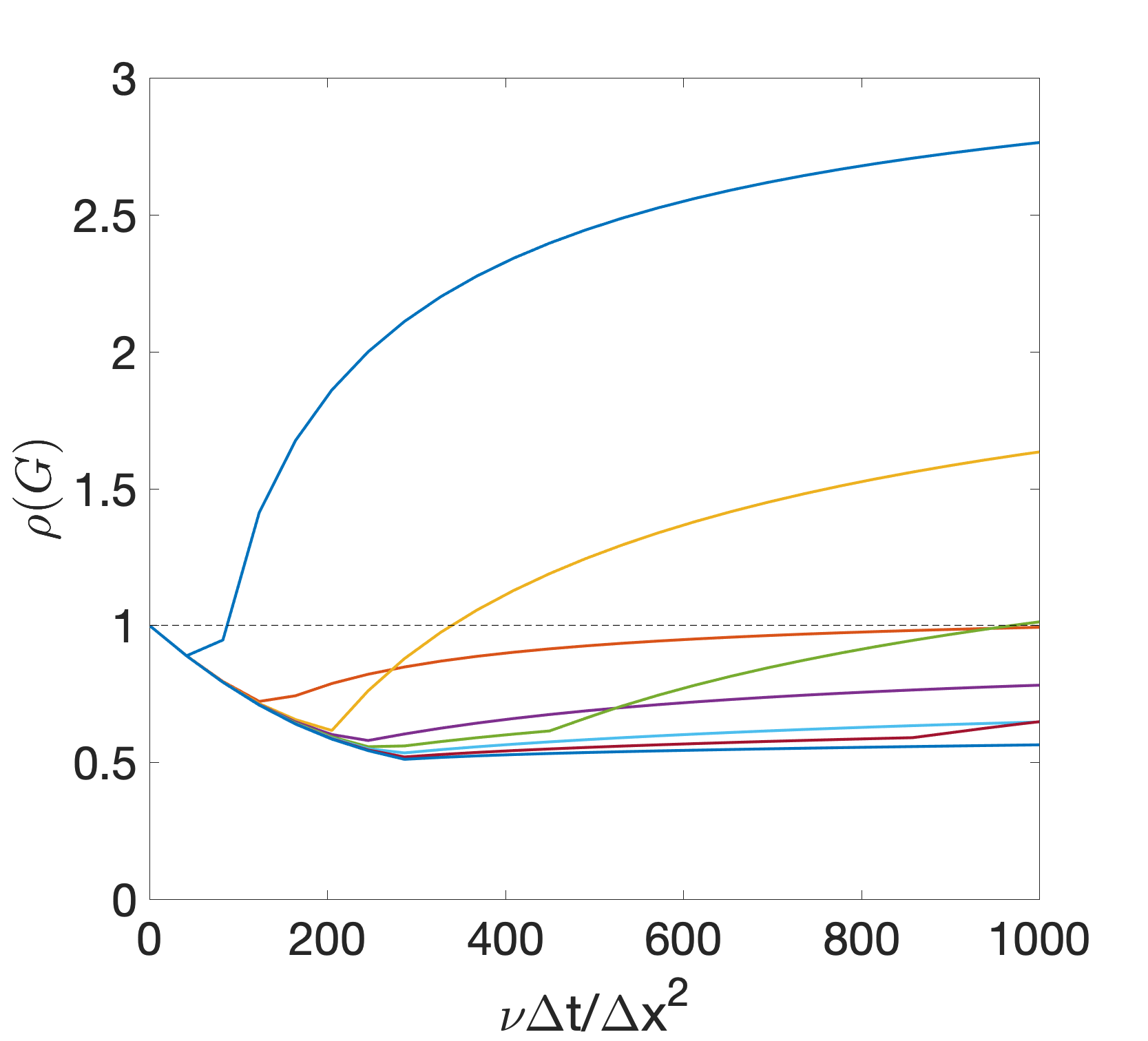} \\
\multicolumn{3}{c}{\textrm{(c) BDF3/EXT3}} \\
\textrm{(i) $N=32,K=3$} & \textrm{(ii) $N=32,K=5$} & \textrm{(iii) $N=32,K=7$} \\
\multicolumn{3}{c}{\includegraphics[width=140mm]{figures/legend}} \\
\end{array}$
\end{center}
\vspace{-6mm}
\caption{Spectral radius $\rho(G)$ versus nondimensional time $\frac{\nu
\dt}{\dx^2}$ for (a) BDF2/EXT2, (b) BDF3/EXT2, and (c) BDF3/EXT3 scheme with $N=32$
and varying $K$: (left to right) (i) $K=3$, (ii) $K=5$, and (iii) $K=7$.}
\label{fig:rhostab_overlap_k}
\end{figure}

Figure \ref{fig:rhostab_overlap_k} show that increasing the overlap width
increases the stability range over which a given scheme is stable for a
specified number of corrector iterations $Q$.  We also notice that high-order
extrapolation schemes are less stable as compared to their low order
counterparts, e.g., the BDF3/EXT3 scheme is less stable than the BDF3/EXT2
scheme, for a given $Q$.

\subsubsection{Effect of increasing grid resolution while keeping overlap fixed}
\begin{figure}[t!]
\begin{center}
$\begin{array}{ccc}
\includegraphics[height=40mm]{figures/rho_vs_s_n132_k5_BDF2_EXT2} &
\includegraphics[height=40mm]{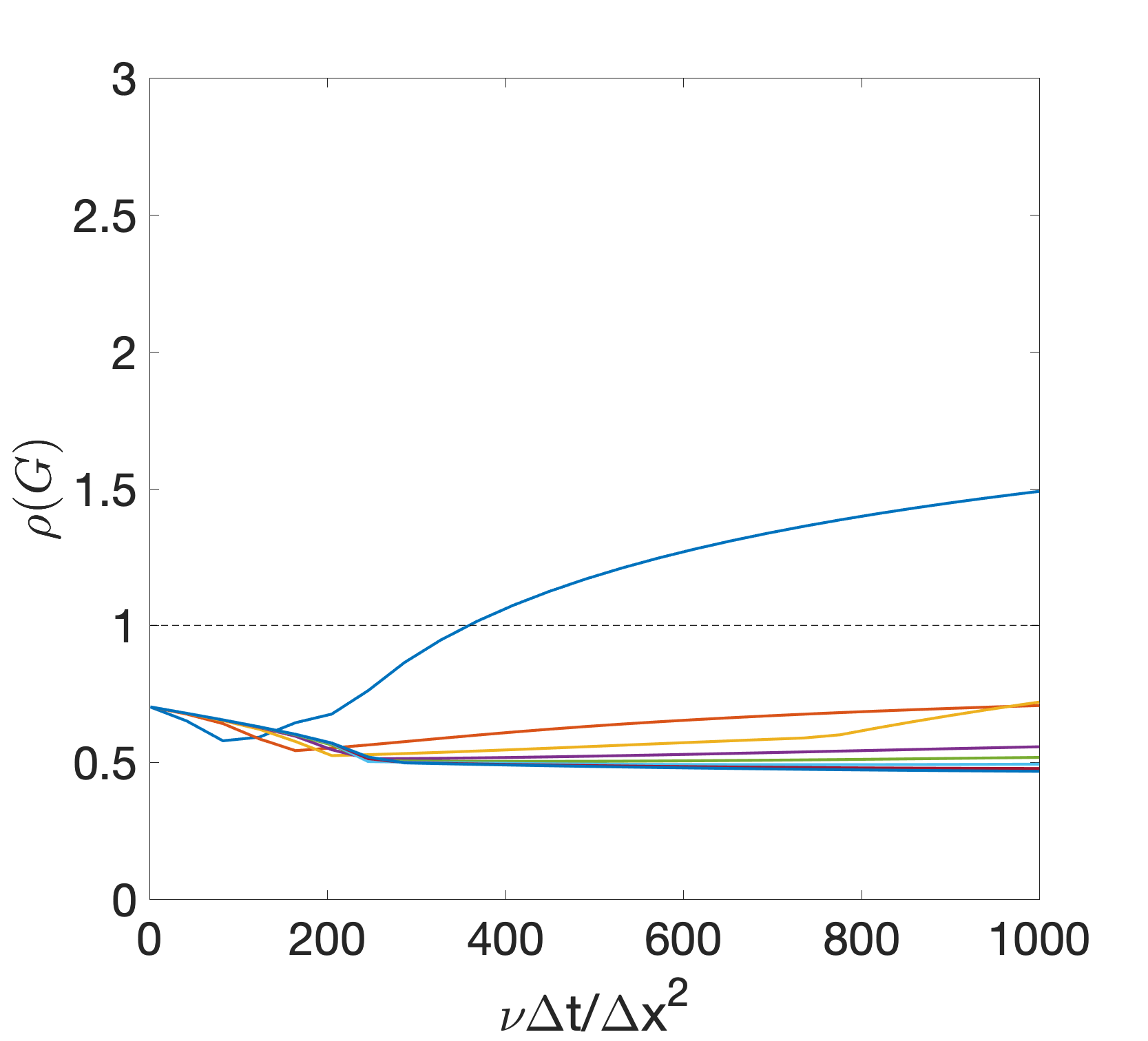} &
\includegraphics[height=40mm]{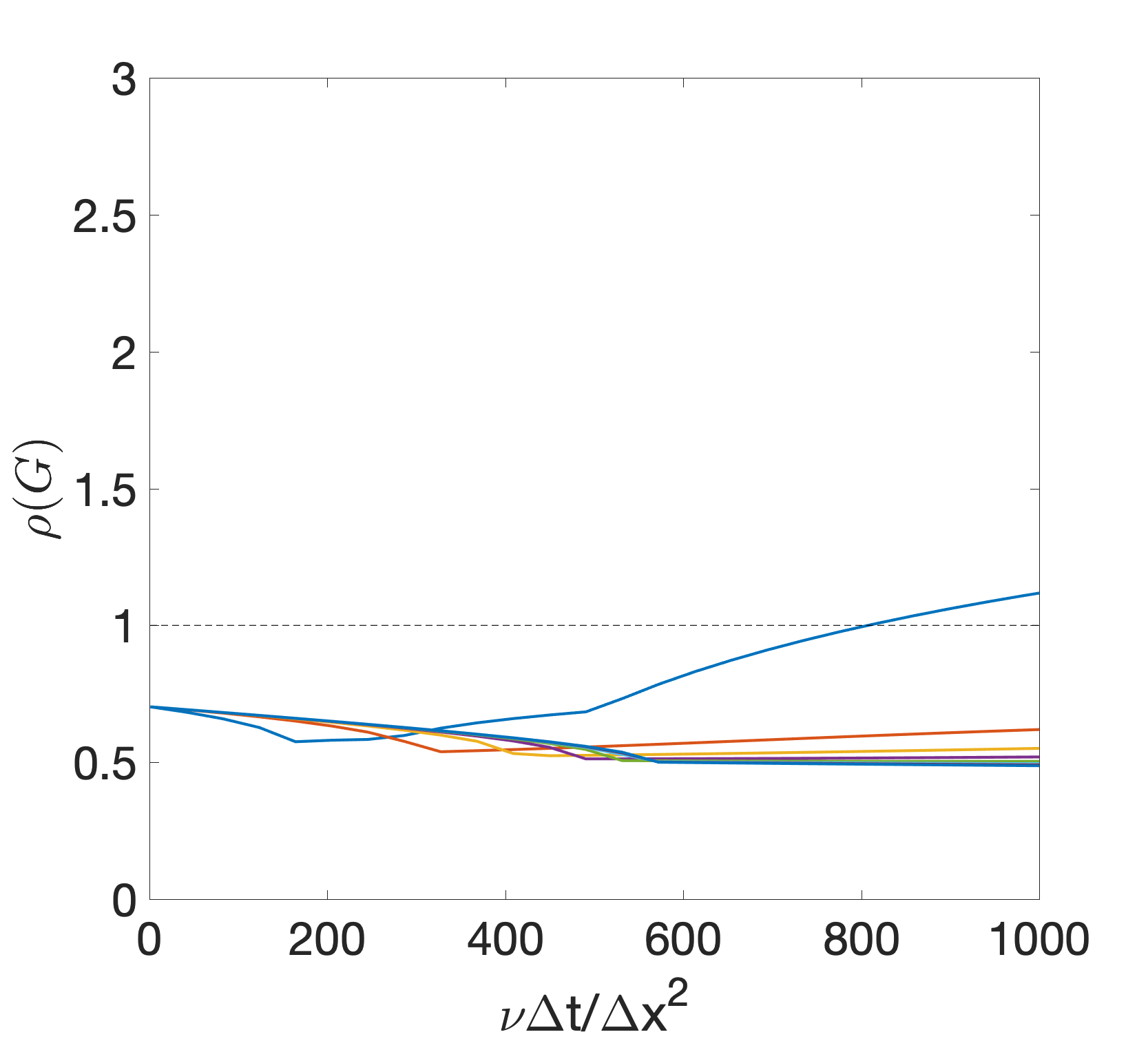} \\
\multicolumn{3}{c}{\textrm{(a) BDF2/EXT2}} \\
\includegraphics[height=40mm]{figures/rho_vs_s_n132_k5_BDF3_EXT2} &
\includegraphics[height=40mm]{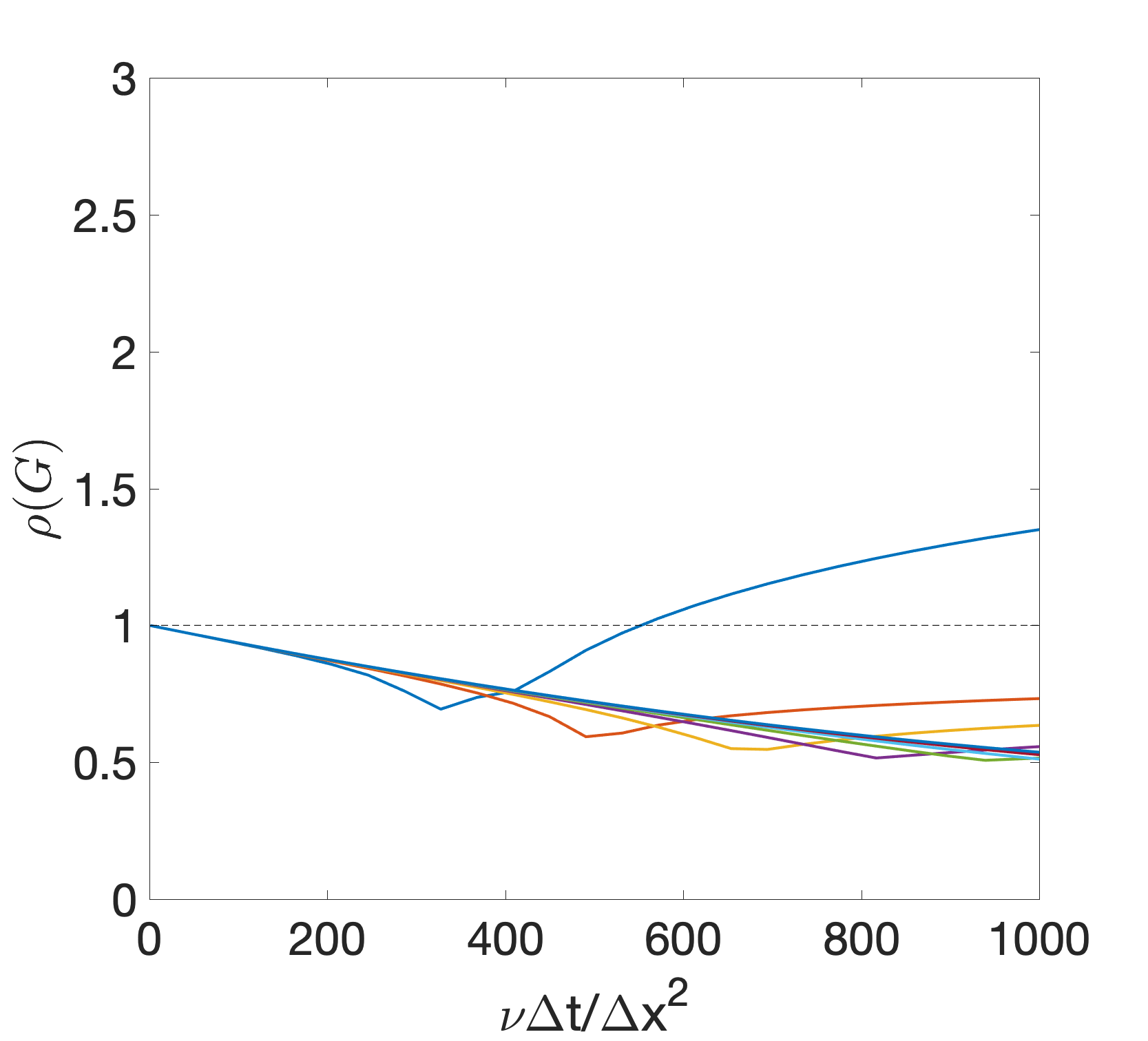} &
\includegraphics[height=40mm]{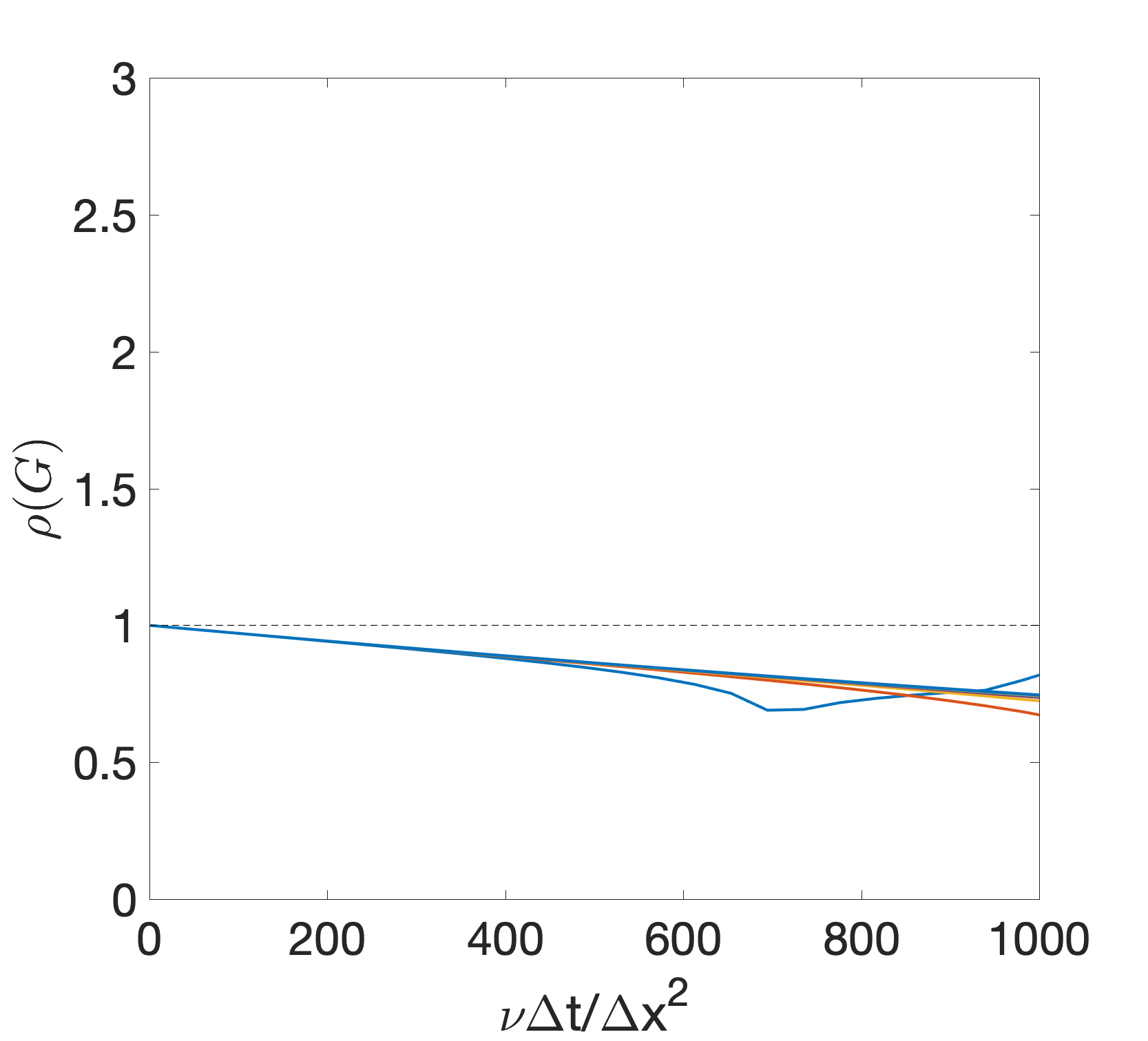} \\
\multicolumn{3}{c}{\textrm{(b) BDF3/EXT2}} \\
\includegraphics[height=40mm]{figures/rho_vs_s_n132_k5_BDF3_EXT3} &
\includegraphics[height=40mm]{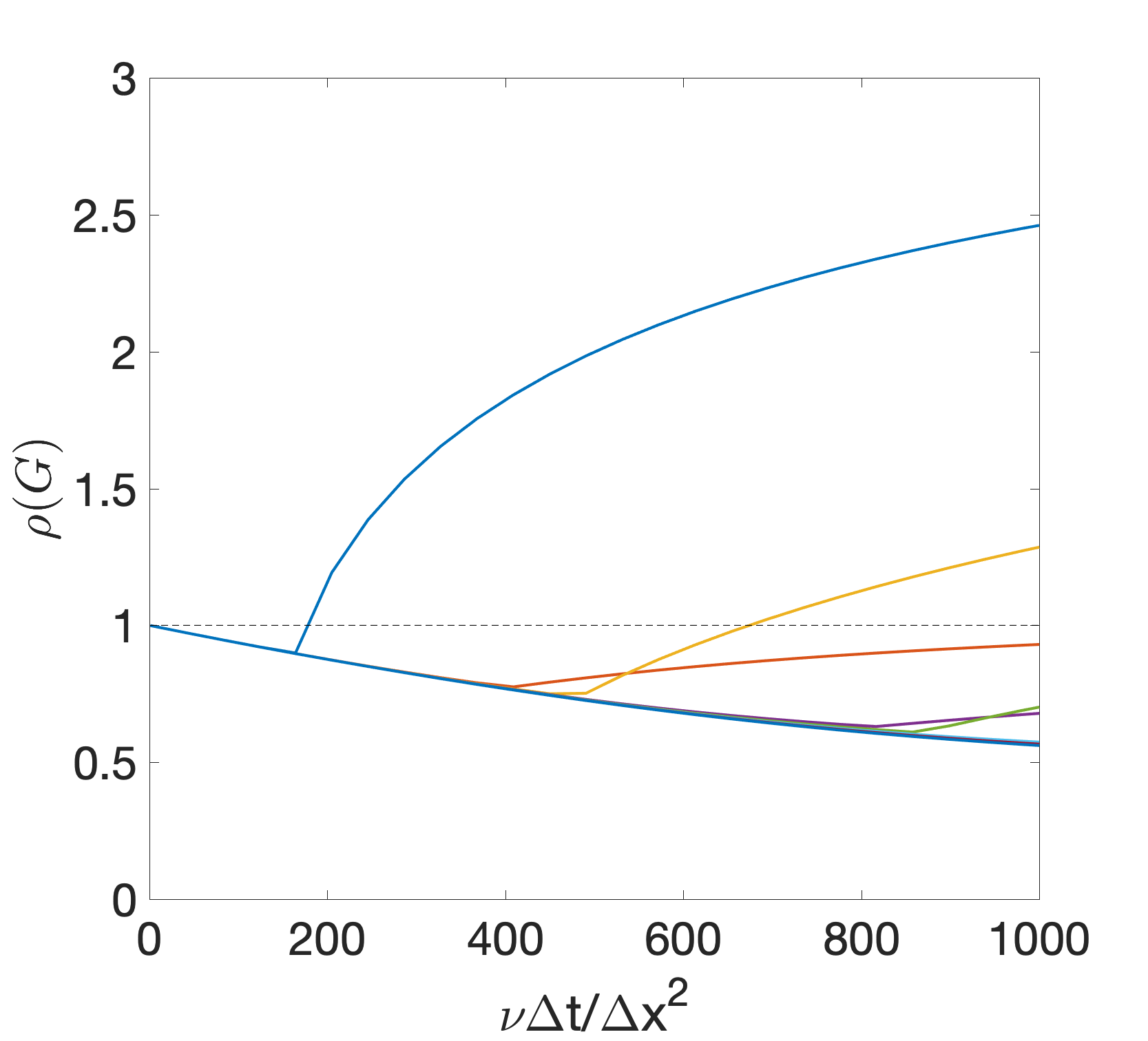} &
\includegraphics[height=40mm]{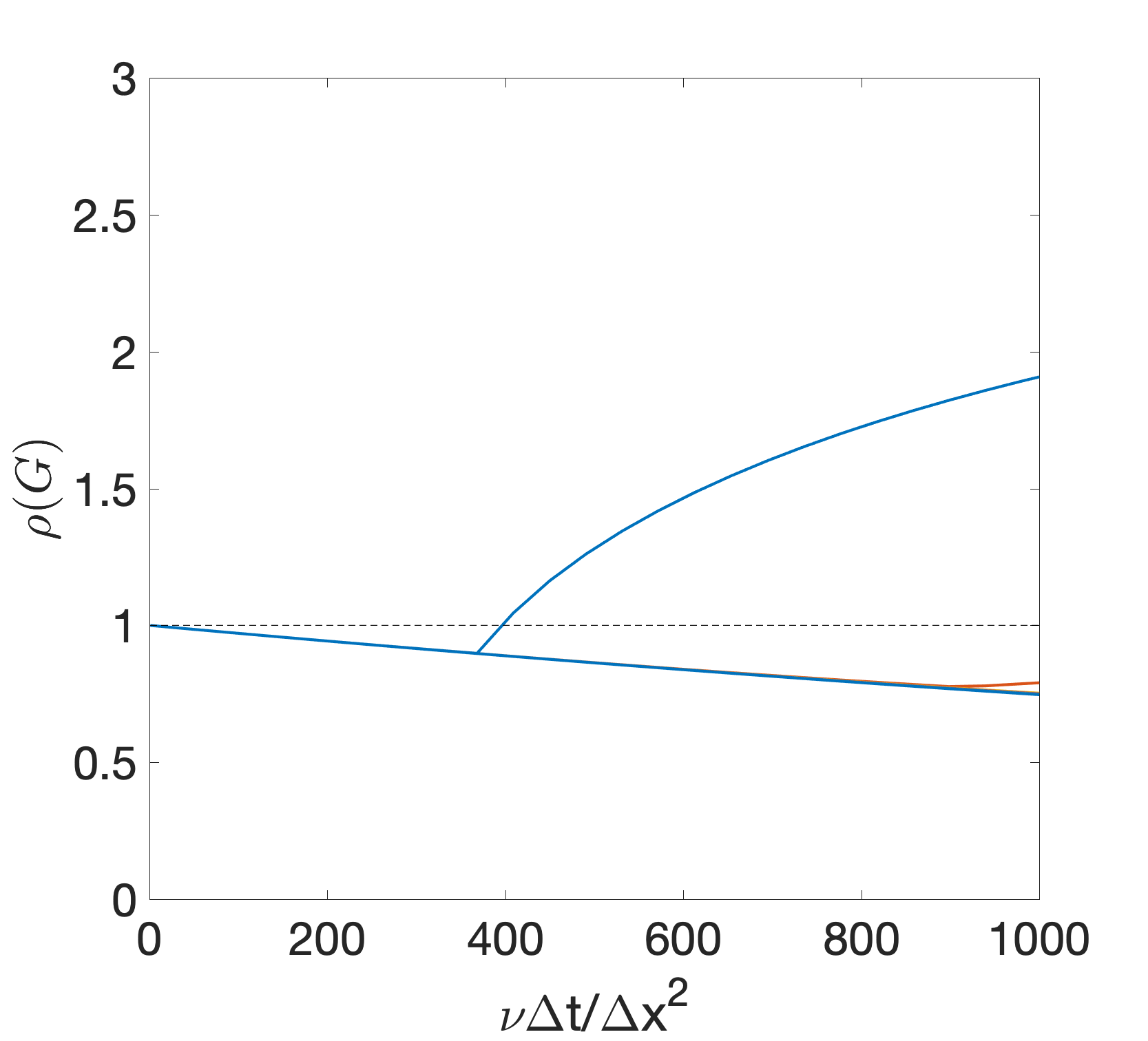} \\
\multicolumn{3}{c}{\textrm{(c) BDF3/EXT3}} \\
\textrm{(i) $N=32,K=5$} & \textrm{(ii) $N=65,K=10$} & \textrm{(iii) $N=98,K=15$} \\
\multicolumn{3}{c}{\includegraphics[width=140mm]{figures/legend}} \\
 \end{array}$
\end{center}
\vspace{-6mm}
\caption{Spectral radius $\rho(G)$ versus nondimensional time $\frac{\nu
\dt}{\dx^2}$ for (a) BDF2/EXT2, (b) BDF3/EXT2, and (c) BDF3/EXT3 scheme
with $N$ and $K$ varying such that $K\dx$ is fixed:
(left to right) (i) $N=32$, $K=5$, (ii) $N=65$, $K=10$, and (iii) $N=98$, $K=15$.}
\label{fig:rhostab_overlap_dx}
\end{figure}

For certain applications (e.g, Fig.  1.1 of \cite{mittaloverlapping}),
geometric constraints can limit the maximum allowable overlap width between
different subdomains.  In these cases, if the overlap width is not enough for a
stable predictor-corrector scheme with $m>1$, the application of the
Schwarz-SEM framework is limited.  Thus, we look at the impact of increasing
the grid resolution, while keeping the overlap width fixed in Fig.
\ref{fig:rhostab_overlap_dx}.  These results show that increasing the grid
resolution has a stabilizing effect on the PC scheme.

\clearpage
\subsection{Difference between stability behavior for odd and even $Q$} \label{sec:oddevendiscussion}
Figure \ref{fig:rhovss} shows that for the high-order BDF$k$/EXT$m$ scheme used
in the Schwarz-SEM framework, odd-$Q$ leads to a relatively more stable
formulation than even-$Q$. Love et al. have
described similar behavior for a second-order predictor-corrector scheme using a
FD-based staggered grid formulation for Lagrangian shock hydrodynamics \cite{love2009stability}.
Similarly, Stetter's stability analysis of a high-order Adam Bashforth- (AB3) and
Adam Moulton-based (AM2) predictor-corrector scheme for ODEs shows a
difference in the stability of odd- and even-$Q$, although this is not a part of
Stetter's discussion.

\begin{figure}[b!]
\begin{center}
$\begin{array}{c}
\includegraphics[height=75mm]{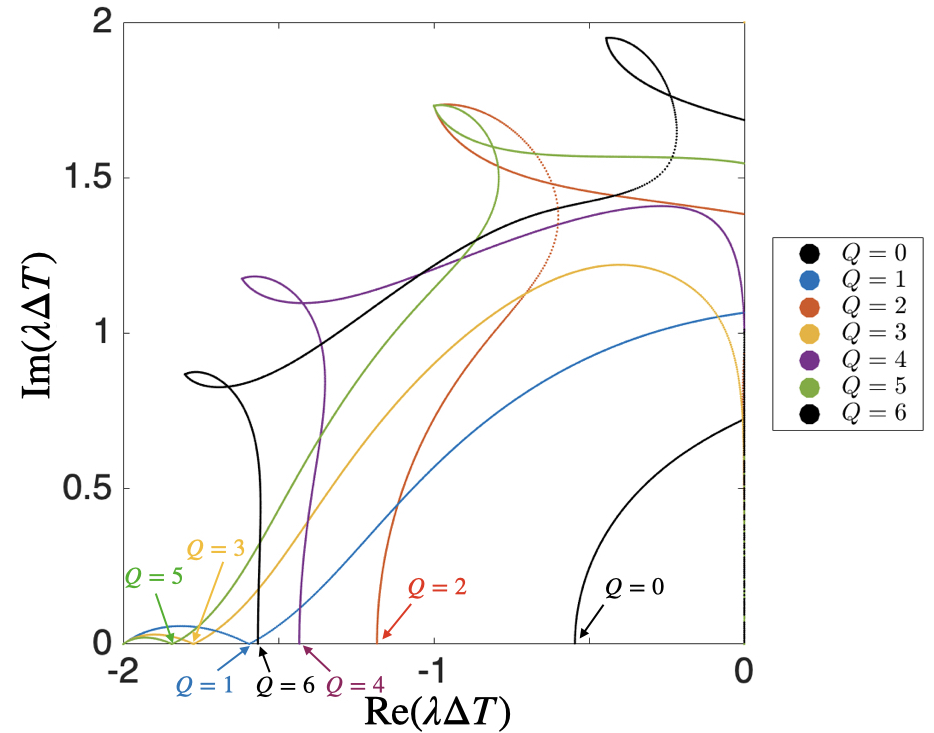}
 \end{array}$
\end{center}
\vspace{-6mm}
\caption{Stability diagram for the AB3/AM2 predictor-corrector scheme for
different number of corrector iterations. The curve corresponding to each
$Q$ is labelled where it intersects with the horizontal-axis to show that
the PC scheme is more stable for odd-$Q$.}
\label{fig:ab3am2stab}
\end{figure}

Figure \ref{fig:ab3am2stab} shows the stability diagram for the AB3/AM2-based
PC scheme used to solve $\frac{du}{dt} = \lambda u$ by Stetter. We observe
that when $\lambda$ is real-valued, the timestepping scheme
is relatively more stable for odd-$Q$ as compared to the even-$Q$.
These stability results for the AB3/AM2-based scheme are similar to the results
of the high-order predictor-corrector scheme in our work because we solve the unsteady heat
equation, where the diffusion operator has negative real eigenvalues. We have also found
that stability analysis of a BDF$k$/EXT$m$-based scheme for $\frac{du}{dt} =
\lambda u$ shows a similar behavior for odd- and even-$Q$.  Finally, this
behavior is also evident in the Schwarz-SEM framework where the PC scheme is
applied to the INSE with the diffusion order replaced by the unsteady Stokes
operator \cite{dfm02}.

\subsubsection{Stability with $Q$ in the Schwarz-SEM framework} \label{sec:oddevenNSE}
The odd-even pattern that we have observed in the Schwarz-SEM framework is
straightforward to demonstrate by considering the exact Navier-Stokes
eigenfunctions by Walsh \cite{walsh1992eddy}, in a periodic domain $\Omega =
[0,2\pi]^2$.

Walsh introduced families of eigenfunctions that can be defined using linear
combinations of $\cos(px) \cos(qy),$ $\sin(px) \cos(qy),$ $\cos(px) \sin(qy),$
and $\sin(px) \sin(qy)$, for all integer pairs ($p, q$) satisfying $\lambda =
-(p^2+q^2)$.  Taking as an initial condition the eigenfunction $\hat{\bu} =
(-{\psi}_y,{\psi}_x)$, a solution to the NSE is $\bu = e^{\nu \lambda t}
\hat{\bu}(\bx)$. Here, ${\psi}$ is the streamfunction resulting from the linear
combinations of eigenfunctions.  Interesting long-time solutions can be
realized by adding a relatively high-speed mean flow $\bu_0$ to the
eigenfunction, in which case the solution is $\bu_{exact} = e^{\nu \lambda t}
\hat{\bu}[\bx-\bu_0t]$, where the brackets imply that the argument is modulo
$2\pi$ in $x$ and $y$.

\begin{figure}[t!]
\begin{center}
$\begin{array}{ccc}
\includegraphics[height=45mm]{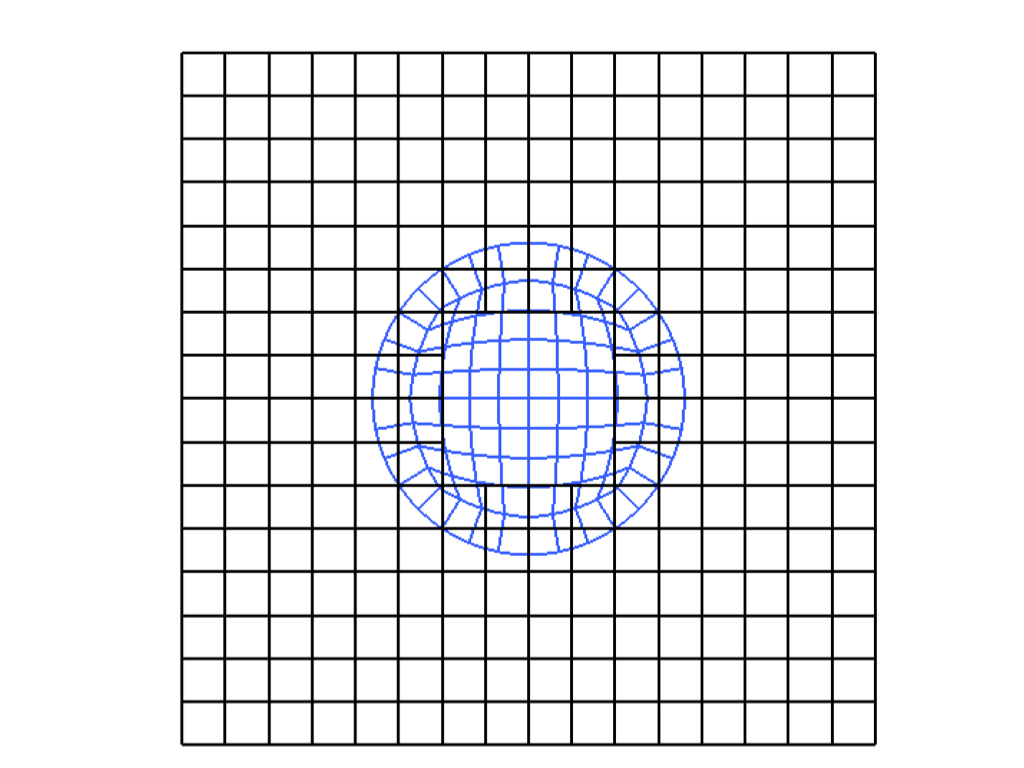} &
\includegraphics[height=45mm]{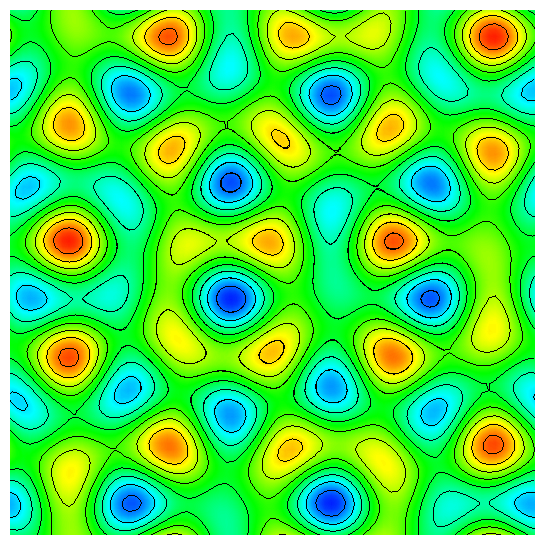} &
\includegraphics[height=45mm]{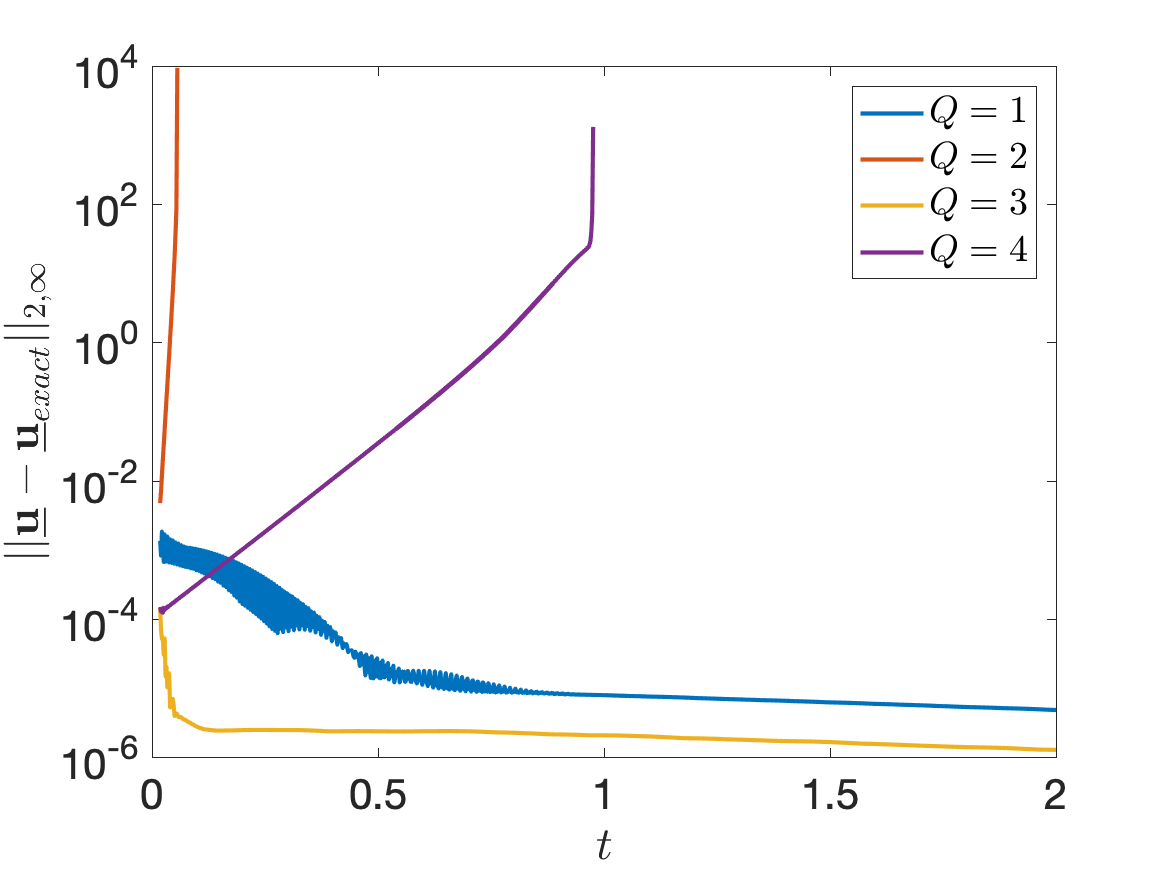}
\end{array}$
\end{center}
\caption{(left to right) (a) Spectral element mesh for each overlapping
subdomain. The background mesh has 240 elements and is covered by a circular
mesh with 96 elements. (b) Vorticity contours at $T_{f}=1$ with $N=7$ and
$Q=3$.  (c) Error variation for the Navier-Stokes eigenfunctions test case with
different $Q$ using $\dt=2\times 10^{-3}$ and $N=7$.} \label{fig:nngrid}
\end{figure}

In the Schwarz-SEM framework, we model this periodic domain using two
overlapping meshes such that the periodic background mesh with a hole in its
center has 240 elements, which is covered with a circular mesh with 96
elements, as shown in Fig. \ref{fig:nngrid}(a).  The flow parameters are
$\nu=0.05$, $\bu_0 = (1,0.3)$, ${\psi} = (1/5)sin(5y) +(1/5)cos(5x) -
(1/4)sin(3x)sin(4y)$, and $\lambda=-25$.  The flow is integrated up to time
$T=1$ convective time units (CTU) with a fixed $\dt=2\times10^{-4}$.  The
polynomial order for representing the solution is set to $N=7$.  Figure
\ref{fig:nngrid}(b) shows the vorticity contours for the solution at $T_f=1$
with $\dt=2\times10^{-4}$.

Since the exact solution is known, we compute the error in time as
$||{\bue}||_{2,\infty}$, where $\bue=\buu-\buu_{exact}$ and the norm is the
2-norm of the point-wise maximum of the vector field.  Figure
\ref{fig:nngrid}(c) shows the error in solution versus time for different $Q$,
and as evident, the solution is stable for $Q=1$ and 3 but unstable for $Q=2$
and 4.  We have also observed similar behavior for odd- and even-$Q$ in highly
turbulent flow calculations using the Schwarz-SEM framework.

Based on the stability analyses presented in this section and the numerical
experiments that we have done with the Schwarz-SEM framework, we conjecture
that for BVP (dominated) with negative real eigenvalues, high-order PC schemes
lead to a difference in the stability between odd- and even-corrector
iterations. In future work, we will continue this work to determine the
fundamental reasoning behind this odd-even pattern, and look at the stability
of high-order PC methods for general boundary value problems, including the
unsteady Stokes problem in our SEM-based framework.

\subsection{Improving the Stability for Even-Corrector Iterations}
\label{sec:stabilityimprovestab}
In the predictor-corrector scheme discussed so far, since odd-$Q$ is more
stable than even-$Q$, we have to increase $Q$ by 2 (e.g., increase $Q=1$ to
$Q=3$) if the number of corrector iterations is not sufficient from a stability
standpoint. Since each corrector iteration requires an additional PDE solve and
thus increases the computational cost the calculation, we modify our PC scheme
to improve the stability for even-$Q$.

In \cite{stetter1968improved}, Stetter improves the stability of the AB3/AM2 PC
scheme by using a linear combination of solution from each corrector iteration
to determine the final solution as
\begin{eqnarray}
\label{eq:ab3am2st}
u^{n} &=& \sum_{q=0}^Q \gamma_q u^{n,q}, \qquad \sum_{q=0}^Q \gamma_q=1,
\end{eqnarray}
where $\gamma_q$ is some weight corresponding to the solution of $q$th iterate
at each timestep.  Using this approach, Stetter shows that the stability
region for the PC scheme can be extended for any $\lambda$.  Since
we are primarily concerned with the inferior stability properties of even-$Q$,
we extend Stetter's idea to modify only the last corrector iteration when $Q$
is even.

The proposed predictor-corrector scheme is
\begin{eqnarray}
\label{uqm1}
q=0: \! \uu_{i}^{n,0} &=& - \sum_{l=1}^{k} \beta_{l} H_i^{-1} \uiiinml
 + \sum_{l=1}^{m} \textw_{l} H_i^{-1} J_{ij} \uu_{j}^{n-l}, \\
\label{uqmodd}
q=1 \dots Q-1: \! \uu_{i}^{n,q} &=& - \sum_{l=1}^{k} \beta_{l} H_i^{-1} \uiiinml
 + H_i^{-1} J_{ij} \uu_{j}^{n,q-1}, \\
\label{uqmeven}
q=Q: \! \uu_{i}^{n,q} &=& - \sum_{l=1}^{k} \beta_{l} H_i^{-1} \uiiinml
 + H_i^{-1} J_{ij} \bigg(\gamma \uu_{j}^{n,Q-1} + (1-\gamma) \uu_{j}^{n,Q-2}\bigg),
\end{eqnarray}
where $\gamma$ is a parameter that uses the combination of the two most recent
solutions (instead of just the most recent solution) at the last corrector
iteration.  We set $\gamma=1$ when $Q$ is odd to recover the original PC scheme
(\ref{eq:uq1}-\ref{eq:uqall}), and $0<\gamma<1$ when $Q$ is even.  The
rationale behind this new predictor-corrector scheme is that we do not want to
modify the convergence properties of the original PC scheme if $Q$ is odd.
Thus, we modify only the last corrector iteration ($u^{n,Q}$), when $Q$ is
even.

\begin{figure}[t!]
\begin{center}
$\begin{array}{c}
\includegraphics[height=60mm]{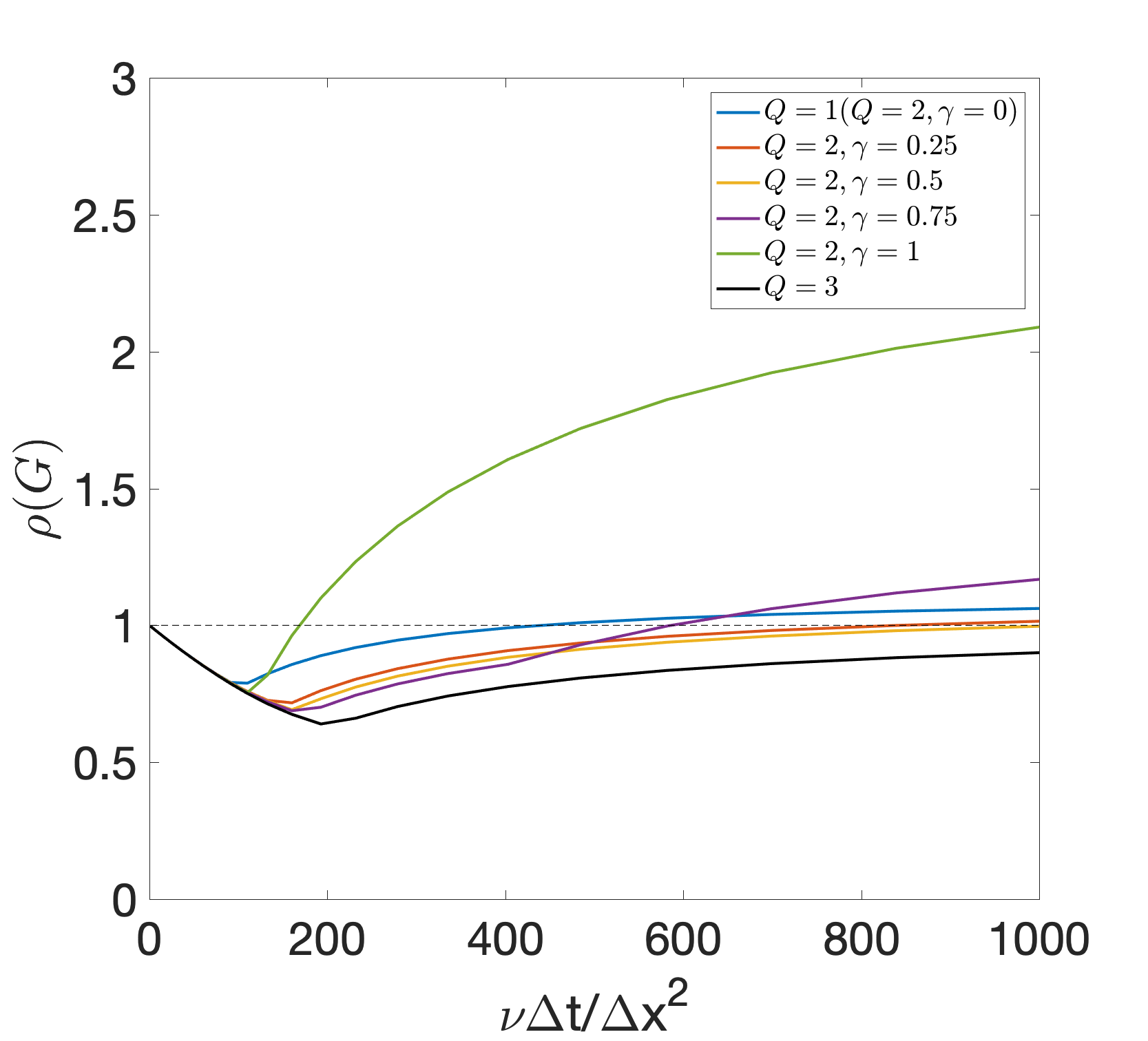}
 \end{array}$
\end{center}
\vspace{-6mm}
\caption{Spectral radius $\rho(G)$ versus nondimensional time $\frac{\nu
\dt}{\dx^2}$ for the BDF3/EXT3 scheme with $N=32$, $K=5$ for different $\gamma$.}
\label{fig:pcgamma}
\end{figure}

\begin{figure}[tb!]
\begin{center}
$\begin{array}{cc}
\includegraphics[height=60mm]{figures/rho_vs_s_n132_k5_BDF3_EXT3} &
\includegraphics[height=60mm]{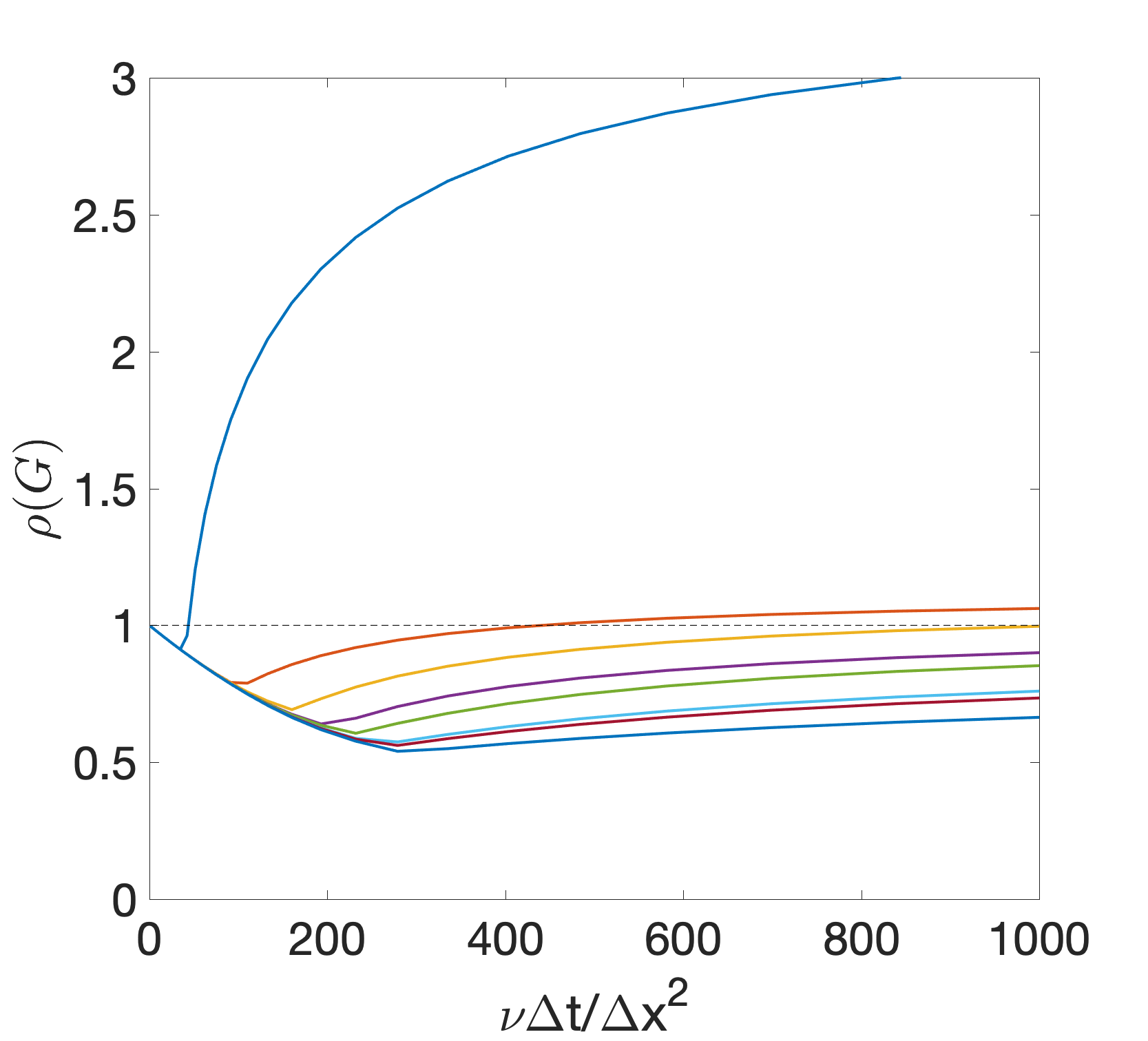} \\
\multicolumn{2}{c}{\includegraphics[width=150mm]{figures/legend}}  \\
\textrm{(a) Original PC scheme} & \textrm{(b) Improved PC scheme}
 \end{array}$
\end{center}
\vspace{-6mm}
\caption{Spectral radius $\rho(G)$ versus nondimensional time $\frac{\nu
\dt}{\dx^2}$ for the BDF3/EXT3 scheme with $N=32$ and $K=5$ comparing the original and improved predictor-corrector scheme. $\gamma=0.5$ for even-$Q$ in the improved PC scheme.}
\label{fig:pcgammak5}
\end{figure}

Figure \ref{fig:pcgamma} shows the stability plot for the modified
predictor-corrector schemes for different values of $\gamma$ for $Q=2$, and
$Q=3$ with $\gamma=1$. The parameters for grid size are $N=32$ and $K=5$ and we
use the BDF3/EXT3 scheme for time-integration.  We see that $Q=2, \gamma=1$
leads to the original scheme and $Q=2, \gamma=0$ restores the behavior of the
original scheme with $Q=1$.  It is also apparent that the stability of the PC
scheme has significantly improved for $Q=2$, in comparison to the original
scheme (Fig. \ref{fig:rhovss}), and it is no longer more unstable than $Q=1$
when $\gamma=0.25$ or 0.5.  Figure \ref{fig:pcgammak5} compares the stability
plots for the original and improved predictor-corrector scheme ($\gamma=0.5$
for even-$Q$), and we see that the proposed formulation leads to a scheme with
monotonically increasing stability with $Q$.

\begin{figure}[tb!]
\includegraphics[width=90mm]{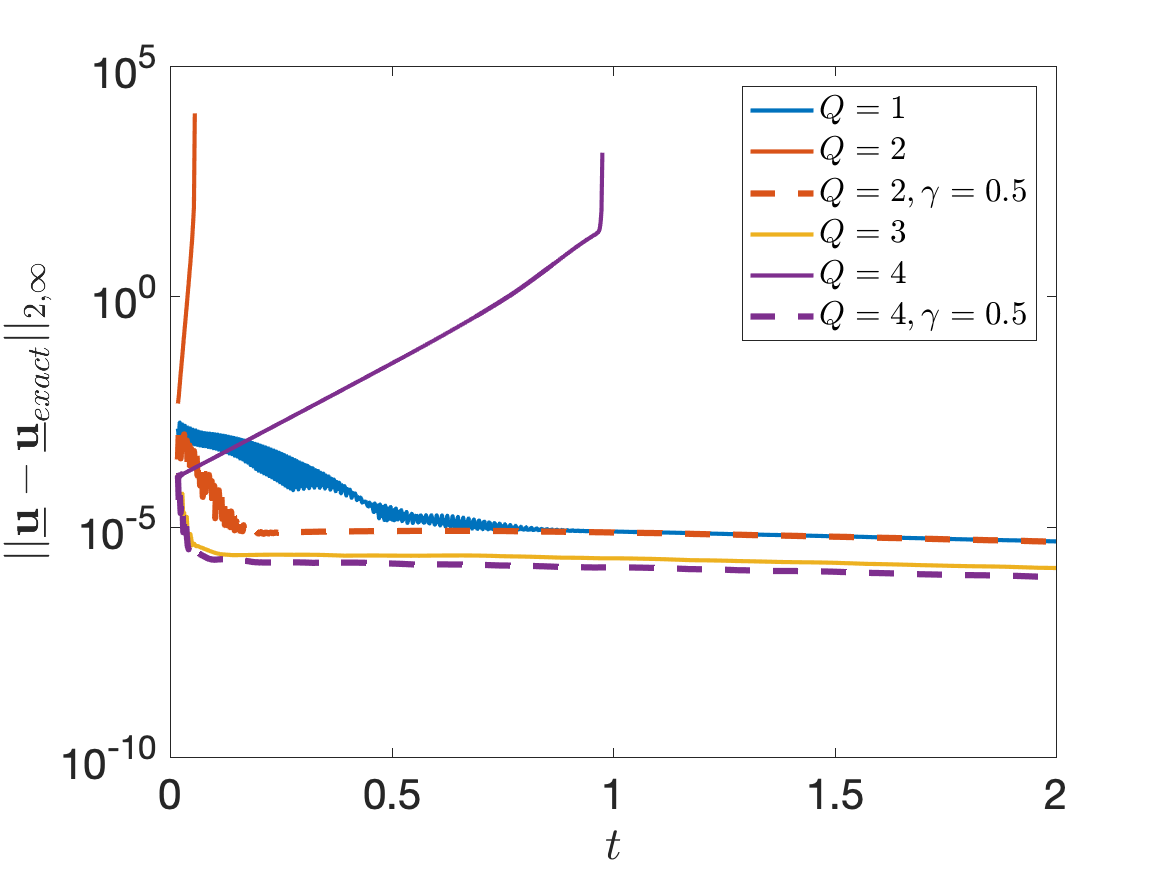}
\vspace{-4mm}
\caption{Error variation for the Navier-Stokes eigenfunctions test case from
Section \ref{sec:oddevenNSE} with different $Q$ using $\dt=2\times10^{-3}$ and
$N=7$ comparing the original (\solidrule) and improved (\protect\dashedrule)
predictor-corrector scheme.}
\label{fig:eddyoddevenfixed} \end{figure}

The improved predictor-corrector scheme for even-$Q$ can be readily extended to
the Schwarz-SEM framework. Using the exact Navier-Stokes eigenfunctions, Fig.
\ref{fig:eddyoddevenfixed} shows the error variation for $Q=1, 2,$ 3 and 4 with
the original scheme, and compares it to the improved PC scheme with $Q=2$ and 4
for $\gamma=0.5$.  In \cite{mittaloverlapping}, we have also verified that the
improved PC scheme maintains $m$th-order temporal accuracy of the underlying
SEM-based solver.

\section{Stability of multirate PC method} \label{sec:multistability}
Next, we analyze the multirate PC scheme described in Section
\ref{sec:multirate} using the matrix method for stability analysis. For a
timestep ratio $\tsr$, advancing the solution from $t^{n-1}$ to $t^n$ requires
the PDE of interest to be solved $\tsr$ times in $\Omega^f$ and once in
$\Omega^c$.  Thus, casting the multirate timestepping scheme into a system
$\uz^n=G\uz^{n-1}$ is not as straightforward as it is for the singlerate
scheme. For simplicity, we describe the method to cast the multirate PC scheme
(\ref{eq:uqmspc})-(\ref{eq:uqmscf}) for $\tsr=2$, and this method readily
extends to arbitrary $\tsr$.

For notational purposes, we define
\begin{eqnarray}
\label{eq:uzmulti}
\uz^{n} = \big[\uu^{c,n^T} \, \uu^{f,n^T}
\, \uu^{f,n-1/2^T} \, \uu^{c,n-1^T} \, \uu^{f,n-1^T} \, \uu^{f,n-3/2^T} \,
\uu^{c,n-2^T} \, \uu^{f,n-2^T} \, \uu^{f,n-5/2^T} \, \uu^{c,n-3^T} \,
\uu^{f,n-3^T} \big]^T,
\end{eqnarray}
where $\uu^f$ and $\uu^c$ are the solution vectors for subdomain $\Omega^f$ and
$\Omega^c$, respectively, and we use $\uz^{n,q}$ to represent the vector with
solutions at $q$th corrector iteration.  To account for the sub-timestep
solution for $\Omega^f$ in  $\uz^{n}$, we modify our methodology of building
the predictor and corrector matrices. The predictor matrix is now a product of
$\tsr$ matrices, of which $\tsr-1$ matrices correspond to the sub-timesteps of
$\Omega^f$, and 1 matrix for the last sub-timestep of $\Omega^f$ and the only
step of $\Omega^c$.  For $\tsr=2$, the predictor matrix is $P=P_2P_1$ where
$P_1$ outputs the solution $\uu^{f,n-\frac{1}{2},0}$ and $P_2$ outputs the
solution $\uu^{c,n,0}$ and $\uu^{f,n,0}$. The matrices $P_2$ and $P_1$ are

\scriptsize
\hspace{-5mm}
$
\left[ \begin{array}{c} u^{f,n-\frac{1}{2},0} \\ u^{c,n-1,Q} \\
u^{f,n-1,Q} \\ u^{f,n-\frac{3}{2},Q} \\ u^{c,n-2,0} \\ u^{f,n-2,Q} \\
u^{f,n-\frac{5}{2},Q} \\ u^{c,n-3,Q} \\ u^{f,n-3,0} \end{array}
\right]=
$
$
\underbrace{
\begin{pmatrix}
\textw_{11} H_{f}^{-1} J_{fc} & -\beta_1 H_{f}^{-1} & -\beta_2 H_{f}^{-1} & \textw_{12} H_{f}^{-1} J_{fc} &
-\beta_3 H_{f}^{-1} & 0 & \textw_{13} H_{f}^{-1} J_{fc} & 0 & 0 & 0 & 0 \\
I_1 & 0 & 0 & 0 & 0 & 0 & 0 & 0 & 0 & 0 & 0 \\
0 & I_2 & 0 & 0 & 0 & 0 & 0 & 0 & 0 & 0 & 0 \\
0 & 0 & I_2 & 0 & 0 & 0 & 0 & 0 & 0 & 0 & 0 \\
0 & 0 & 0 & I_1 & 0 & 0 & 0 & 0 & 0 & 0 & 0 \\
0 & 0 & 0 & 0 & I_2 & 0 & 0 & 0 & 0 & 0 & 0 \\
0 & 0 & 0 & 0 & 0 & I_2 & 0 & 0 & 0 & 0 & 0 \\
0 & 0 & 0 & 0 & 0 & 0 & I_1 & 0 & 0 & 0 & 0 \\
0 & 0 & 0 & 0 & 0 & 0 & 0 & I_2 & 0 & 0 & 0 \\
\end{pmatrix}
}_{P_1}
$
$
\underbrace{\left[ \begin{array}{c} u^{c,n-1,Q} \\ u^{f,n-1,Q} \\
u^{f,n-\frac{3}{2},Q} \\ u^{c,n-2,Q} \\ u^{f,n-2,Q} \\
u^{f,n-\frac{5}{2},Q} \\ u^{c,n-3,Q} \\ u^{f,n-3,Q} \\
u^{f,n-\frac{7}{2},Q} \\ u^{c,n-4,Q} \\ u^{f,n-4,Q} \end{array} \right]}_{\uz^{n,Q}}
$,\\
\\
\\
\hspace{-5mm}
$
\underbrace{\left[ \begin{array}{c} u^{c,n,0} \\ u^{f,n,0} \\
u^{f,n-\frac{1}{2},0} \\ u^{c,n-1,Q} \\ u^{f,n-1,Q} \\
u^{f,n-\frac{3}{2},Q} \\ u^{c,n-2,Q} \\ u^{f,n-2,Q} \\
u^{f,n-\frac{5}{2},Q} \\ u^{c,n-3,Q} \\ u^{f,n-3,Q} \end{array} \right]}_{\uz^{n,0}} =
$
$
\underbrace{
\begin{bmatrix}
0 & -\beta_1 H_{c}^{-1} & \textw_{21} H_{c}^{-1} J_{cf} & \textw_{22} H_{c}^{-1} J_{cf} &
-\beta_2 H_{c}^{-1} & \textw_{23} H_{c}^{-1} J_{cf} & 0 & -\beta_3 H_{c}^{-1} & 0 \\
-\beta_1 H_{f}^{-1} & \textw_1 H_{f}^{-1} J_{fc} & -\beta_2 H_{f}^{-1} & -\beta_3 H_{f}^{-1} &
\textw_2 H_{f}^{-1} J_{fc} & 0 & 0 & \textw_3 H_{f}^{-1} J_{fc} & 0 \\
I_2 & 0 & 0 & 0 & 0 & 0 & 0 & 0 & 0 \\
0 & I_1 & 0 & 0 & 0 & 0 & 0 & 0 & 0 \\
0 & 0 & I_2 & 0 & 0 & 0 & 0 & 0 & 0 \\
0 & 0 & 0 & I_2 & 0 & 0 & 0 & 0 & 0 \\
0 & 0 & 0 & 0 & I_1 & 0 & 0 & 0 & 0 \\
0 & 0 & 0 & 0 & 0 & I_2 & 0 & 0 & 0 \\
0 & 0 & 0 & 0 & 0 & 0 & I_2 & 0 & 0 \\
0 & 0 & 0 & 0 & 0 & 0 & 0 & I_1 & 0 \\
0 & 0 & 0 & 0 & 0 & 0 & 0 & 0 & I_2 \\
\end{bmatrix}
}_{P_2}
$
$
\left[ \begin{array}{c} u^{f,n-\frac{1}{2},0} \\ u^{c,n-1,Q} \\ u^{f,n-1,Q} \\
u^{f,n-\frac{3}{2},Q} \\ u^{c,n-2,Q} \\ u^{f,n-2,Q} \\
u^{f,n-\frac{5}{2},Q} \\ u^{c,n-3,Q} \\ u^{f,n-3,Q} \end{array} \right]
$.
\normalsize
Thus, the predictor step to time-advance the solution in $\Omega^f$ and
$\Omega^c$ from $t^{n-1}$ to $t^n$ is $\uz^{n,0} = P_2 P_1 \uz^{n,q}$.
Similarly, the system for corrector iterations is\\

\scriptsize
\hspace{-5mm}
$
\left[ \begin{array}{c} u^{c,n,q-1} \\ u^{f,n,q-1} \\
 u^{f,n-\frac{1}{2},q} \\ u^{c,n-1,Q} \\ u^{f,n-1,Q} \\
u^{f,n-\frac{3}{2},Q} \\ u^{c,n-2,Q} \\ u^{f,n-2,Q} \\
u^{f,n-\frac{5}{2},Q} \\ u^{c,n-3,Q} \\ u^{f,n-3,Q} \end{array} \right]=
$
$
\underbrace{
\begin{bmatrix}
I_1 & 0 & 0 & 0 & 0 & 0 & 0 & 0 & 0 & 0 & 0 & \\
0 & I_2 & 0 & 0 & 0 & 0 & 0 & 0 & 0 & 0 & 0 & \\
\tintw_{11} H_{f}^{-1} J_{fc} & 0 & 0 & \tintw_{12} H_{f}^{-1} J_{fc} & -\beta_1 H_{f}^{-1} &
-\beta_2 H_{f}^{-1} & \tintw_{13} H_{f}^{-1} J_{fc} & -\beta_3 H_{f}^{-1} & 0 & 0 & 0 \\
0 & 0 & 0 & I_1 & 0 & 0 & 0 & 0 & 0 & 0 & 0 \\
0 & 0 & 0 & 0 & I_2 & 0 & 0 & 0 & 0 & 0 & 0 \\
0 & 0 & 0 & 0 & 0 & I_2 & 0 & 0 & 0 & 0 & 0 \\
0 & 0 & 0 & 0 & 0 & 0 & I_1 & 0 & 0 & 0 & 0 \\
0 & 0 & 0 & 0 & 0 & 0 & 0 & I_2 & 0 & 0 & 0 \\
0 & 0 & 0 & 0 & 0 & 0 & 0 & 0 & I_2 & 0 & 0 \\
0 & 0 & 0 & 0 & 0 & 0 & 0 & 0 & 0 & I_1 & 0 \\
0 & 0 & 0 & 0 & 0 & 0 & 0 & 0 & 0 & 0 & I_2 \\
\end{bmatrix}
}_{C_1}
$
$
\underbrace{\left[ \begin{array}{c} u^{c,n,q-1} \\ u^{f,n,q-1} \\
u^{f,n-\frac{1}{2},q-1} \\ u^{c,n-1,Q} \\ u^{f,n-1,Q} \\
u^{f,n-\frac{3}{2},Q} \\ u^{c,n-2,Q} \\ u^{f,n-2,Q} \\
u^{f,n-\frac{5}{2},Q} \\ u^{c,n-3,Q} \\ u^{f,n-3,Q} \end{array} \right]}_{\uz^{n,q-1}}
$,\\
\\
\\
$
\underbrace{\left[\begin{array}{c} u^{c,n,q} \\ u^{f,n,q} \\
 u^{f,n-\frac{1}{2},q} \\ u^{c,n-1,Q} \\ u^{f,n-1,Q} \\
u^{f,n-\frac{3}{2},Q} \\ u^{c,n-2,Q} \\ u^{f,n-2,Q} \\
u^{f,n-\frac{5}{2},Q} \\ u^{c,n-3,Q} \\ u^{f,n-3,Q} \end{array} \right]}_{\uz^{n,q}}=
$
$
\underbrace{
\begin{bmatrix}
0 & H_{c}^{-1} J_{cf} & 0 & -\beta_1 H_{c}^{-1} & 0 & 0 & -\beta_2 H_{c}^{-1} &
0 & 0 & -\beta_3 H_{c}^{-1} & 0 \\
H_{f}^{-1} J_{fc} & 0 & -\beta_1 H_{f}^{-1} & 0 & -\beta_2 H_{f}^{-1} &
-\beta_3 H_{f}^{-1} & 0 & 0 & 0 & 0 & 0  \\
0 & 0 & I_2 & 0 & 0 & 0 & 0 & 0 & 0 & 0 & 0 \\
0 & 0 & 0 & I_1 & 0 & 0 & 0 & 0 & 0 & 0 & 0 \\
0 & 0 & 0 & 0 & I_2 & 0 & 0 & 0 & 0 & 0 & 0 \\
0 & 0 & 0 & 0 & 0 & I_2 & 0 & 0 & 0 & 0 & 0 \\
0 & 0 & 0 & 0 & 0 & 0 & I_1 & 0 & 0 & 0 & 0 \\
0 & 0 & 0 & 0 & 0 & 0 & 0 & I_2 & 0 & 0 & 0 \\
0 & 0 & 0 & 0 & 0 & 0 & 0 & 0 & I_2 & 0 & 0 \\
0 & 0 & 0 & 0 & 0 & 0 & 0 & 0 & 0 & I_1 & 0 \\
0 & 0 & 0 & 0 & 0 & 0 & 0 & 0 & 0 & 0 & I_2 \\
\end{bmatrix}
}_{C_2}
$
$
\left[ \begin{array}{c} u^{c,n,q-1} \\ u^{f,n,q-1} \\
 u^{f,n-\frac{1}{2},q} \\ u^{c,n-1,Q} \\ u^{f,n-1,Q} \\
u^{f,n-\frac{3}{2},Q} \\ u^{c,n-2,Q} \\ u^{f,n-2,Q} \\
u^{f,n-\frac{5}{2},Q} \\ u^{c,n-3,Q} \\ u^{f,n-3,Q} \end{array} \right]
$,
\normalsize
and thus, the corrector step is $\uz^{n,q} = C_2 C_1 \uz^{n,q-1}$. We note that
the solution $u^{f,n-1/2,q}$ is effected via $C_1$, and $u^{f,n,q}$ and
$u^{c,n,q}$ are determined using $C_2$. For $Q$ corrector iterations, thus, the
multirate PC scheme is $\uz^{n}=G\uz^{n-1}$, where $G=C^QP$, $C=C_2C_1$, and
$P=P_2P_1$, for $\tsr=2$.

This methodology can readily be extended for arbitrary $\tsr$ where the
predictor matrix is $P=P_\tsr \dots P_1$ and the corrector matrix is $C=C_\tsr
\dots C_1$. For example, for $\tsr=3$, the predictor matrix $P$ is
$P=P_3P_2P_1$, where $P_1$ determines $u^{f,n-1/3,0}$, $P_2$ determines
$u^{f,n-2/3,0}$, and $P_3$ determines $u^{f,n,0}$ and $u^{c,n,0}$.  Similarly
the corrector matrix $C$ is $C=C_3C_2C_1$ where $C_1$ determines
$u^{f,n-1/3,q}$, $C_2$ determines $u^{f,n-2/3,q}$, and $C_3$ determines
$u^{f,n,q}$ and $u^{c,n,q}$.

Using this approach, we determine the growth matrix $G=C^QP$ for any arbitrary
$\tsr$. The spectral radius of $G$ is used to understand the stability behavior
of the multirate timestepping scheme for different parameters such as the
extrapolation order of interdomain boundary data during the predictor step
($m$), grid resolution ($N$), overlap width ($K$), and the number of corrector
iterations ($Q$).

\subsection{Stability results for different BDF$k$/EXT$m$ schemes}
In this section, we present the spectral radius ($\rho(G)$) versus
nondimensional time plots ($\nu\dt_c/\dx^2$) for different timestep size
ratios. We start with $\tsr=2$, and then look at increasing values of $\tsr$.

\begin{figure}[t!]
\begin{center}
$\begin{array}{ccc}
\includegraphics[height=45mm]{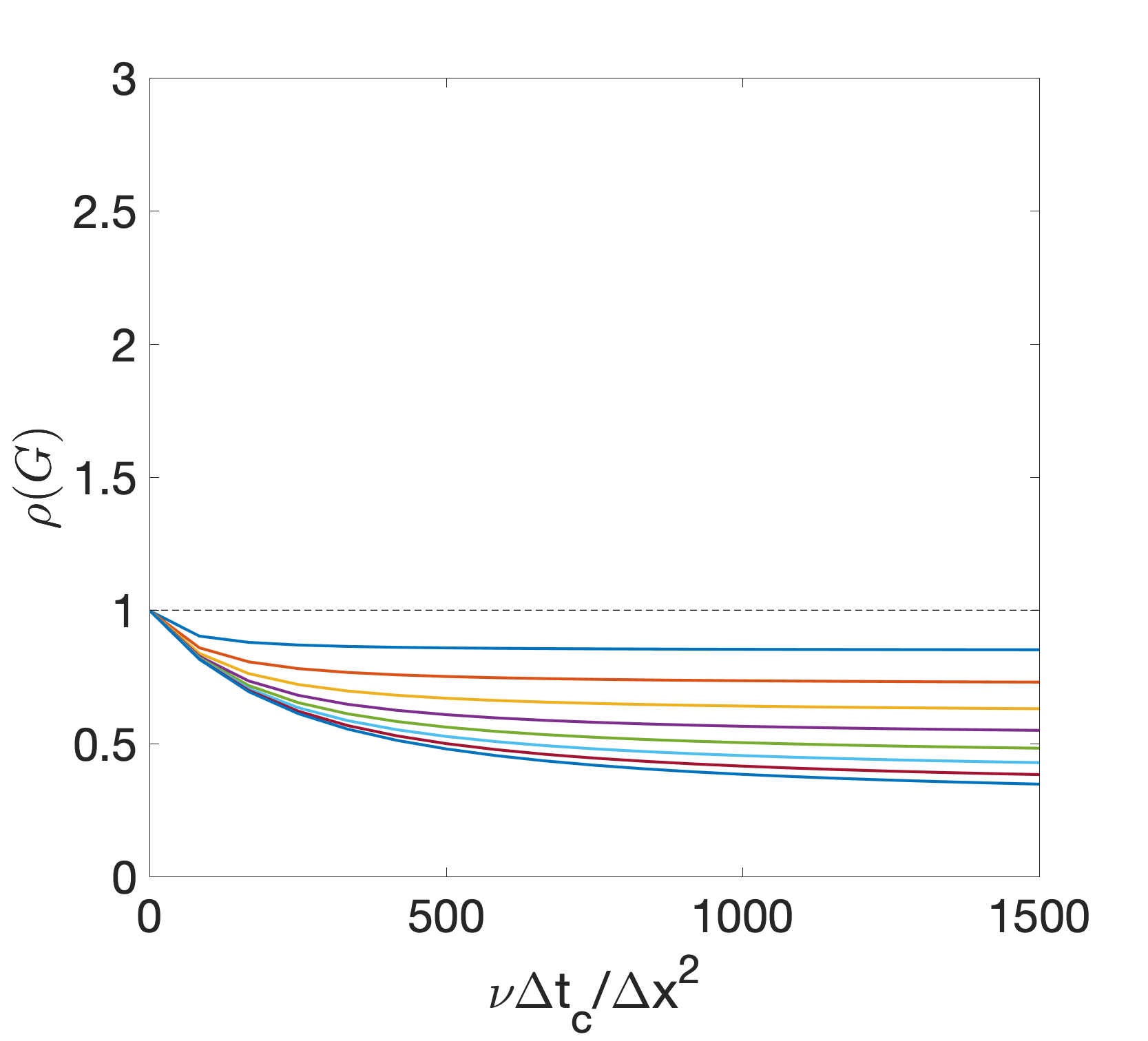} &
\includegraphics[height=45mm]{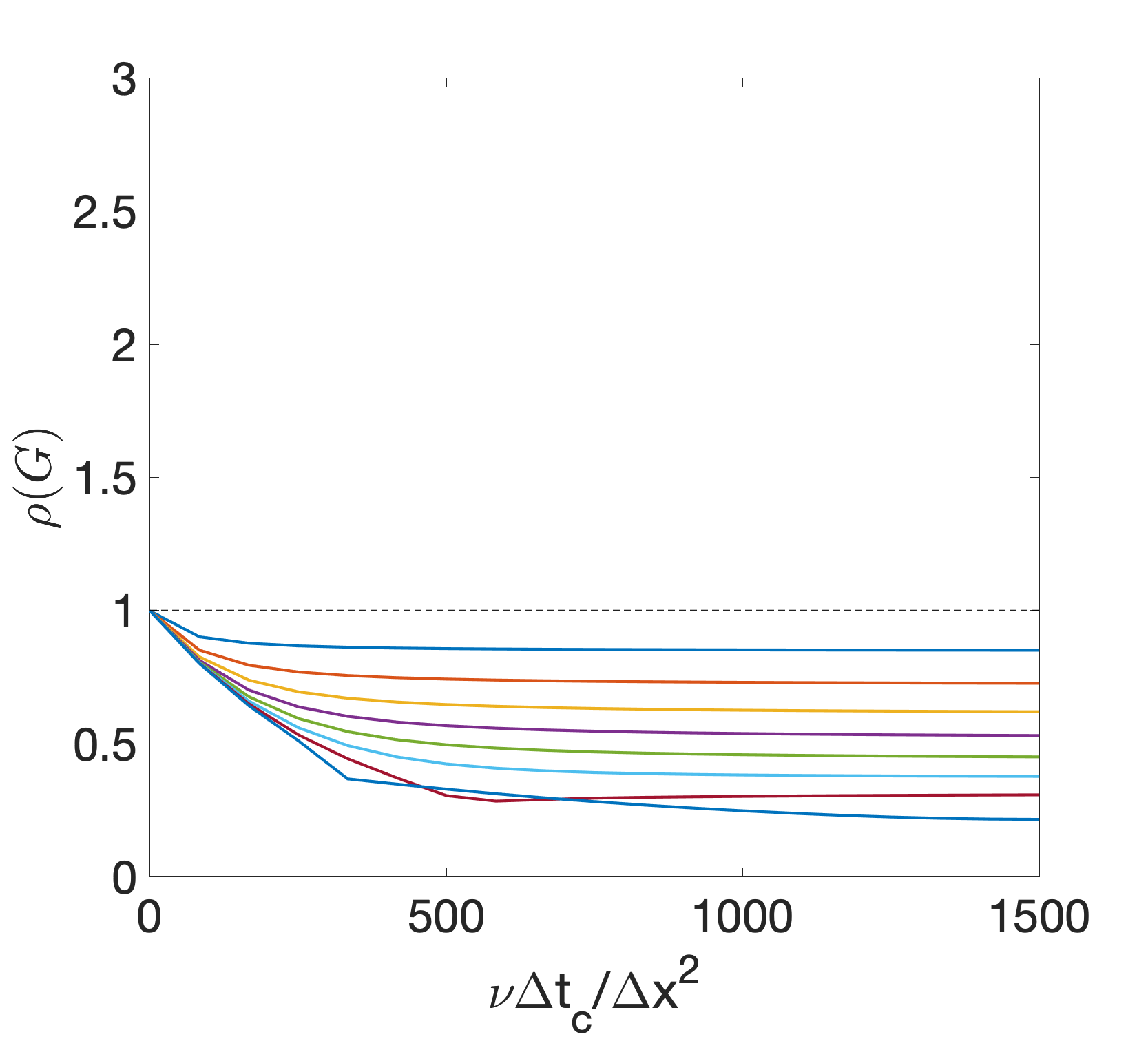} &
\includegraphics[height=45mm]{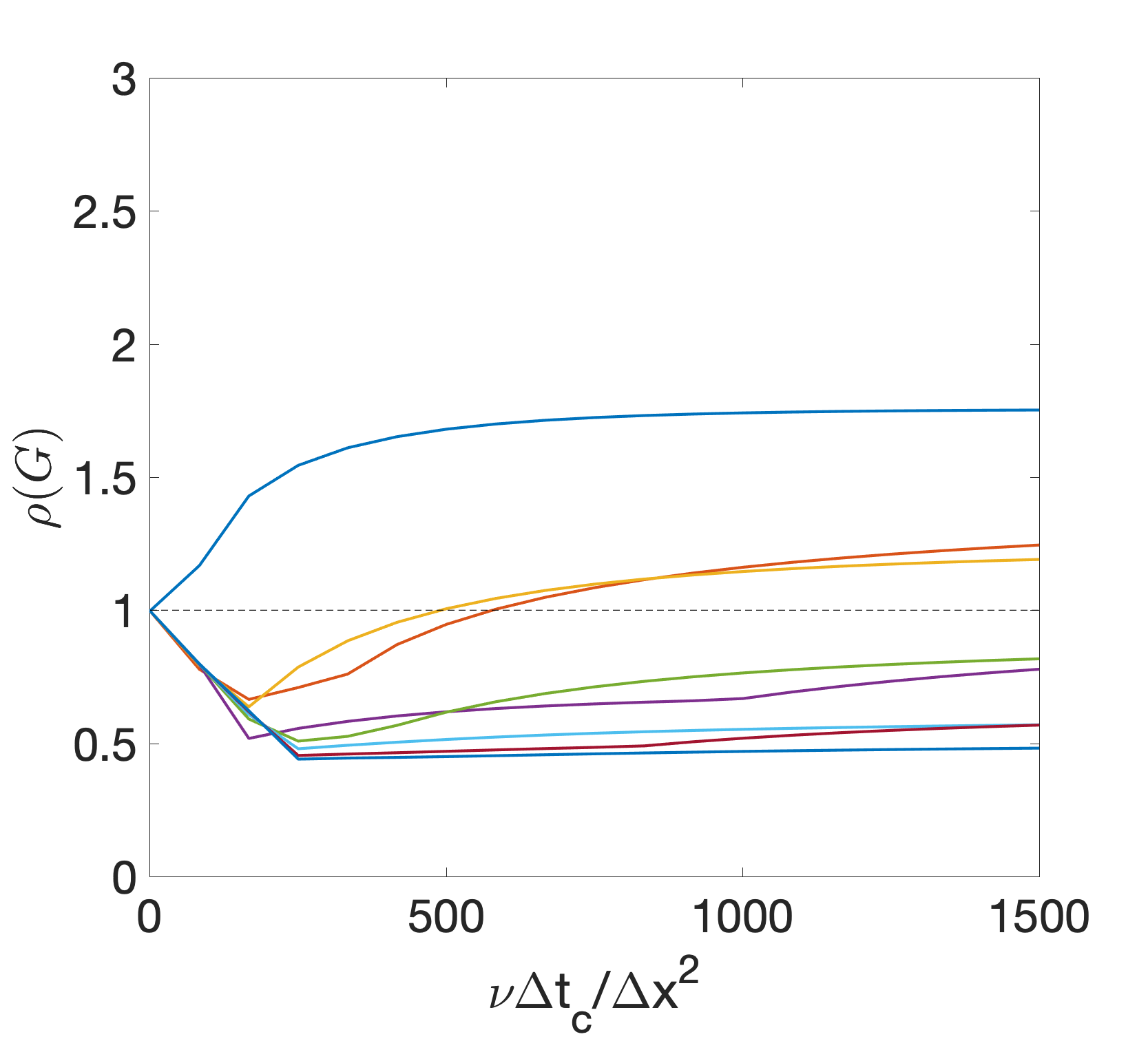} \\
\textrm{(a) BDF$1$/EXT$1$} &
\textrm{(b) BDF$2$/EXT$1$} &
\textrm{(c) BDF$2$/EXT$2$} \\
\includegraphics[height=45mm]{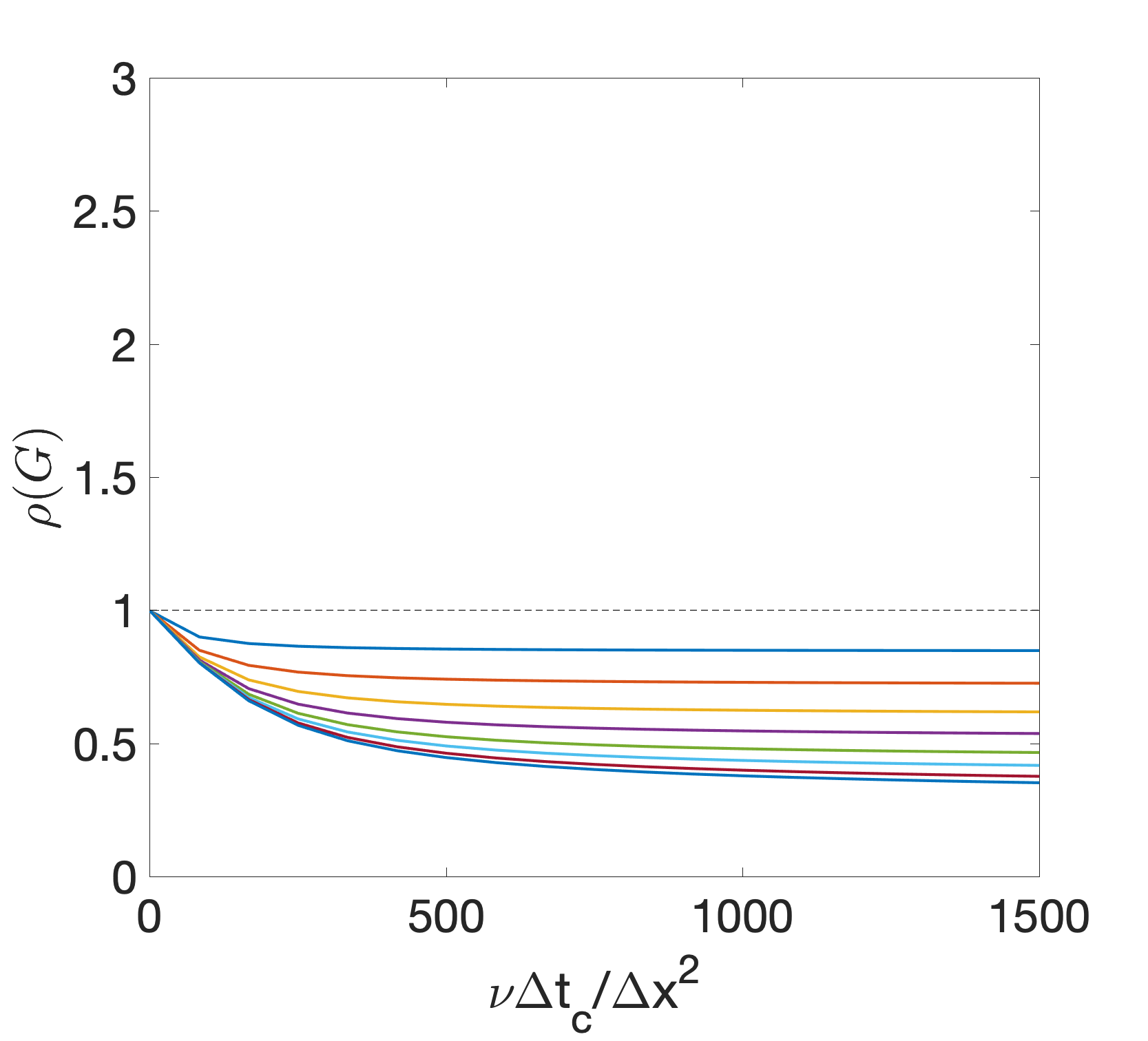} &
\includegraphics[height=45mm]{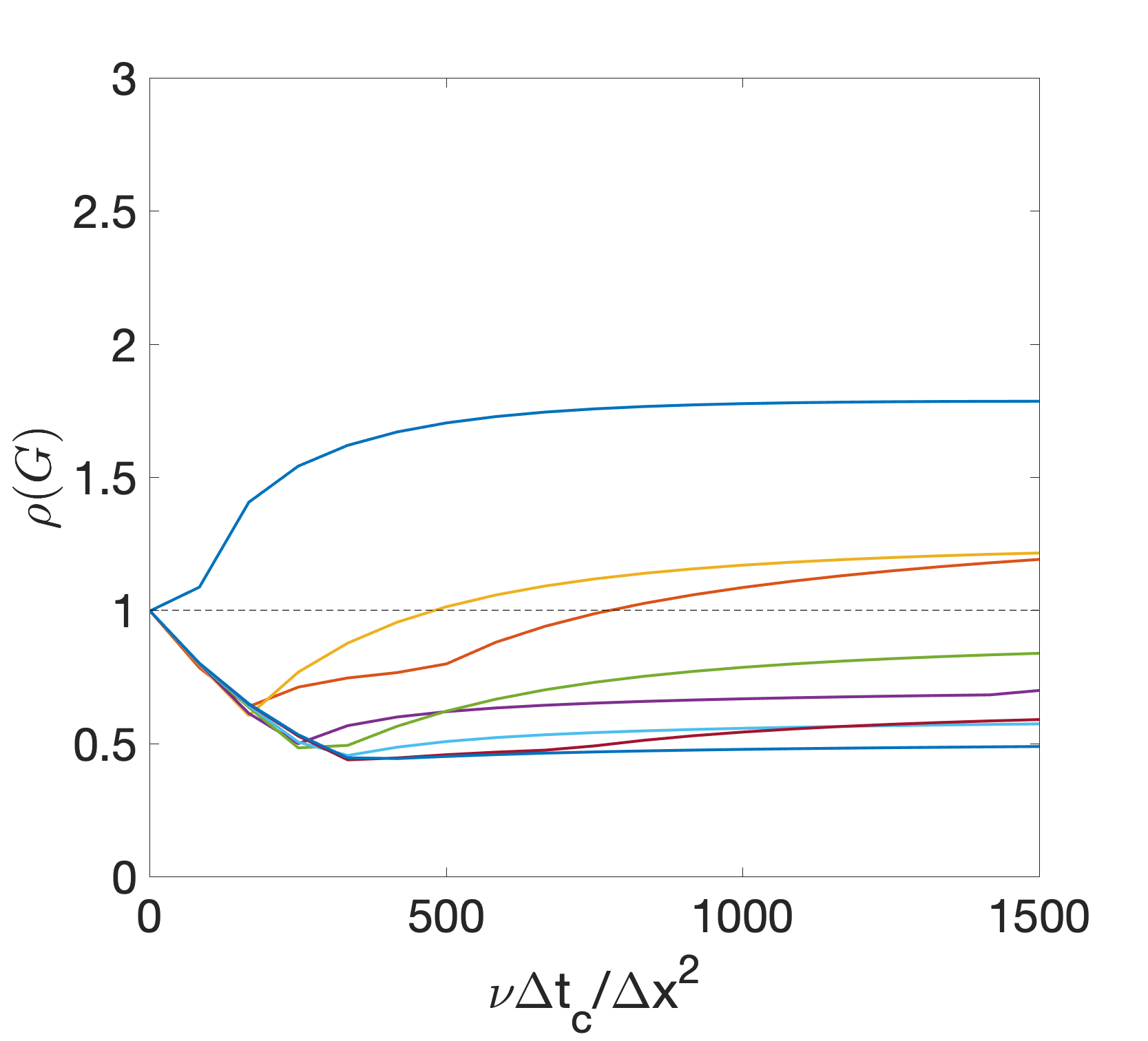} &
\includegraphics[height=45mm]{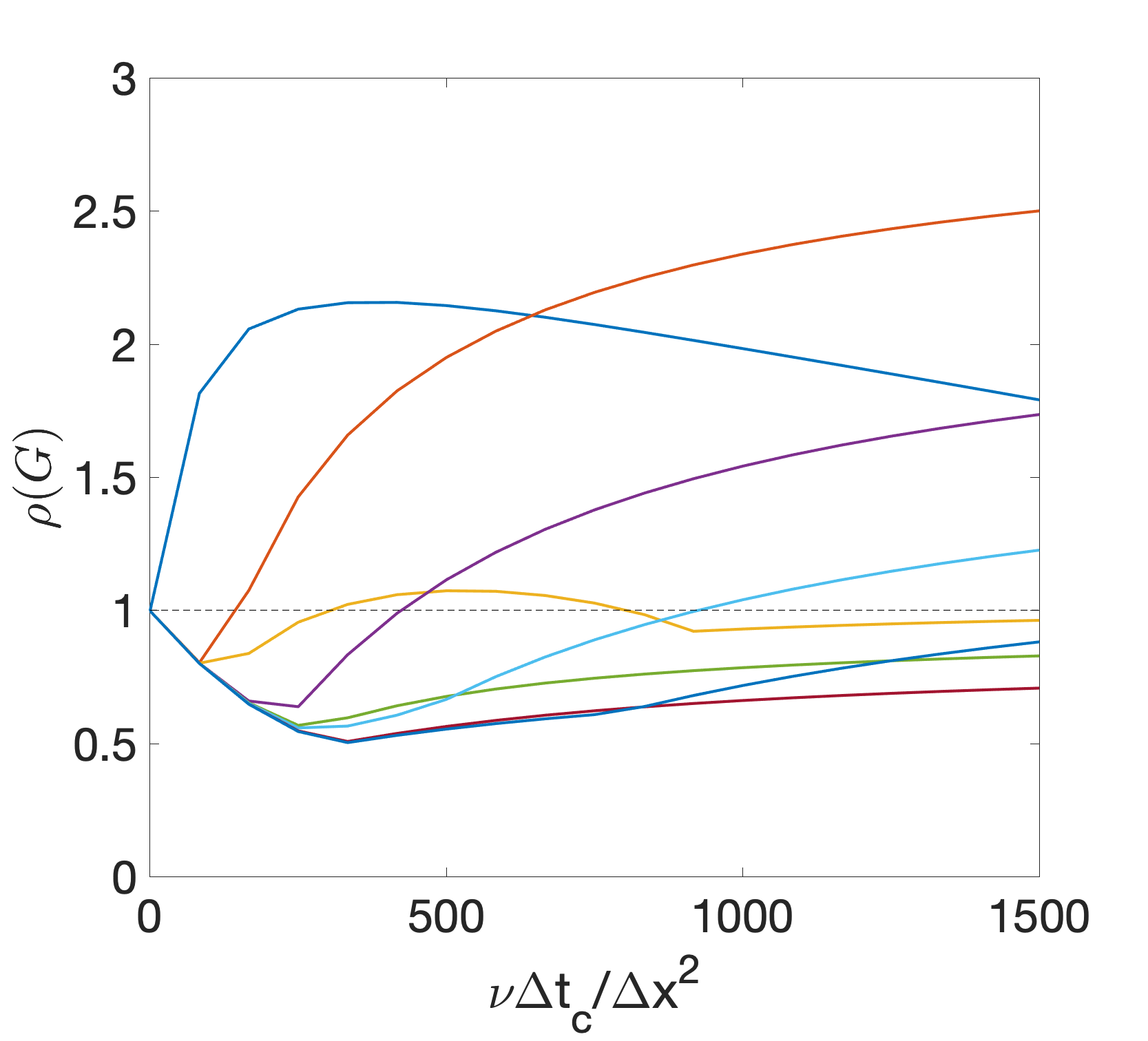} \\
\textrm{(d) BDF$3$/EXT$1$} &
\textrm{(e) BDF$3$/EXT$2$} &
\textrm{(f) BDF$3$/EXT$3$} \\
\multicolumn{3}{c}{\includegraphics[width=150mm]{figures/legend}} \\
 \end{array}$
\end{center}
\vspace{-6mm}
\caption{Spectral radius $\rho(G)$ versus nondimensional time $\frac{\nu
\dt_c}{\dx^2}$ for different BDF$k$/EXT$m$ schemes for $\tsr=2$ with $Q=0\dots7$, with grid parameters $N=32$ and $K=5$.}
\label{fig:rho_multi_eta2}
\end{figure}

Figure \ref{fig:rho_multi_eta2} shows the stability plots for BDF$k$/EXT$m$
scheme for different number of corrector iterations ($Q$) with $\tsr=2$ and
grid parameters $N=32$ and $K=5$.  We conclude from Fig.
\ref{fig:rho_multi_eta2} that using high-order extrapolation for interdomain
boundary data ($m$) decreases the nondimensional time at which the spectral
radius of $G$ is greater than 1. This result is similar to the stability
results for the singlerate timestepping scheme (Fig.  \ref{fig:rhovss}). We
also see that using an {odd}-$Q$ decreases the stability of the scheme, which
is in contrast to the decreased stability for even-$Q$ with the singlerate
timestepping method ($\tsr=1$). To determine if the behavior of odd- and
even-$Q$ is universal for all $\tsr$, we look at the results for $\tsr>2$.
Additionally, since we are interested in third-order temporal accuracy, we
focus on the high-order scheme with $k=3$ and $m=3$.

\begin{figure}[t!]
\begin{center}
$\begin{array}{ccc}
\includegraphics[height=45mm]{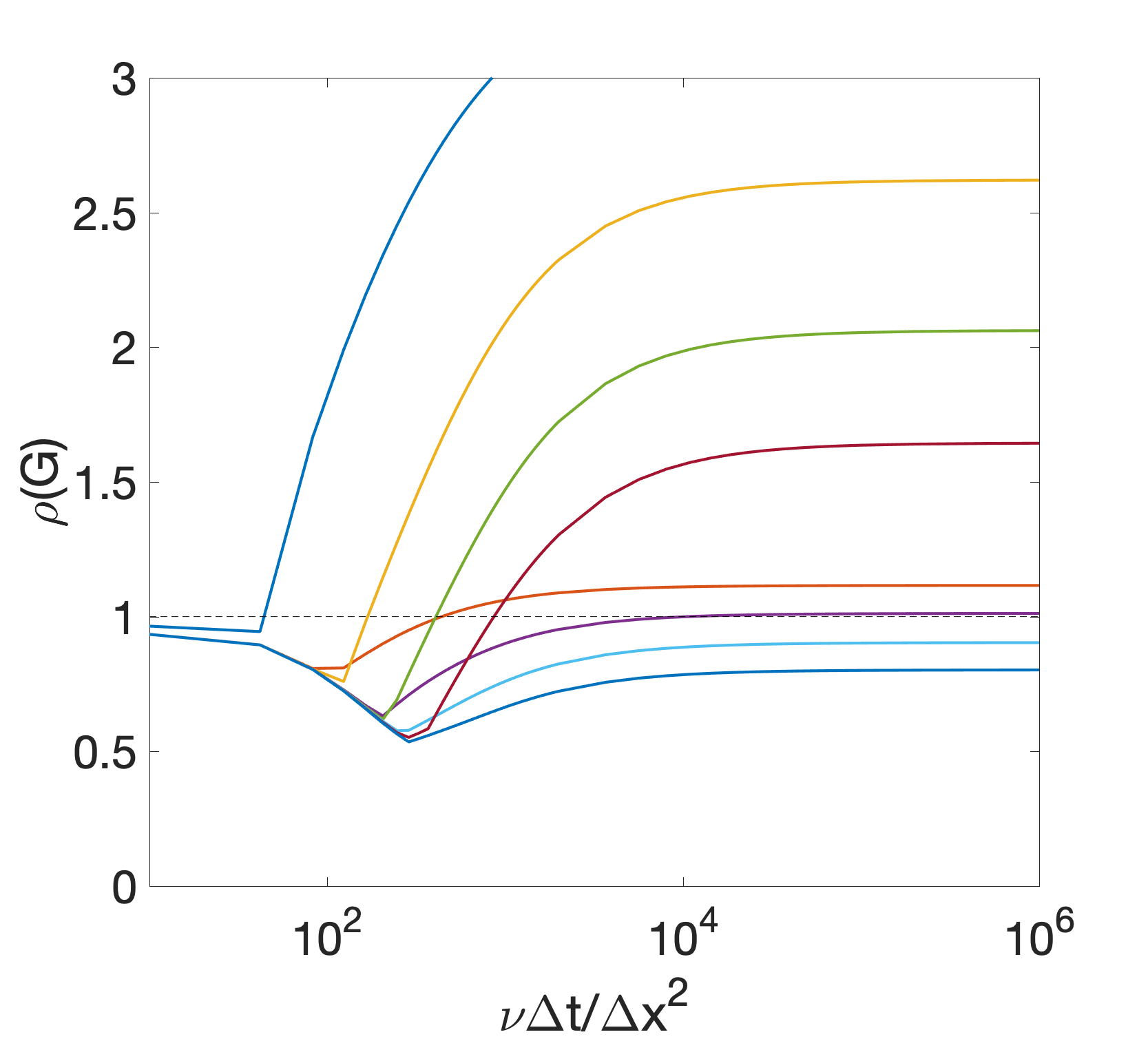} &
\includegraphics[height=45mm]{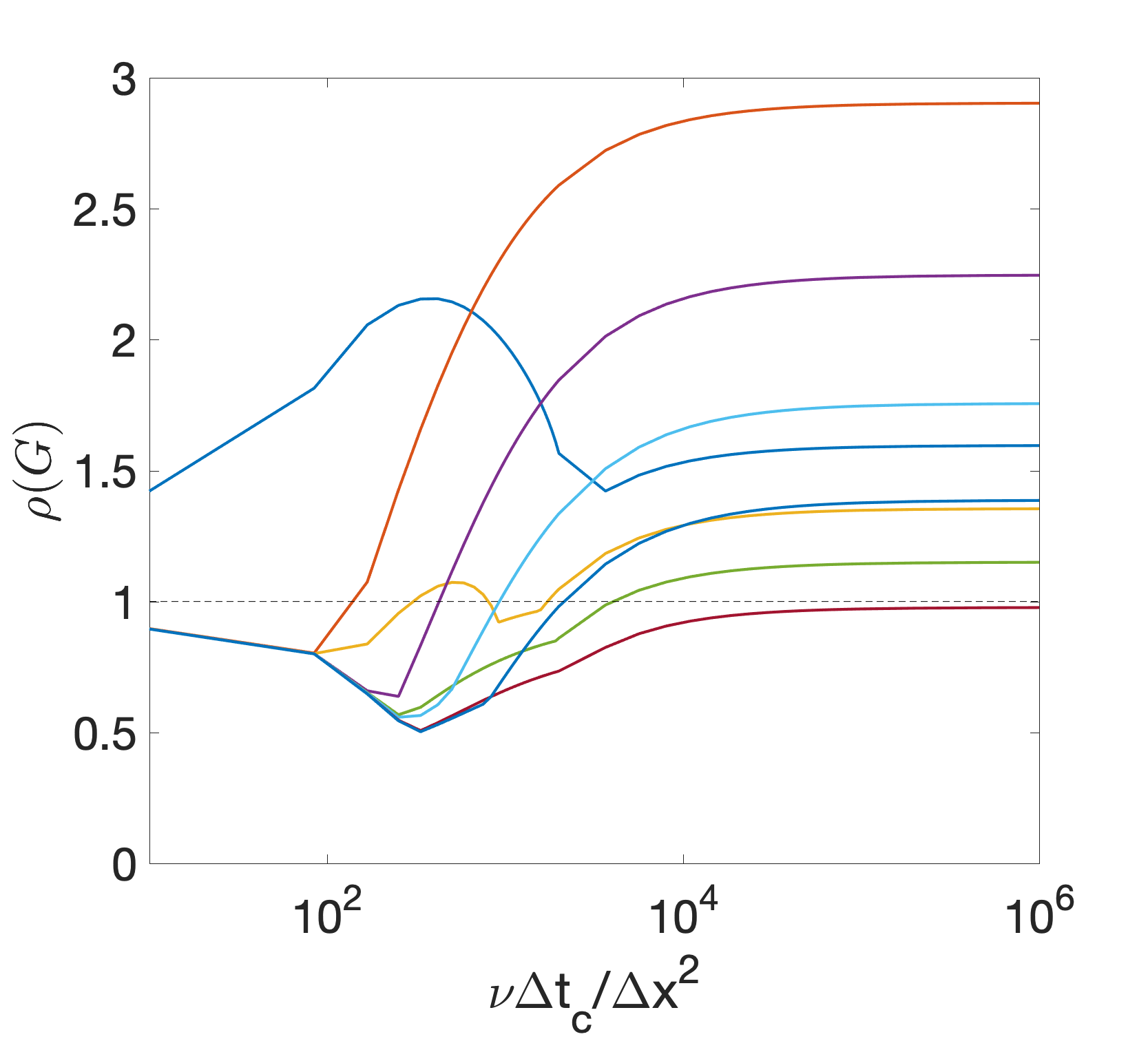} &
\includegraphics[height=45mm]{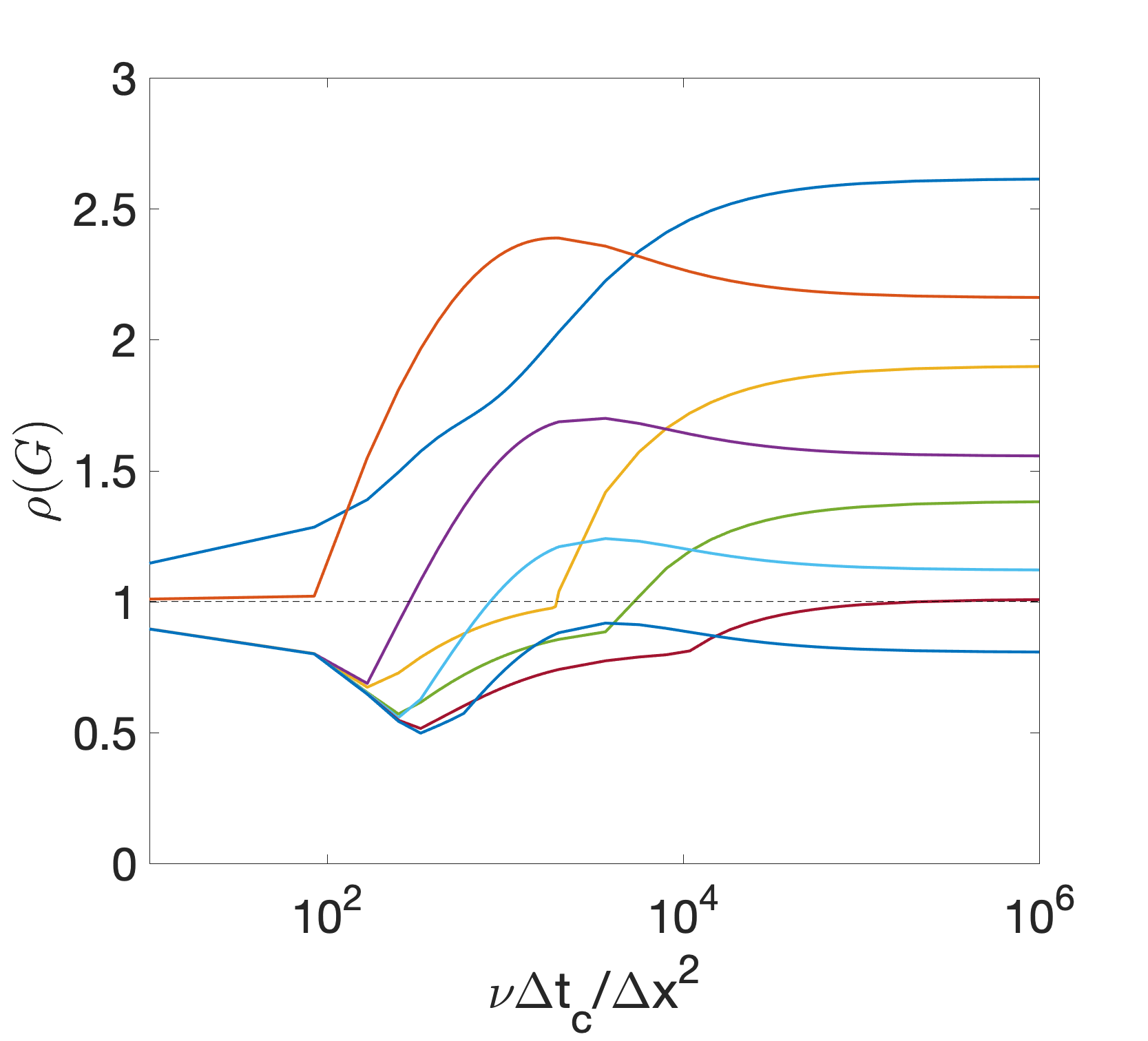} \\
\textrm{(a) $\tsr=1$} &
\textrm{(b) $\tsr=2$} &
\textrm{(c) $\tsr=3$} \\
\includegraphics[height=45mm]{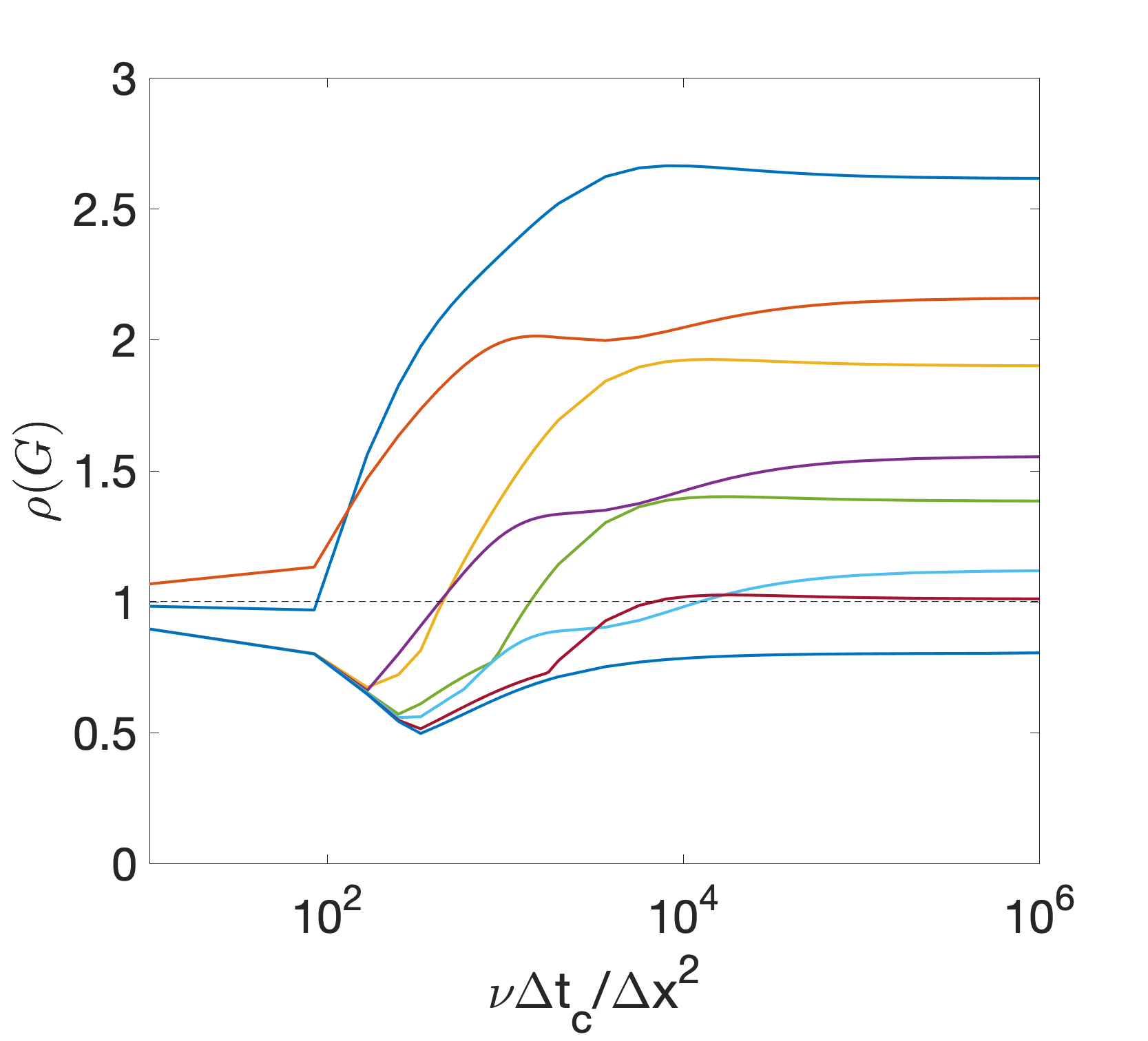} &
\includegraphics[height=45mm]{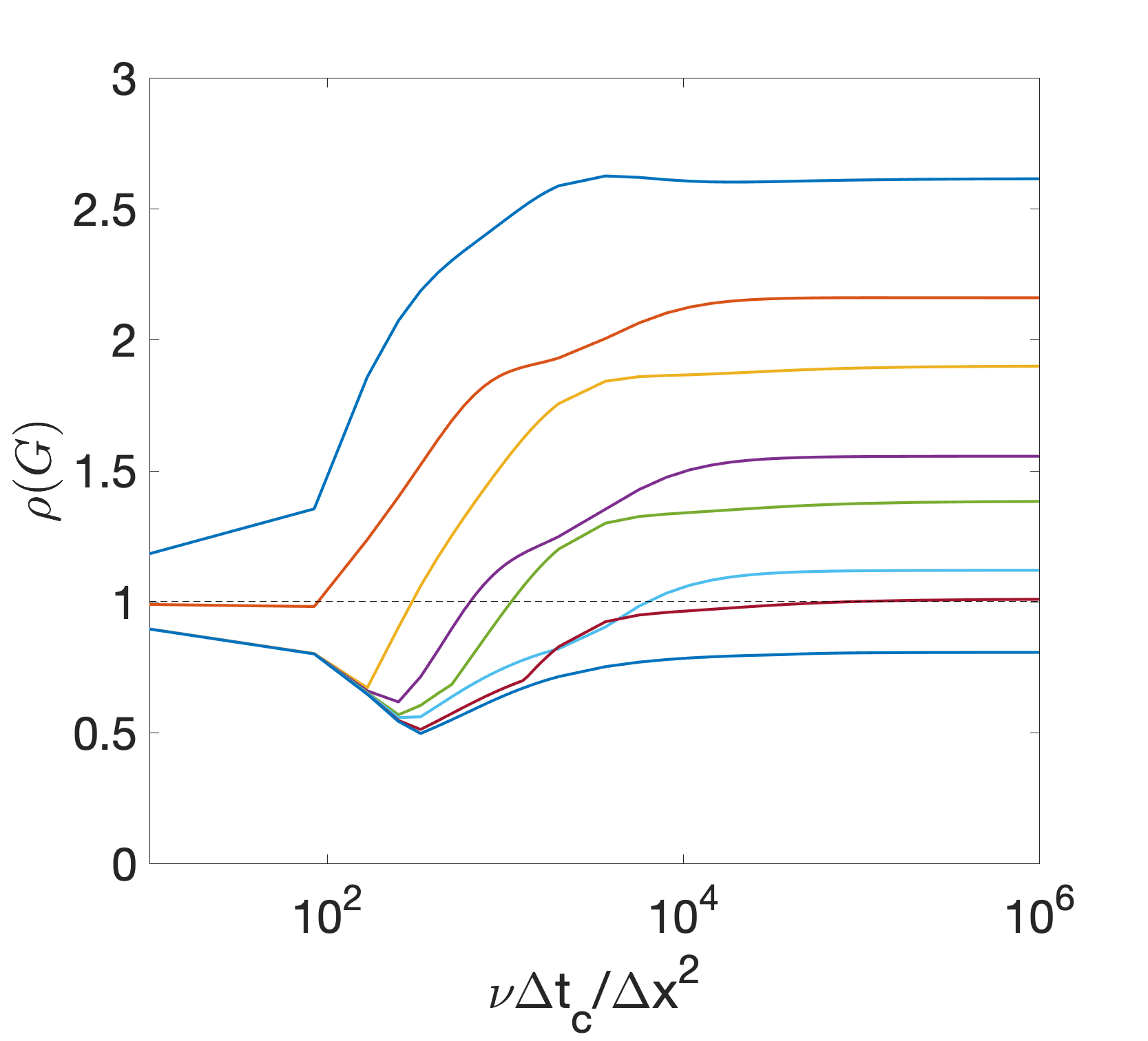} &
\includegraphics[height=45mm]{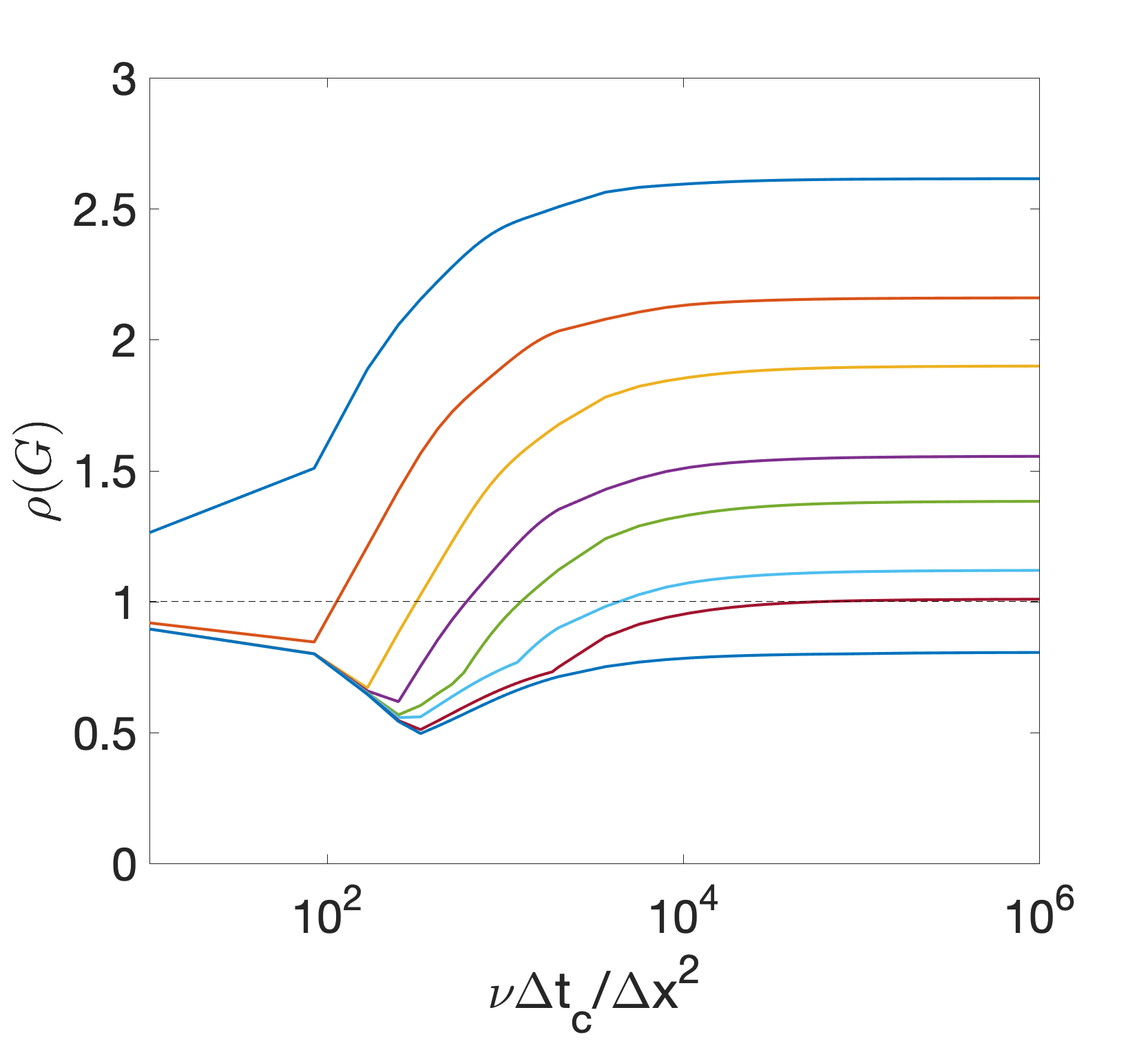} \\
\textrm{(d) $\tsr=4$} &
\textrm{(e) $\tsr=5$} &
\textrm{(f) $\tsr=10$} \\
\multicolumn{3}{c}{\includegraphics[width=150mm]{figures/legend}} \\
 \end{array}$
\end{center}
\vspace{-6mm}
\caption{Spectral radius $\rho(G)$ versus nondimensional time $\frac{\nu
\dt_c}{\dx^2}$ for the BDF$3$/EXT$3$ scheme with different $\tsr$ and $Q=0\dots7$, with grid
parameters $N=32$ and $K=5$.} \label{fig:rho_multi_int2_etamany2}
\end{figure}

Figure \ref{fig:rho_multi_int2_etamany2} shows the spectral radius versus
nondimensional time plot for $\tsr=1, 2, 3, 4, 5,$ and 10 with different $Q$.
The plots in Fig. \ref{fig:rho_multi_int2_etamany2} use a semi-log scale for
the x-axis to show the stability behavior of the multirate timestepping scheme
for a large range of nondimensional timestep size $\nu\dt_c/\dx^2$.  We observe
that for $\tsr=1$, odd-$Q$ is more stable than even-$Q$, and for $\tsr=2$,
even-$Q$ is more stable than odd-$Q$. For $\tsr\ge3$, however, the odd-even
pattern goes away for a large nondimensional timestep size
($\nu\dt_c/\dx^2>2\times10^4$).  We also observe that the singlerate
timestepping scheme ($\tsr=1$) requires fewer corrector iterations to guarantee
unconditional stability in comparison to the multirate timestepping scheme.
Here, unconditional stability means that $\rho(G)<1$ irrespective of the
nondimensional timestep size. Figure \ref{fig:rho_multi_int2_etamany2} shows
that for the grid parameters considered here ($N=32$ and $K=5$), $\tsr=1$
requires $Q=3$ and $\tsr\ge2$ requires $Q=6$ for unconditional stability to
solve the unsteady heat equation using overlapping grids with third-order
temporal accuracy in the FD-based framework.

\begin{figure}[t!]
\begin{center}
$\begin{array}{ccc}
\includegraphics[height=45mm]{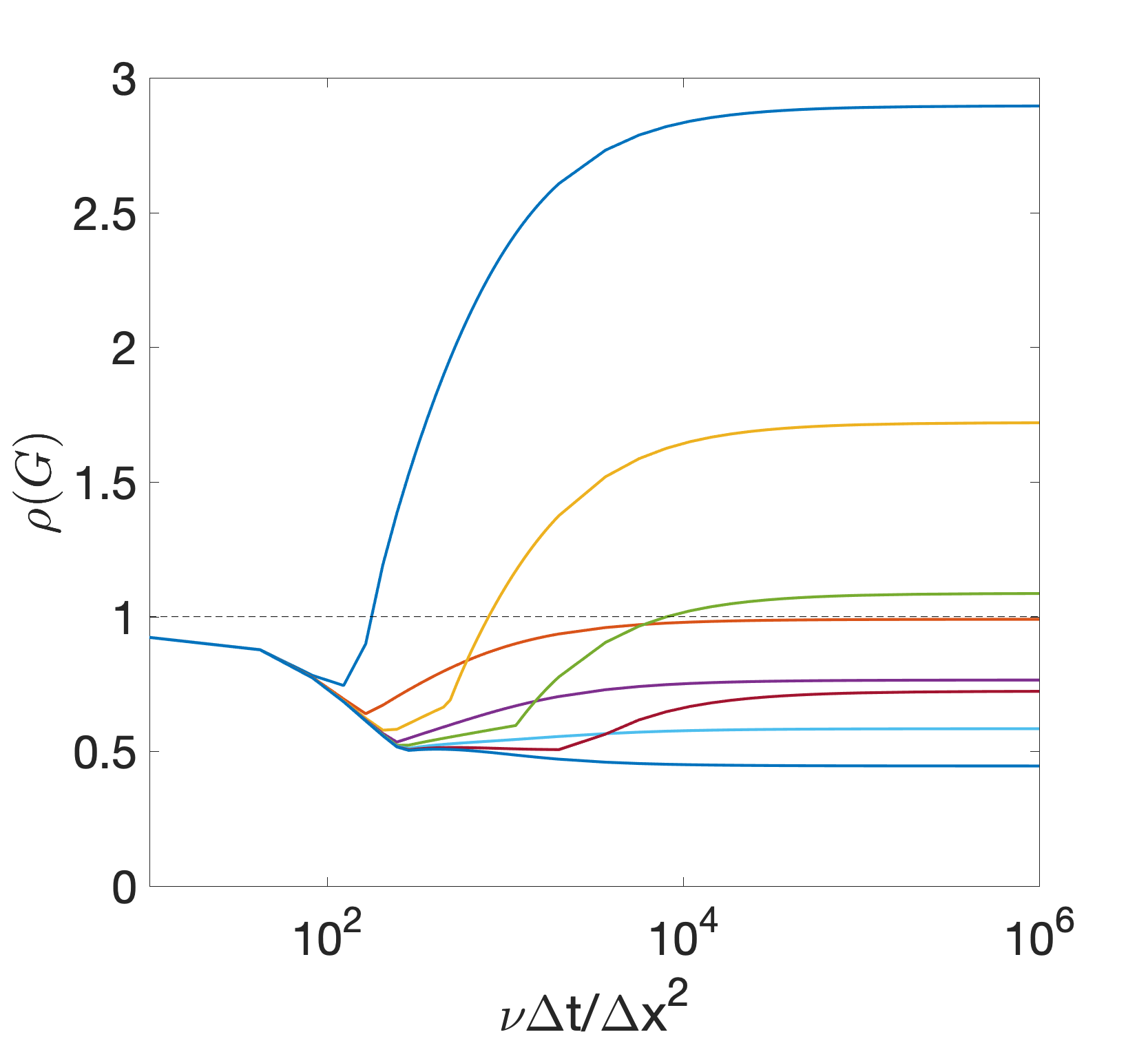} &
\includegraphics[height=45mm]{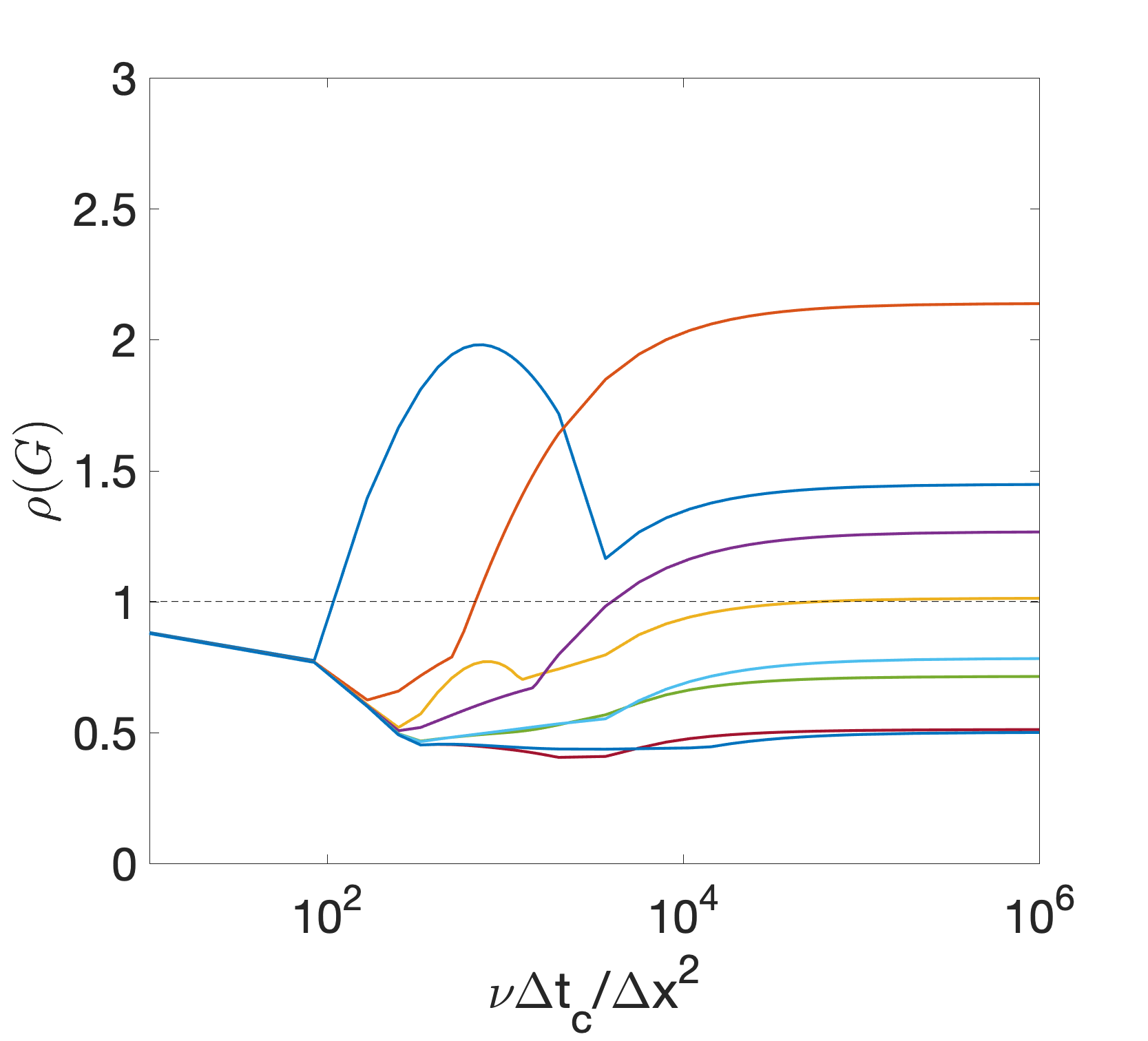} &
\includegraphics[height=45mm]{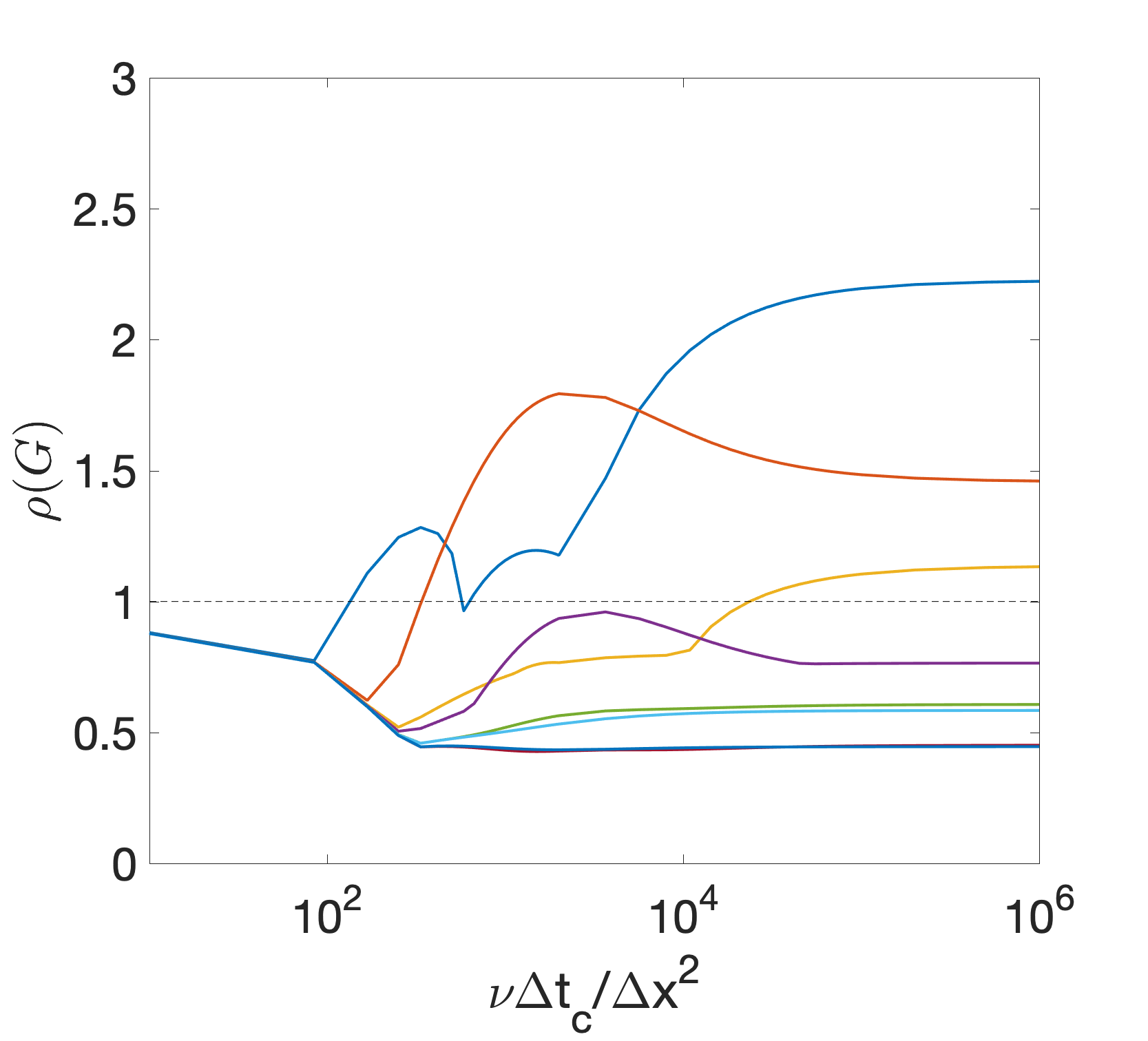} \\
\textrm{(a) $\tsr=1$} &
\textrm{(b) $\tsr=2$} &
\textrm{(c) $\tsr=3$} \\
\includegraphics[height=45mm]{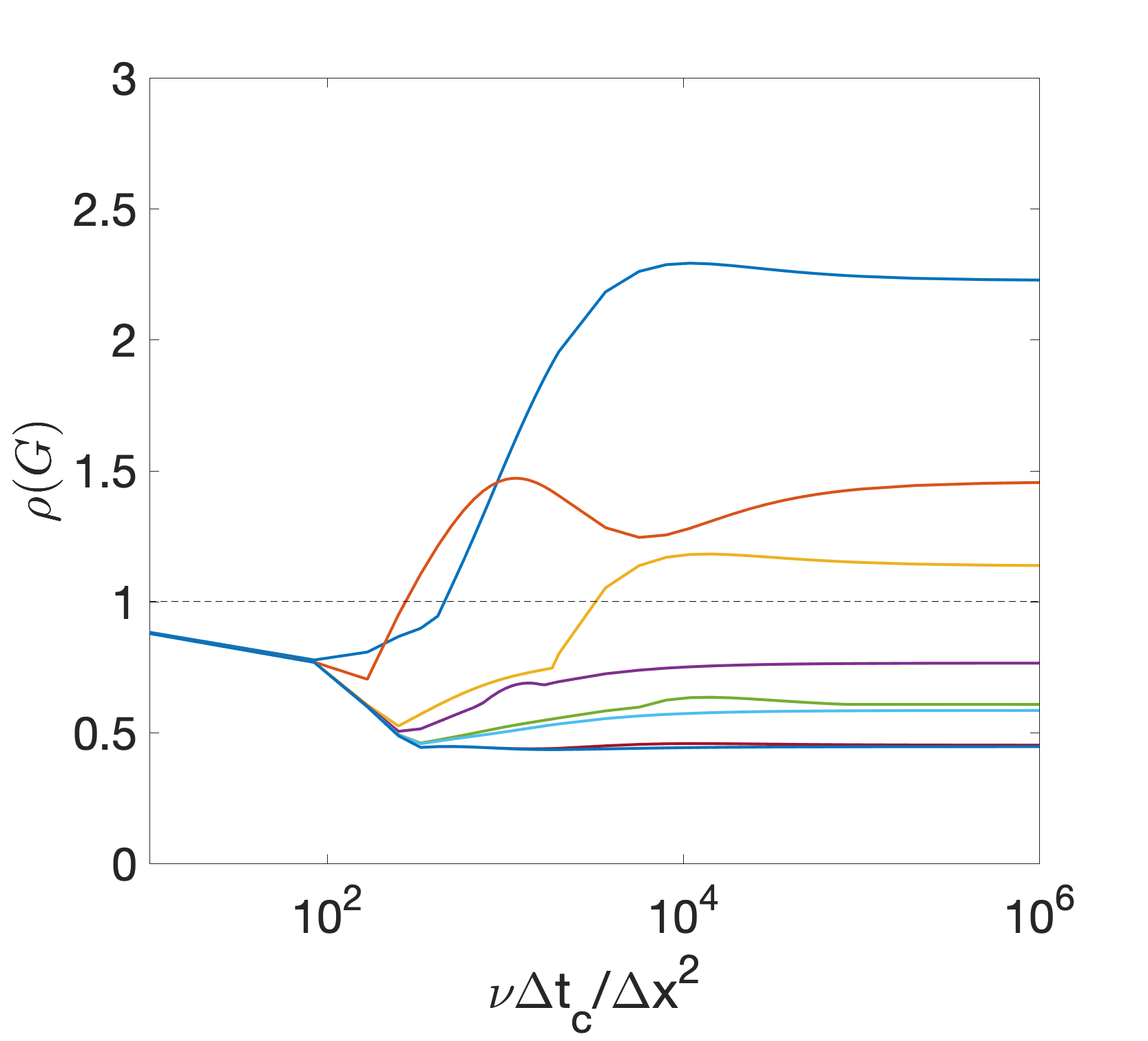} &
\includegraphics[height=45mm]{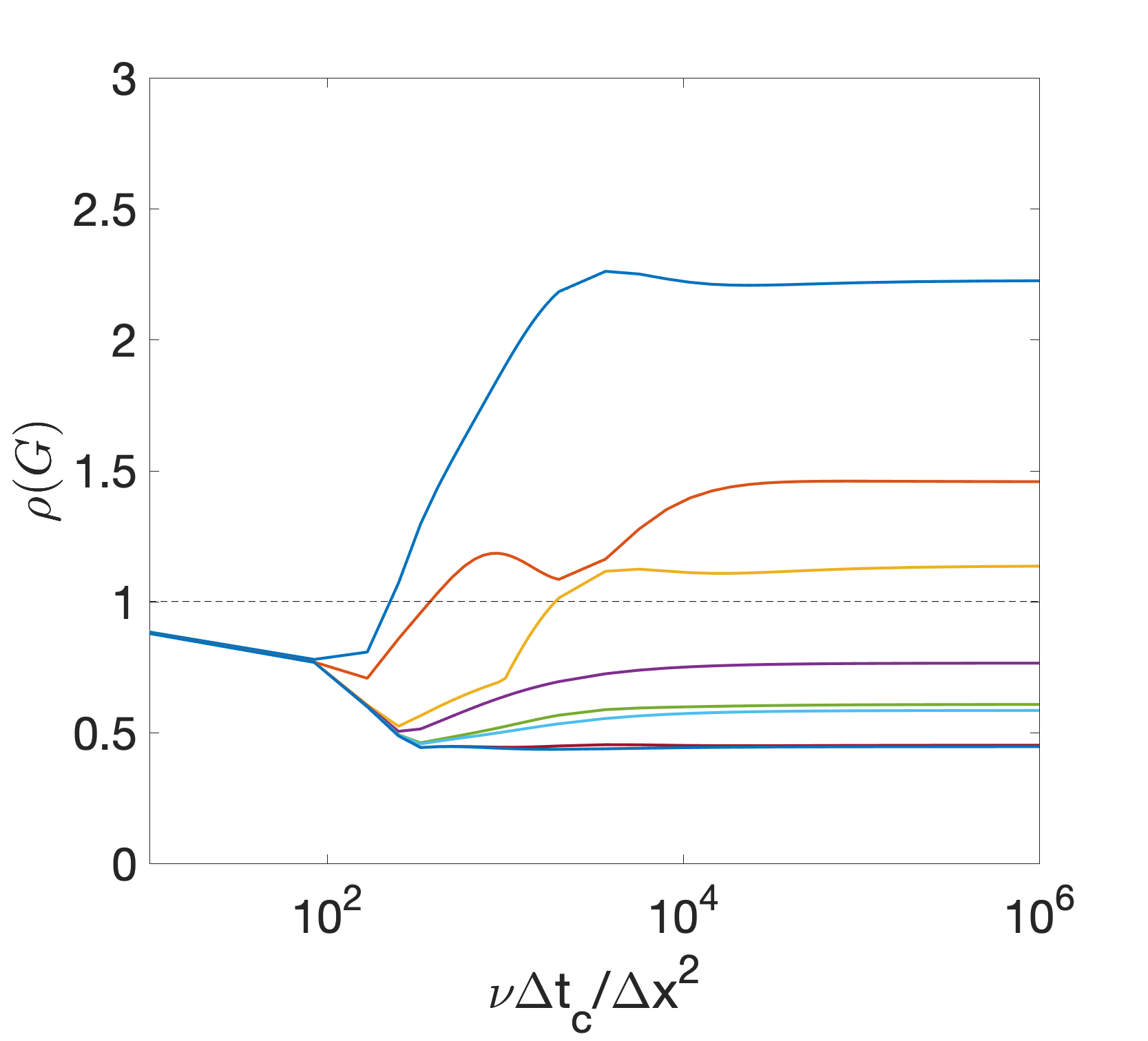} &
\includegraphics[height=45mm]{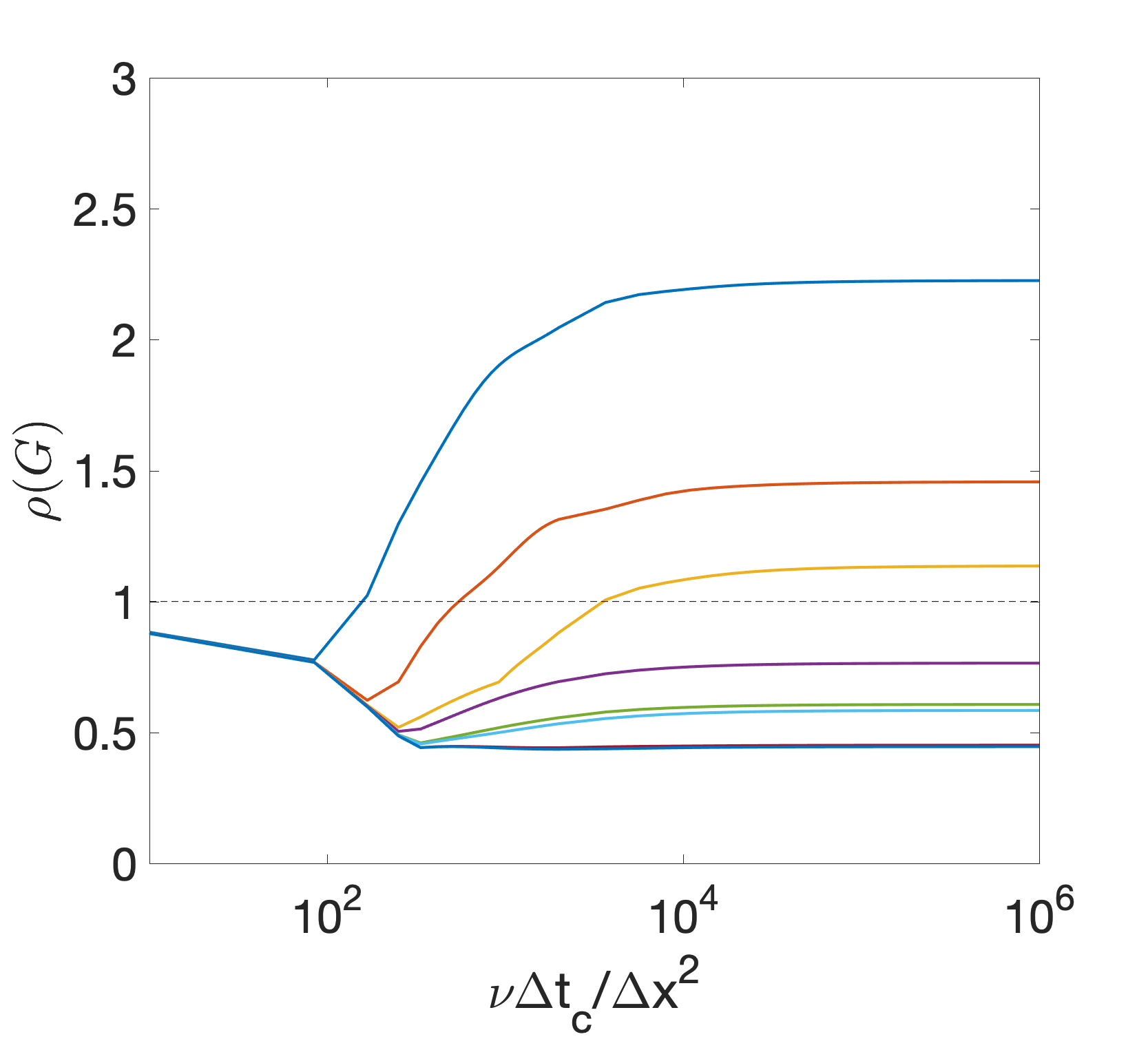} \\
\textrm{(d) $\tsr=4$} &
\textrm{(e) $\tsr=5$} &
\textrm{(f) $\tsr=10$} \\
\multicolumn{3}{c}{\includegraphics[width=150mm]{figures/legend}} \\
 \end{array}$
\end{center}
\vspace{-6mm}
\caption{Spectral radius $\rho(G)$ versus nondimensional time $\frac{\nu
\dt_c}{\dx^2}$ for the BDF$3$/EXT$3$ scheme with different $\tsr$ and $Q=0\dots7$, with grid
parameters $N=32$ and $K=10$.} \label{fig:rho_multi_int2_etamany3}
\end{figure}

We notice similar behavior in stability if we increase the grid overlap by
changing the grid parameter $K=5$ to $10$. Fig.
\ref{fig:rho_multi_int2_etamany3}(b), similar to Fig.
\ref{fig:rho_multi_int2_etamany2}(b), shows that even-$Q$ is more stable than
odd-$Q$ for $\tsr=2$. For $\tsr\ge3$, we observed that the odd-even pattern in
stability goes away for a large nondimensional timestep size, and as expected,
increasing the grid overlap makes the PC scheme more stable.

In the following section, we will show that the stability behavior that we have
observed for large nondimensional timestep size in the 1D model problem
qualitatively captures the general stability behavior of the PC-based multirate
timestepping scheme for solving the INSE in the Schwarz-SEM framework. In
future work, we will extend this analysis to make more rigorous predictions and
establish theoretical bounds on the stability of the PC-based multirate
timestepping scheme. This will require us to understand how the nondimensional
timestep size of the 1D model problem is related to the timestep size for the
unsteady Stokes problem in the Schwarz-SEM framework. We will also investigate
why the odd-even stability pattern manifests for $\tsr=1$ and 2, and not for
$\tsr\ge3$.

\subsection{Validation with Schwarz-SEM Framework}
\begin{figure}[t!]
\begin{center}
$\begin{array}{c}
\includegraphics[height=45mm]{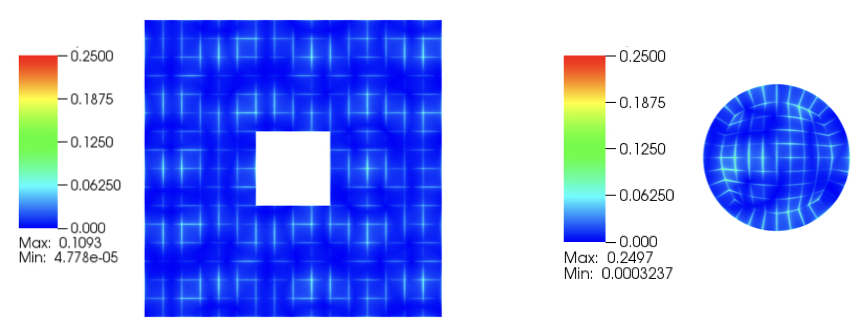}
\end{array}$
\end{center}
\caption{CFL comparison for the overlapping spectral element mesh for same
timestep size. (left) $\Omega^c$ and (right) $\Omega^f$. The CFL is maximum
($0.11$ in $\Omega^c$ and $0.25$ in $\Omega^f$) where
the ratio of flow velocity to grid spacing is highest in each subdomain.}
\label{fig:nnmultieddycfl}
\end{figure}

In this section, we consider the Navier-Stokes eigenfunctions by Walsh from
Section \ref{sec:oddevenNSE} using multirate timestepping in the Schwarz-SEM
framework.  Figure \ref{fig:nnmultieddycfl} shows a  snapshot of the local grid
CFL for the overlapping grids (Fig. \ref{fig:nngrid}) used to model the
periodic domain. Due to the difference in the grid sizes, the maximum CFL
number in the circular mesh (CFL=0.2497) is twice as much as the maximum CFL in
the background mesh (CFL=0.1093). Thus, we can use at least a timestep ratio of
$\eta = \dt_c/\dt_f = 2$ with the larger timestep size ($\dt_c$) for the
background mesh and the smaller timestep size ($\dt_f$) for the circular mesh.

\begin{figure}[t!]
\begin{center}
$\begin{array}{ccc}
\includegraphics[height=40mm]{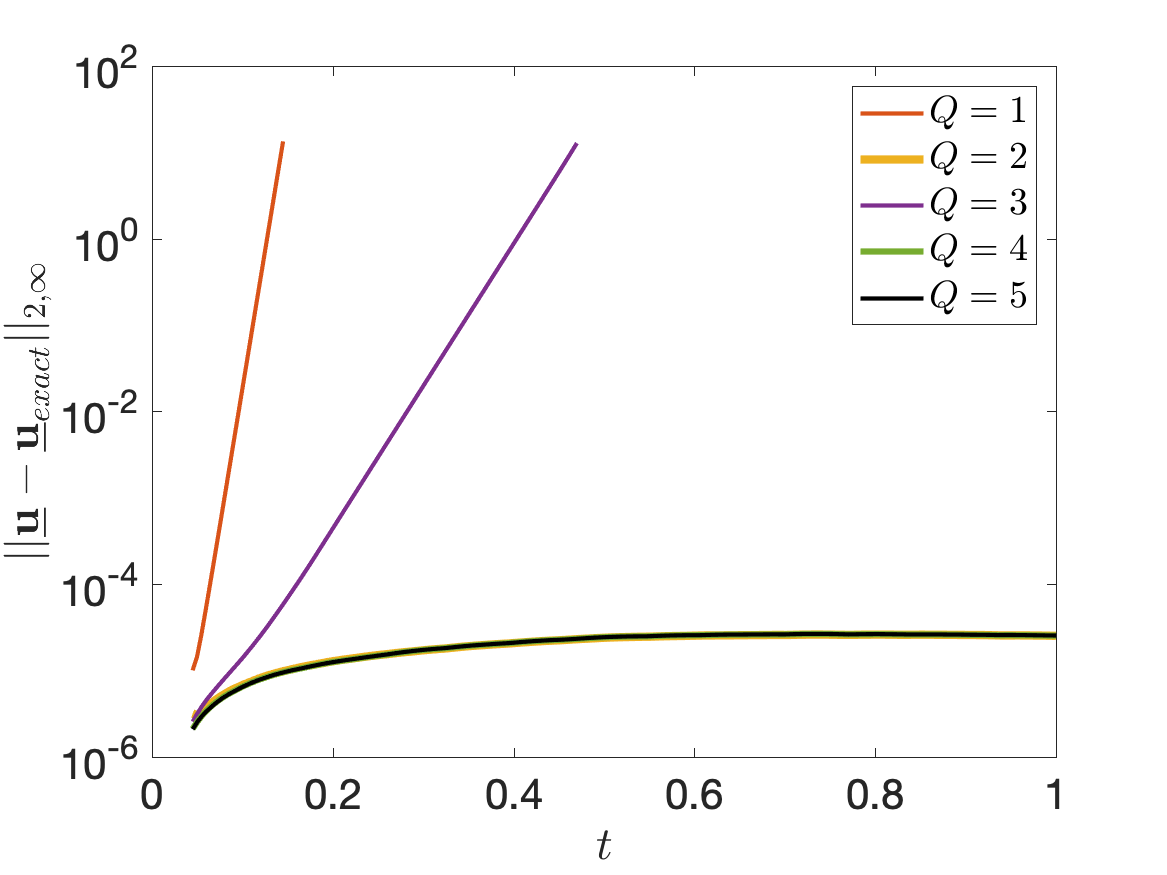} &
\includegraphics[height=40mm]{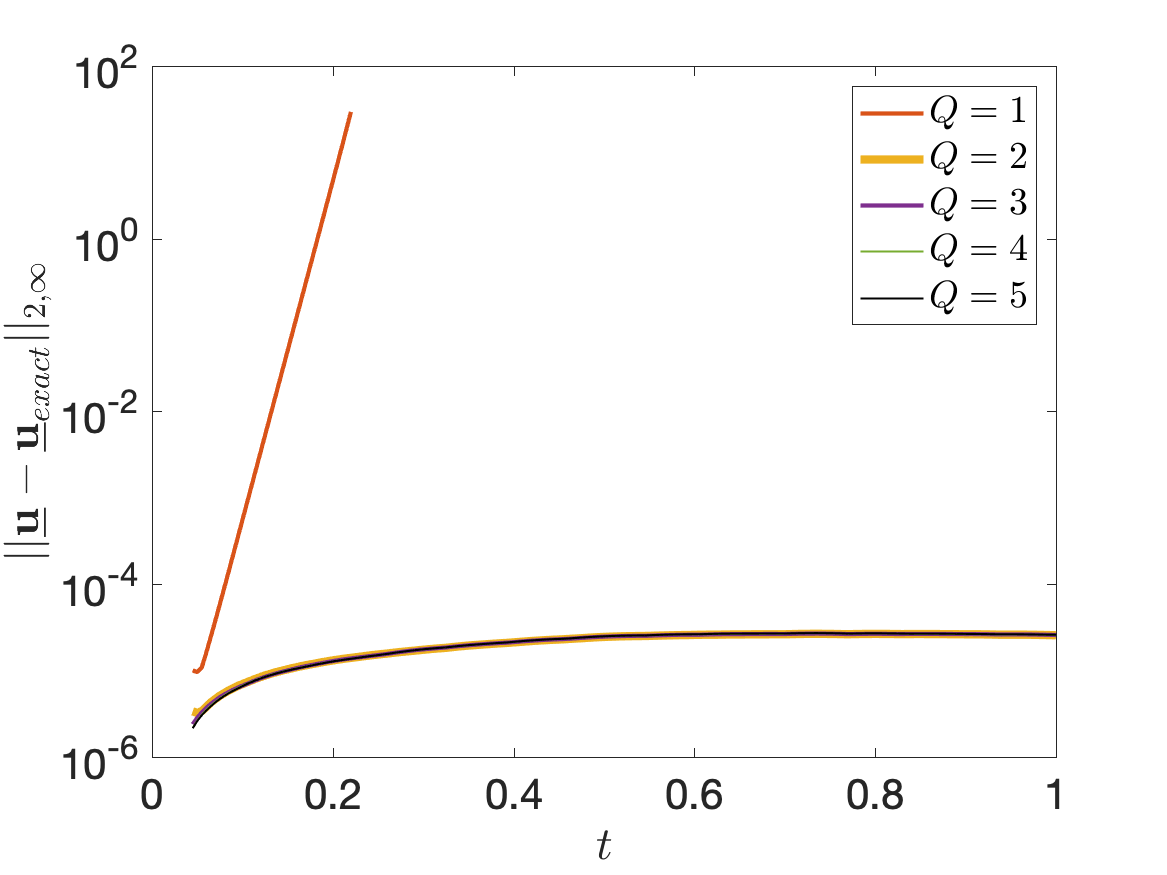} &
\includegraphics[height=40mm]{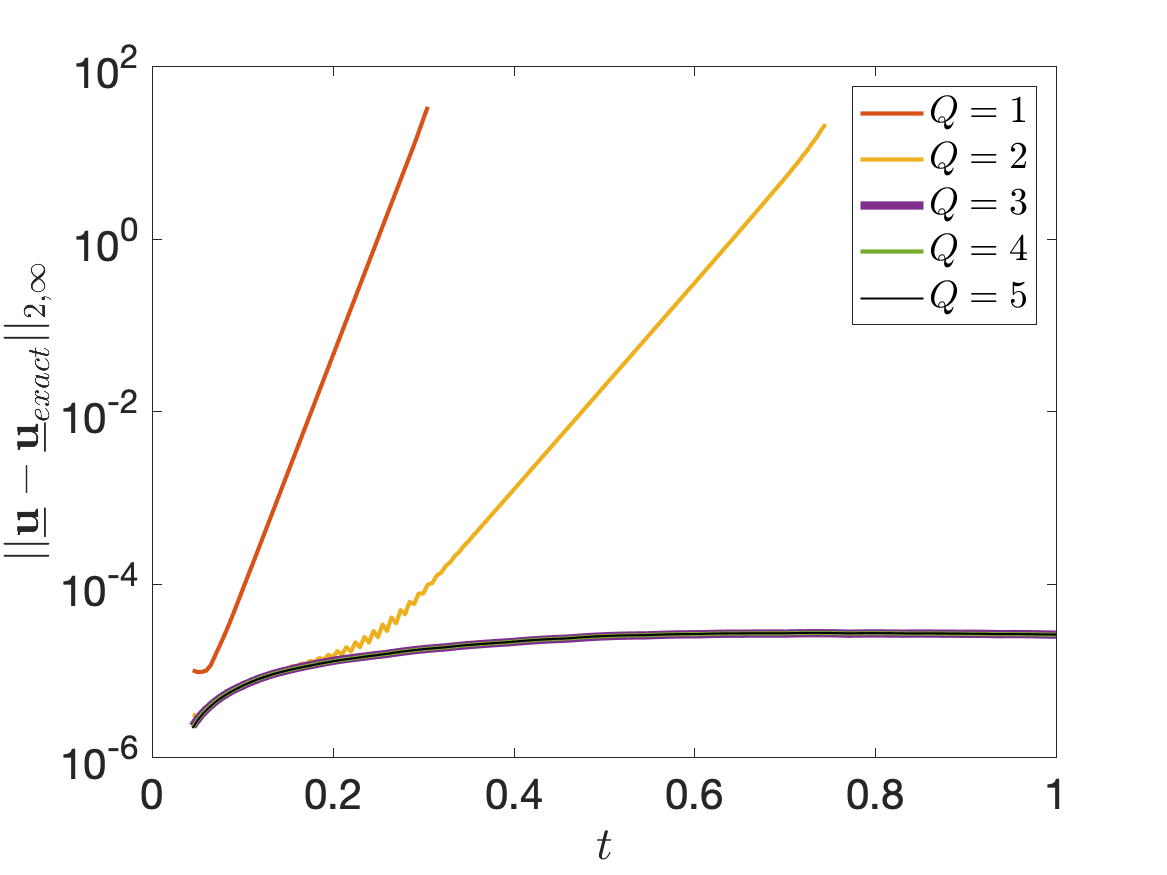} \\
\textrm{(a) } \tsr =2 &
\textrm{(b) } \tsr =3 &
\textrm{(c) } \tsr =4
\end{array}$
\end{center}
\caption{Error variation for the Navier-Stokes eigenfunctions test case from
Section \ref{sec:oddevenNSE} with different $Q$ for $\tsr=2-4$. We set $k=3$, $m=3$, $\dt_c=5\times
10^{-3}$, and $N=7$ for the results presented here.}
\label{fig:eddyoddeveneta}
\end{figure}

To understand the stability properties of the PC-based multirate timestepping
scheme in the Schwarz-SEM framework, we set $\dt_c=5\times 10^{-3}$, $N=7$,
$k=3$, and $m=3$ with different timestep ratio ($\tsr=2-4$) and corrector
iterations ($Q=1-5$). Preliminary results using the Schwarz-SEM framework show
that we observe the stability behavior that we had expected from the analysis
of the 1D model problem.  Figure \ref{fig:eddyoddeveneta}(a) shows that for
$\tsr=2$, odd-$Q$ is less stable than even-$Q$, as we had expected from the
results in Fig.  \ref{fig:rho_multi_int2_etamany2}(b) and
\ref{fig:rho_multi_int2_etamany3}(b).  Figure \ref{fig:eddyoddeveneta}(b) and
(c) shows that for $\tsr=3$ and 4, respectively, we do not observe the odd-even
stability pattern in the Schwarz-SEM framework, which we had observed for large
nondimensional timestep size in the FD-based framework (Fig.
\ref{fig:rho_multi_int2_etamany2}(c) and (d)).
We note that we have observed this same behavior for similar numerical experiments
in the Schwarz-SEM framework for different $N$.

Based on numerical experiments in the Schwarz-SEM framework
with the singlerate and multirate timestepping PC scheme, we conclude that the
asymptotic behavior (in terms of the nondimensional timestep size) that we
observe for different $Q$ and $\tsr$ in the 1D model problem, qualitatively
captures the stability behavior that we observe in the Schwarz-SEM framework.
We note that from a practical standpoint, the results in the current work are
sufficient because the odd-even behavior seems to vanish for $\tsr >= 10$,
which is the limit in which we are interested as high timestep ratios help us
realize the maximum potential of MTS. For example, we have used $\tsr=50$ to
model a thermally buoyant plume with two overlapping grids in
\cite{mittal2020multirate} and demonstrated the computational savings
associated with MTS in comparison to STS. Nonetheless, in future work we will
look at methods that can allow us to make more rigorous predictions on the
impact of $Q$ and $\tsr$ on the stability properties of the multirate
timestepping scheme.

\section{Conclusion \& future work}
We have used the matrix stability method to analyze the
stability of a singlerate and multirate predictor-corrector scheme used in the
Schwarz-SEM framework for the incompressible Navier-Stokes equations.
We simplify the analysis by considering the unsteady heat equation in 1D with a finite-difference-based
spatial discretization, and results indicate that the stability
of our BDF$k$/EXT$m$-based timestepping scheme increases with increase in
resolution and overlap of the subdomains. We also observe that for singlerate timestepping,
the PC scheme is more stable when odd number of corrector iteration ($Q$) are used in
comparison to an even-$Q$. Based on empirical analyses and the results
in the literature, it appears that this difference in the stability
of odd- and even-$Q$ is a universal behavior of
PC schemes for ODEs of the form $y' = \lambda y$ if $\lambda$ has a negative
real part. In future work, we will further investigate this stability
behavior and look at the stability of our PC scheme with the unsteady Stokes
operator in the context of SEM. For multirate timestepping,
we have observed that the difference in stability
of odd- and even-$Q$ varies with the timestep ratio$\eta$.
For timestep ratio $\eta=2$, even-$Q$ is more stable than odd-$Q$. For $\eta\ge3$, even-$Q$ is
more stable than odd-$Q$ for a small nondimensional timestep size, and the
odd-even behavior vanishes as we increase the timestep size. In future work,
we will further explore the stability behavior of
the multirate PC scheme in the context of the Schwarz-SEM framework. We will also
consider the multirate scheme for a system of ODEs of the form $y' = \lambda y$ to determine
if the stability behavior can be generalized for ODEs and PDEs in general, similar
to the results for the singlerate timestepping scheme.

\bibliographystyle{elsarticle-num}
\bibliography{lit}

\begin{thebibliography}{10}
\expandafter\ifx\csname url\endcsname\relax
  \def\url#1{\texttt{#1}}\fi
\expandafter\ifx\csname urlprefix\endcsname\relax\def\urlprefix{URL }\fi
\expandafter\ifx\csname href\endcsname\relax
  \def\href#1#2{#2} \def\path#1{#1}\fi

\bibitem{mittal2019nonconforming}
K.~Mittal, S.~Dutta, P.~Fischer, Nonconforming {S}chwarz-spectral element
  methods for incompressible flow, Computers \& Fluids (2019) 104237.

\bibitem{henshaw1994}
W.~D. Henshaw, A fourth-order accurate method for the incompressible
  {N}avier-{S}tokes equations on overlapping grids, Journal of computational
  physics 113~(1) (1994) 13--25.

\bibitem{merrill2016}
B.~E. Merrill, Y.~T. Peet, P.~F. Fischer, J.~W. Lottes, A spectrally accurate
  method for overlapping grid solution of incompressible {N}avier--{S}tokes
  equations, Journal of Computational Physics 307 (2016) 60--93.

\bibitem{rogers1991steady}
S.~E. Rogers, D.~Kwak, C.~Kiris, Steady and unsteady solutions of the
  incompressible {N}avier-{S}tokes equations, AIAA journal 29~(4) (1991)
  603--610.

\bibitem{cd2012v7}
S.-C. CD-adapco, V7. 02.008, User Manual.

\bibitem{ahmad1996helicopter}
J.~Ahmad, E.~P. Duque, Helicopter rotor blade computation in unsteady flows
  using moving overset grids, Journal of Aircraft 33~(1) (1996) 54--60.

\bibitem{cambier2013onera}
L.~Cambier, S.~Heib, S.~Plot, The {Onera elsA CFD} software: input from
  research and feedback from industry, Mechanics \& Industry 14~(3) (2013)
  159--174.

\bibitem{eberhardt1985overset}
S.~Eberhardt, D.~Baganoff, Overset grids in compressible flow, in: 7th
  Computational Physics Conference, 1985, p. 1524.

\bibitem{nicholsoverflowmanual}
R.~H. Nichols, P.~G. Buning, User’s manual for {OVERFLOW} 2.1, University of
  Alabama and NASA Langley Research Center.

\bibitem{saunier2008third}
O.~Saunier, C.~Benoit, G.~Jeanfaivre, A.~Lerat, Third-order {C}artesian overset
  mesh adaptation method for solving steady compressible flows, International
  journal for numerical methods in fluids 57~(7) (2008) 811--838.

\bibitem{blake1996overset}
D.~Blake, T.~Buter, Overset grid methods applied to a finite-volume time-domain
  {M}axwell equation solver, in: 27th Plasma Dynamics and Lasers Conference,
  1996, p. 2338.

\bibitem{angel2018}
J.~B. Angel, J.~W. Banks, W.~D. Henshaw, A high-order accurate {FDTD} scheme
  for {M}axwell's equations on overset grids, in: Applied Computational
  Electromagnetics Society Symposium (ACES), 2018 International, IEEE, 2018,
  pp. 1--2.

\bibitem{meng2017stable}
F.~Meng, J.~W. Banks, W.~D. Henshaw, D.~W. Schwendeman, A stable and accurate
  partitioned algorithm for conjugate heat transfer, Journal of Computational
  Physics 344 (2017) 51--85.

\bibitem{kao1997application}
K.-H. Kao, M.-S. Liou, Application of chimera/unstructured hybrid grids for
  conjugate heat transfer, AIAA journal 35~(9) (1997) 1472--1478.

\bibitem{henshaw2009composite}
W.~D. Henshaw, K.~K. Chand, A composite grid solver for conjugate heat transfer
  in fluid-structure systems, Journal of Computational Physics 228~(10) (2009)
  3708--3741.

\bibitem{koblitz2017}
A.~Koblitz, S.~Lovett, N.~Nikiforakis, W.~Henshaw, Direct numerical simulation
  of particulate flows with an overset grid method, Journal of Computational
  Physics 343 (2017) 414--431.

\bibitem{mittal2020multirate}
K.~Mittal, S.~Dutta, P.~Fischer, Multirate timestepping for the incompressible
  navier-stokes equations in overlapping grids, arXiv preprint
  arXiv:2003.00347.

\bibitem{mittal2020direct}
K.~Mittal, S.~Dutta, P.~Fischer, Direct numerical simulation of rotating
  ellipsoidal particles using moving nonconforming schwarz-spectral element
  method, Computers \& Fluids (2020) 104556.

\bibitem{chatterjee2019towards}
T.~Chatterjee, S.~S. Patel, M.~M. Ameen, Towards improved mesh-designing
  techniques of spark-ignition engines in the framework of spectral element
  methods, Tech. rep., Argonne National Lab.(ANL), Argonne, IL (United States)
  (2019).

\bibitem{patera84}
A.~T. Patera, A spectral element method for fluid dynamics: {laminar} flow in a
  channel expansion, Journal of computational Physics 54~(3) (1984) 468--488.

\bibitem{dfm02}
M.~O. Deville, P.~F. Fischer, E.~H. Mund, High-order methods for incompressible
  fluid flow, Vol.~9, Cambridge University Press, 2002.

\bibitem{varga1999matrix}
R.~S. Varga, Matrix iterative analysis, Vol.~27, Springer Science \& Business
  Media, 1999.

\bibitem{hamming1959stable}
R.~W. Hamming, Stable predictor-corrector methods for ordinary differential
  equations, Journal of the ACM (JACM) 6~(1) (1959) 37--47.

\bibitem{chase1962stability}
P.~Chase, Stability properties of predictor-corrector methods for ordinary
  differential equations, Journal of the ACM (JACM) 9~(4) (1962) 457--468.

\bibitem{hall1967stability}
G.~Hall, The stability of predictor-corrector methods, The Computer Journal
  9~(4) (1967) 410--412.

\bibitem{peet2012}
Y.~T. Peet, P.~F. Fischer, Stability analysis of interface temporal
  discretization in grid overlapping methods, SIAM Journal on Numerical
  Analysis 50~(6) (2012) 3375--3401.

\bibitem{mathew2003maximum}
T.~Mathew, G.~Russo, Maximum norm stability of difference schemes for parabolic
  equations on overset nonmatching space-time grids, Mathematics of Computation
  72~(242) (2003) 619--656.

\bibitem{wu2012convergence}
S.-L. Wu, C.-M. Huang, T.-Z. Huang, Convergence analysis of the overlapping
  schwarz waveform relaxation algorithm for reaction-diffusion equations with
  time delay, IMA Journal of Numerical Analysis 32~(2) (2012) 632--671.

\bibitem{love2009stability}
E.~Love, W.~J. Rider, G.~Scovazzi, Stability analysis of a
  predictor/multi-corrector method for staggered-grid lagrangian shock
  hydrodynamics, Journal of Computational Physics 228~(20) (2009) 7543--7564.

\bibitem{stetter1968improved}
H.~J. Stetter, Improved absolute stability of predictor-corrector schemes,
  Computing 3~(4) (1968) 286--296.

\bibitem{schwarz1870}
H.~A. {S}chwarz, Ueber einen Grenz{\"u}bergang durch alternirendes Verfahren,
  Z{\"u}rcher u. Furrer, 1870.

\bibitem{gslibrepo}
P.~Fishcer, \href{https://github.com/gslib/gslib}{gslib/gslib} (2017).
\newline\urlprefix\url{https://github.com/gslib/gslib}

\bibitem{fornberg1998practical}
B.~Fornberg, A practical guide to pseudospectral methods, Vol.~1, Cambridge
  university press, 1998.

\bibitem{smith2004domain}
B.~Smith, P.~Bjorstad, W.~Gropp, Domain decomposition: parallel multilevel
  methods for elliptic partial differential equations, Cambridge university
  press, 2004.

\bibitem{mittaloverlapping}
K.~Mittal, S.~Dutta, P.~F. Fischer, Overlapping schwarz based spectral element
  method for incompressible flow in complex domains, Letter from Lead-Chairs.

\bibitem{walsh1992eddy}
O.~Walsh, Eddy solutions of the {N}avier--{S}tokes equations, in: The
  {N}avier-{S}tokes Equations II—Theory and Numerical Methods, Springer,
  1992, pp. 306--309.

\end{thebibliography}
\end{document}